\documentclass[mnsc,nonblindrev]{informs3_arXiv} 

\DoubleSpacedXI 

\usepackage{endnotes}
\let\footnote=\endnote
\let\enotesize=\normalsize
\def\notesname{Endnotes}%
\def\makeenmark{$^{\theenmark}$}
\def\enoteformat{\rightskip0pt\leftskip0pt\parindent=1.75em
  \leavevmode\llap{\theenmark.\enskip}}

\usepackage{natbib}
 \bibpunct[, ]{(}{)}{,}{a}{}{,}%
 \def\bibfont{\small}%
\TheoremsNumberedThrough     
\ECRepeatTheorems

\EquationsNumberedThrough    

\def\ie{{\it i.e.,\ \/}}
\def\eg{{\it e.g.,\ \/}}

\def\defeq{{\,\stackrel{\Delta}{=}}\,}

\def\nn{{\nonumber}}
\newcommand{\Tau}{\mathcal{T}}

\makeatletter
\newcommand\notsotiny{\@setfontsize\notsotiny\@viipt\@viipt}
\makeatother

\usepackage{bbm}
\usepackage{extarrows}
\usepackage{algorithmicx,algorithm}
\usepackage[noend]{algpseudocode}
\usepackage{diagbox}
\usepackage[hidelinks]{hyperref}
\usepackage[T1]{fontenc}
\usepackage{cleveref}

\begin{document}


\RUNAUTHOR{Liu}

\RUNTITLE{Index Policy for Partially Observable RMAB}

\TITLE{Relaxed Indexability and Index Policy for Partially Observable Restless Bandits}

\ARTICLEAUTHORS{%
\AUTHOR{Keqin Liu}
\AFF{Department of Financial and Actuarial Mathematics, Xi'an Jiaotong-Liverpool University, Suzhou, China, 215123; Jiangsu National Center for Applied Mathematics, Nanjing, China, 210093. \EMAIL{Keqin.Liu@xjtlu.edu.cn}} 
} 

\ABSTRACT{%
This paper addresses an important class of restless multi-armed bandit (RMAB) problems that finds broad application in operations research,
stochastic optimization, and reinforcement learning. There are $N$ independent Markov processes that may be operated, observed and offer
rewards. Due to the resource constraint, we can only choose a subset of $M~(M<N)$ processes to operate and accrue reward determined
by the states of selected processes. We formulate the problem as a partially observable RMAB with an infinite state space and design an
algorithm that achieves a near-optimal performance with low complexity. Our algorithm is based on a generalization of Whittle's original idea of indexability. Referred to as the relaxed indexability, the extended definition leads to the efficient online verifications and computations of the approximate Whittle index under the proposed algorithmic framework.
}%



\maketitle

%


\section{Introduction}

The first multi-armed bandit (MAB) problem was proposed in 1933 in the context of a clinical trial for adaptively selecting
the best treatment over time \citep{T1933}. Specifically, given two new medicines just invented for curing some disease, we want to find out which medicine has the better effect in long-run. When a new patient arrives, the doctor needs to decide which medicine to use based on past observations on the recovery processes of the previous patients after being treated with one of the medicines. Because the effect of each medicine is often modeled as a {\em random} variable, the decision problem involves a famous dilemma between ``Exploitation'' and ``Exploration'' often appeared in reinforcement learning: choosing the medicine that seems to be the best versus choosing the one less frequently used. In other words, the choice of the medicine determines not only the immediate effect of the treatment but also which medicine to observe for better estimation in the future. In the following subsection, we formally state the classical MAB under the Bayesian framework.

\subsection{The Classical MAB and Gittins Index}
In the classical Bayesian model of MAB, there are~$N$ arms and a single player. At each discrete time (decision epoch), a player chooses one arm to operate and accrues certain amount of reward determined by the state of the arm. The state of the chosen arm transits to a new one according to
a {\em known} Markovian rule while the states of other arms remain frozen. The observation model is assumed to be complete, \ie the states of all arms can be observed before deciding which
arm to choose. The objective is to maximize the total discounted reward over the infinite horizon \citep{GGW2011}. About~40 years later, \cite{G1979} solved the problem by showing that the optimal policy has an index structure, \ie at each time one can compute an index (a real number) solely based on the current state of an arm and choosing the arm associated with the highest index is optimal. Besides Gittins' original proof of optimality based on an interchange argument, \cite{W1980} gave a proof by introducing {\em retirement option} which was further generalized to the restless MAB model. \cite{W1992} gave a beautiful proof without any mathematical equation by an argument of fair charge, while \cite{BN1996} took the achievable region approach for a proof based on linear programming and the duality theory. These four classical proofs of the optimality of Gittins index were elegantly summarized and extended by \cite{FW2016}.

\subsection{Whittle's Generalization to Restless MAB}

\cite{W1988} generalized the classical MAB to the {\em restless} bandit model, where each unselected arm can also change state (accordingly to another known Markovian rule) and offer reward. Furthermore, the player is not restricted to select only one arm but can choose $M~(M<N)$ of them at each time. Either extension of the above makes Gittins index suboptimal in general. Whittle introduced an index policy based on the idea of {\em subsidy}, \ie by focusing on a single-armed bandit one can attach a fixed amount of reward (subsidy) to the arm when it is unselected (made passive) and Whittle index is defined as the minimum subsidy that makes selecting (activating) the arm or not {\em equally} optimal at its current state. This subsidy decouples arms for computing Whittle index and is reduced to Gittins index in the classical MAB model. Whittle showed that the subsidy is essentially the Lagrangian multiplier associated with a relaxed constraint on the {\em expected} number of arms to activate over the infinite horizon, thus providing an upper bound for the original problem. However, there is a great challenge before we can apply Whittle index policy, namely, the {\em indexability} condition. In other words, we require that the subsidy, which makes actions indifferent, exists and is {\em uniquely} defined for each state of each arm. In this case, we call the RMAB is {\em indexable} in which the Whittle index is well-defined. However, proving indexability is generally difficult even for RMAB with finite state spaces \citep{NinoMora2001}. Furthermore, even the indexability is proved to hold, solving for the Whittle index in closed-form is again a difficult problem in the design of an implementable policy \citep{LZ2010,LWZ2011}. See Sec.~\ref{sec:relate} for more details.

\subsection{Resource Constraint and Partial Observability}

The restless MAB (RMAB) is a special class of Markov Decision Processes (MDP) where system state vector is completely observed at the beginning of each decision epoch. However, many problems do not possess such a perfect observation model. Instead, only the selected arms will reveal their states to the player {\em after} arm selection is determined. This category of problems belongs to the class of {\em Partially Observable} MDP (POMDP), which encompasses a much wider application range than MDP \citep{S1978}. In this paper, the~$N$ processes (arms) are modeled as Markov chains evolving over time, according to potentially different rules for state transitions and reward offering. At each time, the player chooses only~$M~(M<N)$ arms to observe and obtain reward determined by the observed states of chosen arms. The states of other unchosen arms remain {\em unknown}. To formulate the problem as an RMAB, we can use {\em information state} as a sufficient statistics for optimal control that characterizes the probability distribution of arm states based on past observations. For the case that each Markov chain has only~$2$ states, the problem was solved near-optimally by Whittle index policy \citep{LZ2010,LWZ2011}. This paper extends those results to the case of $K$-state Markov chains for~$K>2$. As shown in the rest of this paper, this extension makes the problem fundamentally more complex. Our approach is to embrace a family of {\em threshold policies} that significantly simplifies the system dynamics while keeping the major benefits from the fundamental structure of Whittle's relaxation. We summarize the main results of this paper in the next subsection.

\subsection{Main Results}
First, we formulate the problem as a partially observable RMAB with an infinite state space. Second, we establish an equivalent condition for indexability in our problem which further leads to a proof of indexability when the discount factor $\beta\le0.5$. Third, we extend the classical indexability proposed by Whittle to the relaxed indexability. With this generalization, we propose a threshold policy on a single arm that linearizes the original decision boundary and leads to a closed-form expression of the approximate Whittle index under the relaxed indexability. Meanwhile we show that the relaxed indexability relative to the linearized threshold function is reduced to the classical indexability with zero approximation error of the Whittle index for $K=2$. Fourth, we establish an efficient algorithm based on the relaxed indexability and the approximate Whittle index for general $K>2$. Last, we consider the special case of $K=3$ and further optimize the implementation of our algorithm with its near-optimal performance demonstrated by numerical experiments.

\subsection{Related Work}\label{sec:relate}

By considering a large deviation theory applied on Markov jump processes under the time-average reward criterion, \cite{WW1990} showed that Whittle index policy implemented under the strict constraint (\ie choosing exactly~$M$ arms with the highest indices at each time) converges to the upper bound with the relaxed constraint per-arm-wise as~$N\rightarrow\infty$ with~$M/N$ fixed under a sufficient condition. This sufficient condition requires the global stability of a deterministic fluid dynamic system approximating the stochastic state evolution processes of all arms. \cite{WW1991} further showed that the sufficient condition is satisfied when the cardinality of arm states is~$2$ or~$3$. \cite{Verloop2016} extended these results to a wider class of indexable and also non-indexable restless bandits with finite state spaces. However, verifying this sufficient condition is very difficult without a general theoretical approach. The existence of such an index policy (\ie indexability) is part of the sufficient condition and is itself without a general way to verify. For RMAB with finite state spaces, some sufficient conditions for indexability were established (see, \eg \citealt{WW1990,WW1991,NinoMora2001}) as well as some necessary ones (see, \eg \citealt{WW1990,NinoMora2007}). For the indexable RMAB problems studied so far, Whittle index policy has been shown a near-optimal performance in different application areas (see, \eg \citealt{NinoMora2001,GKO2009,LZ2010,HG2011,Verloop2016}). Furthermore, based on Whittle's original idea of arm-decoupling, various index policies have been proposed for restless bandits with finite state spaces with asymptotic optimality proved under certain conditions and a strong performance numerically demonstrated in finite regimes (see, \eg \citealt{BN2000,HF2017,ZJW2019,BS2020,GGY2021}). It is worth noticing that the general restless MAB with a {\it finite} state space is PSPACE-HARD \citep{PT1999}, making it unlikely to discover an efficient optimal algorithm in general.

The partially observable restless bandit for the case of $K=2$ was first formulated in the context of communications networks by \cite{LZ2008} and in the context of unmanned aerial vehicles by \cite{LENY2008}, where the indexability and the closed-form Whittle index function were established under the total discounted reward criterion in both of the two independent papers. \cite{LZ2010} extended these results to the time-average reward criteria and proved the structure, optimality and equivalence of the Whittle index policy to the myopic policy for homogeneous arms (\ie arms with the same $2\times2$ transition probability matrix and reward function). Following these results, various partial observation and state transition models for $K=2$ were studied in different application areas with the strong performance of such an index policy successfully demonstrated (see, \eg \citealt{LZ2010Two,LA2011,LZ2012,W2014,EL2018,Z2019,L2024}). This motivates us to consider the general case of $K>2$ in this paper.

\section{RMAB Formulation and Classical Indexability}

In this section, we will formulate the multi-armed bandit problem as a partially observable Markov decision process and introduce the concept of Whittle Index. Consider a bandit machine
with totally~$N$ independent arms, each of which is modelled as a Markov process. For the $n$-th arm $(n\in\{1,...,N\})$, let $\textbf{P}^{(n)} = \{p_{i,j}^{(n)}\}_{i,j\in\{0,1,2,\ldots,K_n-1\}}$ denote its state transition matrix~and $B_{n,i}~(i\in\{0,1,2,\ldots,K_n-1\})$ the reward that can be obtained when the arm is observed in state~$i$. Let $B_n = [B_{n,0},B_{n,1},B_{n,2},\ldots,B_{n,K_n-1}]$ be the reward vector for arm~$n$. At each discrete time~$t$, $M$ arms will be selected for observation (activated). Let $U(t)\subseteq\{1,...,N\}~(|U(t)|=M)$ be the set of arms that are observed at time~$t$. The (random) reward obtained at time~$t$ is given by
\begin{equation}
	R_{U(t)}(t) = \sum_{n\in U(t)}B_{n,S_n(t)},
	\label{def: immediate reward}
\end{equation}
where $S_n(t)\in\{0,1,2,\ldots,K_n-1\}$ denotes the state of arm~$n$ at time~$t$. Our objective is to decide an optimal policy~$\pi^*$ of choosing~$M$ arms at each time such that the long-term reward is maximized in expectation. In this paper, we will focus on the expected total discounted reward objective function:
\begin{equation}
	\pi^* = \arg\max_{\pi\in\Pi}\mathbb{E}_{\pi}[\sum_{t=1}^{\infty}\beta^{t-1}R_{U(t)}(t)],
	\label{def: discounted}
\end{equation}
where $\beta\in(0,1)$ is the discount factor for the convergence of the sum in the right-hand side of~\eqref{def: discounted} and~$\Pi$ the set of all feasible policies satisfying~$|U(t)|=M$ at each time~$t$.

\subsection{Belief Vector as System State}

Since no arm state is observable before~$U(t)$ is decided at time~$t$, we need an alternative representation of information for decision making. According to the general POMDP theory, the conditional probability distribution of the Markovian state given all past knowledge is a {\em sufficient statistics} for decision making~\citep{S1978}. Specifically, in our problem, the past knowledge consists of the initial (a prior) probability distribution of the state of each arm at $t=1$, the time of last observation of each arm, and the observed state at the last observation of each arm. Then the conditional probability distribution of each arm's state given the past knowledge can be written in the following equation and is referred to as the {\em belief state} (or belief vector) of the arm. The belief states from all arms thus form a sufficient statistics for our decision making process and are {\em fully} observable.

Denoted by $\omega_n(t)$ the belief vector of arm~$n$ at time~$t$, we have
$$\omega_n(t) =
\left(
\begin{matrix}
	{\rm Pr}(S_n(t) = 0|\omega_n(1), \tau_n, S_n(\tau_n))\\
	{\rm Pr}(S_n(t) = 1|\omega_n(1), \tau_n, S_n(\tau_n))\\
    \vdots\\
	{\rm Pr}(S_n(t) = K_n-1|\omega_n(1), \tau_n, S_n(\tau_n))
\end{matrix}\right)',
\quad \Omega(t) = \left(
\begin{matrix}
	\omega_1(t)\\\vdots\\\omega_N(t)
\end{matrix}
\right),$$
where~$A'$ denotes the transpose of~$A$ and~$\tau_n$ the time of last observation on arm~$n$. If the arm has never been observed, we can set~$\tau_n=-\infty$ and remove~$S_n(\tau_n)$ from the condition. Thus the initial belief vector~$\omega_n(1)$ can be set as the stationary distribution~$\bar{\omega}_n$ of the internal Markov chain (corresponding to the case of~$\tau_n=-\infty$)\footnotemark\footnotetext{Here we assume the Markov chain with transition matrix~$\textbf{P}^{(n)}$ is irreducible and aperiodic.}:
\begin{eqnarray}
\omega_n(1) = \bar{\omega}_n=\lim_{k\rightarrow\infty}  \textbf{p}(\textbf{P}^{(n)})^{k},\label{eq:limDist}
\end{eqnarray}
where~$\bar{\omega}_n$ is the unique solution to~$\omega \textbf{P}^{(n)}=\omega$ and $\textbf{p}$ an {\em arbitrary} probability distribution of the state of arm~$n$. The limit in~\eqref{eq:limDist} can be taken under any norm since belief vectors are in a finite-dimensional vector space. It is also convenient to update the belief vector of each arm at each time according to the following {\em Markovian} rule:
\begin{equation}\label{eq: state dist update}
	\omega_n(t + 1) =
	\begin{cases}
		[p^{(n)}_{S_n(t),0},p^{(n)}_{S_n(t),1},\ldots,p^{(n)}_{S_n(t),K_n-1}], &n\in U(t)\\
		\omega_n(t)\textbf{P}^{(n)}, &n \notin U(t)
	\end{cases}\quad.
\end{equation}
So the POMDP problem is reduced to an MDP one by treating all belief vectors of all arms as the system state of the decision problem. However, the state space becomes {\em infinite} as a function space (consisting of probability measures).

Note that the belief update is deterministic if the arm is not chosen for observation at the time. For the case where the arm is not being observed for a consecutive sequence of time, we define the following operator for updating the belief vector continuously over~$k$ consecutive slots {\em without} any observation:
\begin{align}\label{eq: k-step state dist update}
	\Tau^k_n(\omega_n(t)) &= \left(
	\begin{matrix}
		{\rm Pr}(S_n(t+k) = 0 | \omega_n(t))\\
		{\rm Pr}(S_n(t+k) = 1 | \omega_n(t))\\
        \vdots\\
		{\rm Pr}(S_n(t+k) = K_n-1 | \omega_n(t))
	\end{matrix}
	\right)'
	\nonumber\\ &= \omega_n(t)(\textbf{P}^{(n)})^k.
\end{align}
Now the decision problem has a {\em countable} state space as modelled by the belief vector for a {\em fixed} initial~$\Omega(1)$ and an {\em uncountable } state space for an {\em arbitrarily} chosen~$\Omega(1)$. This infinite-dimensional optimization problem can be formulated as
\begin{eqnarray}
&\max_{\pi\in\Pi}&\mathbb{E}_{\pi}[\sum_{t=1}^{\infty}\beta^{t-1}\sum_{n=1}^N\mathbbm{1}(n\in U(t))B_{n,S_n(t)}|\Omega(1)]\label{max:dsctStrict}\\
&\mbox{s. t.}&\sum_{n=1}^N\mathbbm{1}(n\in U(t)) = M,\quad \forall~t\ge1.\label{ys:dsctStrict}
\end{eqnarray}
It is clear that as the number of arms increases, the number of choices at each time grows geometrically. Furthermore, different choices lead to different updates of the belief vector, yielding a high complexity in solving the problem. In the following, we will extend Whittle's original idea of arm-decoupling for an index policy to our model which has an infinite state space consisting of belief vectors.

\subsection{Definition of Indexability and Whittle Index}
Whittle relaxed the strict constraint on the exact number of arms to choose at each time to requiring only~$M$ arms are chosen in {\em expectation}. Particularly, we consider the following {\em relaxed} form of problem~\eqref{max:dsctStrict}:
\begin{eqnarray}
&\max_{\pi\in\Pi}&\mathbb{E}_{\pi}[\sum_{t=1}^{\infty}\beta^{t-1}\sum_{n=1}^N\mathbbm{1}(n\in U(t))B_{n,S_n(t)}|\Omega(1)]\label{max:dsctRelax}\\
&\mbox{s. t.}&\mathbb{E}_{\pi}[\sum_{t=1}^{\infty}\beta^{t-1}\sum_{n=1}^N\mathbbm{1}(n\notin U(t))|\Omega(1)]=\frac{N-M}{1-\beta}.\label{ys:dsctRelax}
\end{eqnarray}
{\bf Remark.} Note that~\eqref{ys:dsctRelax} indeed relaxes~\eqref{ys:dsctStrict} since any policy that selects exactly $M$ arms to activate (for those $n\in U(t)$) at each time necessarily excludes the remaining $N-M$ arms (for those $n\notin U(t)$). In other words, \eqref{ys:dsctRelax} holds even without the expectation operator if~\eqref{ys:dsctStrict} holds because then $\sum_{n=1}^N\mathbbm{1}(n\notin U(t))=N-M$ for any~$t\ge1$. For RMAB with {\em finite} state spaces and the {\em time-average} reward criterion, \cite{WW1990} showed that the performance gap from \eqref{max:dsctStrict}-\eqref{ys:dsctStrict} to \eqref{max:dsctRelax}-\eqref{ys:dsctRelax} asymptotically tends to zero per-arm-wise as $N\rightarrow\infty$ with~$M/N$ fixed (Theorem~1 in \citealt{WW1990}). They also showed that the performance gap induced by the Whittle index policy is determined by the stability of a high-dimensional nonlinear dynamic system (the fluid approximation), which is still an open problem in general for arm state number greater than~$3$ (\citealt{WW1991,Verloop2016}). For RMAB with {\em finite} state spaces and the {\em discounted} reward criterion, a general LP (linear programming) relaxation with the performance region approach was proposed by \cite{BN2000} to numerically demonstrate the small performance gap of the primal-dual index heuristic; while other index heuristics under various relaxation methods for finite time horizons were proposed with performance gap tending to zero (in the same sense as in \citealt{WW1991}) under certain conditions (\citealt{HF2017,ZJW2019,BS2020,GGY2021}). However, these approaches {\em cannot} be directly applied to analyze our problem which has an {\em infinite} state space.

Whittle's relaxation from~\eqref{ys:dsctStrict} to~\eqref{ys:dsctRelax} allows us to analyze the dual problem with arms decoupled as detailed below. Applying the Lagrangian multiplier~$\lambda$ to~\eqref{ys:dsctRelax}, we arrive at the following unconstrained optimization problem:
\begin{eqnarray}
&\max_{\pi\in\Pi}&\mathbb{E}_{\pi}[\sum_{t=1}^{\infty}\beta^{t-1}\sum_{n=1}^N[\mathbbm{1}(n\in U(t))B_{n,S_n(t)}+\lambda\mathbbm{1}(n\notin U(t))]|\Omega(1)].\label{max: decouple}
\end{eqnarray}
The above unconstrained optimization is equivalent to~$N$ independent optimization problem as shown below:
\begin{eqnarray}
&\max_{\pi\in\Pi}&\mathbb{E}_{\pi}[\sum_{t=1}^{\infty}\beta^{t-1}[\mathbbm{1}(n\in U(t))B_{n,S_n(t)}+\lambda\mathbbm{1}(n\notin U(t))]|\omega_n(1)],\quad n=1,2,\ldots,N.\label{max: single}
\end{eqnarray}
Therefore, it is sufficient to consider a single arm for solving problem~\eqref{max: decouple}. Note that the action applied on a single arm is either ``selected (activated)'' or ``unselected (made passive)'' at each time. We can thus focus on the single-armed problem (with Lagrangian multiplier~$\lambda$) with state space consisting of all probability measures on the Markov chain and a binary action space.

For simplicity, we will drop the subscript~$n$ in consideration of a single-armed bandit without loss of generality. Now the Lagrangian multiplier $\lambda$ in~\eqref{max: single} can be seen as a reward of making the single arm passive (unselected). More generally, we define a reward function~$m(t)$ as the reward obtained at time~$t$ when the arm is made passive. With this new definition added into our model, we have a single-armed bandit problem with {\em subsidy for passivity} (\ie~$m(t)$) originally introduced by Whittle. The main reason that we avoid to use the symbol~$\lambda$ in defining~$m(t)$ is because~$\lambda$ is a {\em global} variable in the context of multiple arms in~\eqref{max: decouple}. In other words, if we set~$m(t)=\lambda^*$ for all arms where~$\lambda^*$ is the optimal Lagrangian multiplier achieving the relaxed constraint~\eqref{ys:dsctRelax}, the relaxed multi-armed bandit problem~\eqref{max: decouple} is solved after the solution to each single-armed bandit with subsidy. However, we now focus on the single-armed bandit problem with subsidy~$m(t)$ as a {\em local} reward function for this arm only. For generality of the single-armed model, it is thus better to use a different symbol in distinction to the global Lagrangian multiplier shared by all arms. For our purpose, it is sufficient to consider~$m(t)$ independent of~$t$ and we will simply denote it by~$m$. A minor reason to use the letter~$m$ specifically is to make our notations consistent with those adopted in \cite{LZ2010} for $K=2$. Let $V_{\beta,m}(\omega)$ denote the value of~\eqref{max: single} with $\omega_n(1)=\omega$. It is straightforward to write out the dynamic equation of the single-armed bandit problem as follows:
\begin{equation}\label{eq: long time value function}
	V_{\beta,m}(\omega) = \max\{V_{\beta,m}(\omega;u=1);V_{\beta,m}(\omega;u=0)\},
\end{equation}
where $V_{\beta,m}(\omega;u=1)$ and $V_{\beta,m}(\omega;u=0)$ denote, respectively, the maximum expected total discounted reward that can be obtained if the arm is activated or made passive at the current belief state~$\omega$, followed by an optimal policy in subsequent slots. Since we consider the infinite-horizon problem, a stationary optimal policy can be chosen and the time index~$t$ is not needed in~\eqref{eq: long time value function}. Let~$p_{i\cdot}=[p_{i0},p_{i1},\ldots,p_{i(K-1)}],(i=0,1,\ldots,K-1)$ denote the~$i$-th row of~$\textbf{P}$, we have
\begin{align}\label{eq: long time value function with active action}
	&V_{\beta,m}(\omega;u=1) = \omega B' + \beta \omega
	\left(
	\begin{matrix}
		V_{\beta,m}(p_0)\\
        V_{\beta,m}(p_1)\\
        \vdots\\
        V_{\beta,m}(p_{K-1})
	\end{matrix}
	\right),\\
    &V_{\beta,m}(\omega;u=0) = m + \beta  V_{\beta,m}(\Tau^1(\omega)),\label{eq: long time value function with passive action}
\end{align}
where~$\Tau^1(\omega)$ is the one-step belief update as defined in (\ref{eq: k-step state dist update}). Without loss of generality, we assume~$0= B_0\le B_1\le\cdots\le B_{K-1}$. Note that $V_{\beta,m}(\omega;u=1)$ is linear in~$\omega$ while $V_{\beta,m}(\omega;u=0)$ is convex in~$\omega$ as shown by Lemma~\ref{lm:vf} in Sec.~\ref{sec:threshOpt}.

Define {\em passive set}~$P(m)$ as the set of all belief states such that taking the passive action~$u=0$ is optimal:
\begin{eqnarray}
P(m)\defeq \{\omega:~V_{\beta,m}(\omega;u=1)\le V_{\beta,m}(\omega;u=0)\}.\label{def:passive}
\end{eqnarray}
It is clear that~$P(m)$ changes from the empty set to the whole space of probability measures as~$m$ increases from~$-\infty$ to~$\infty$. However, such change may not be monotonic as~$m$ increases (see Sec.~\ref{sec:threshOpt} for more discussions). If the passive set~$P(m)$ increases monotonically with~$m$, then for each value~$\omega$ of the belief state, one can define the {\em unique}~$m$ that makes it join~$P(m)$ and stay in the set forever. Intuitively, this~$m$ measures how attractive it is to activate the arm at the belief state~$\omega$ compared to other belief states in a well-ordered manner: the larger~$m$ required for it to be passive, the more incentives to activate at the belief state without~$m$. This value of $m$ (if well-defined) thus yields a priority index of the belief state. In the following, we present the formal definition of {\em indexability} and {\em Whittle index} \citep{W1988}.

\begin{definition}
A restless multi-armed bandit is {\em indexable} if for each single-armed bandit with subsidy, the passive set of arm states increases {\em monotonically} from~$\emptyset$ to the whole state space as~$m$ increases from ~$-\infty$ to~$+\infty$. Under indexability, the {\em Whittle index} of an arm state is defined as the infimum subsidy~$m$ such that the state remains in the passive set.
\end{definition}

Note that if the indexability condition is verified and the Whittle index solved as a function of the state of each arm, the Lagrangian relaxation problem~\eqref{max: decouple} may be solved with the optimal~$\lambda^*$: for each arm at each time, we choose to activate the arm if its current Whittle index is greater than~$\lambda^*$ or make it passive otherwise. There is some {\em randomization} technique involved to ensure the satisfaction of constraint~\eqref{ys:dsctRelax} when the Whittle index is equal to~$\lambda^*$. But that is not the main focus of this paper and we will give some brief discussions following Theorem~\ref{prop:indexability} in Sec.~\ref{sec:eqvIdx}.

\subsection{Threshold Structure of The Optimal Policy}\label{sec:threshOpt}

For our model in which the arm state is given by the belief vector, the indexability is equivalent to the following:
\begin{eqnarray}
\mbox{If}~V_{\beta,m}(\omega;u=1)\le V_{\beta,m}(\omega;u=0),~\mbox{then}~\forall~m'>m,~V_{\beta,m'}(\omega;u=1)\le V_{\beta,m'}(\omega;u=0).
\end{eqnarray}
Under indexability, the Whittle index~$W(\omega)$ of arm state~$\omega$ is defined as
\begin{eqnarray}
W(\omega)\defeq \inf\{m:~V_{\beta,m}(\omega;u=1)\le V_{\beta,m}(\omega;u=0)\}.\label{def:whittleIdx}
\end{eqnarray}

Before we proceed, it helps to emphasize on the {\em recursive} nature in defining the value functions given in \eqref{eq: long time value function with active action} and \eqref{eq: long time value function with passive action} conditional on the active and passive actions, respectively. We know that the indexability condition essentially requires a {\em one-time only} rank change of the two value functions as $m$ increases. Although it is intuitive that the larger subsidy causes more states to join the passive set, we cannot conclude this by merely comparing the immediate rewards obtained by active and passive actions ($\omega B'$ vs. $m$) respectively: the future total expected reward is again in the form of value functions that are {\em dependent on} our current action (which affects the belief update) {\em and} the subsidy~$m$. To evaluate indexability, we need to have sufficient knowledge about the value functions \eqref{eq: long time value function with active action} and \eqref{eq: long time value function with passive action} to determine their rank (as functions of the current belief state and the subsidy). In general, the value functions are hard to solve due to the dilemma between exploitation and exploration mentioned at the beginning of this paper. However, for the problem at hand, we can show that the value function \eqref{eq: long time value function} implies a threshold structure of the problem, which generalizes the case of $K=2$ and further inspires the development of an efficient algorithm, as detailed below.

Now we prove a crucial lemma that gives some fundamental properties of the value function~$V_{\beta,m}(\omega)$.
\begin{lemma}\label{lm:vf}
The value function~$V_{\beta,m}(\omega)$ for the single-armed bandit with subsidy is convex and Lipschitz continuous in both~$\omega$ and~$m$.
\end{lemma}

The proof will be given in the e-companion to this paper.

\noindent\textbf{Remark}
\begin{itemize}
  \item Note that if $m\le 0$, it is optimal to always activate the arm (since all extreme points of a convex function under the passive action are below those of a linear one under the active action) and $V_{\beta,m}(\omega)$ does not depend on~$m$ and is thus Lipschitz continuous in~$m$. If $m\ge B_{K-1}$, it is optimal to always make the arm passive so $V_{\beta,m}(\omega)=\frac{m}{1-\beta}$ and is thus Lipschitz continuous as well. The interesting case is when $0< m< B_{K-1}$ as focused in the rest of the paper. The monotonic property of $V_{\beta,m}(\omega)$ as a nondecreasing function of~$m$ is clear.
  \item Since $V_{\beta,m}(\omega)$ is Lipschitz continuous in~$m$, it is also absolutely continuous and differentiable almost everywhere in~$m$. Assume~$m_0$ is a point where the derivative exists, a small increase to $m_0+\Delta m$ should cause $V_{\beta,m_0}(\omega)$ to boost at a ratio at least equal to the expected total discounted time of being passive, since the subsidy~$m_0$ for passivity is being paid for such a duration of time (passive time in short). The passive time is not necessarily unique and we will give a rigorous formulation of its relation to the (right) derivative of $V_{\beta,m}(\omega)$ in Theorem~\ref{prop:indexability} in Sec.~\ref{sec:eqvIdx}.
  \item Since $V_{\beta,m}(\omega)$ is also Lipschitz continuous in~$\omega$, for sufficiently small~$\beta$, a change of~$\omega$ that makes the immediate reward~$\omega B'$ vary may play a dominating role in determining the order of~\eqref{eq: long time value function with active action} and~\eqref{eq: long time value function with passive action} as the value function $V_{\beta,m}(\omega)$ varies in bounded ratios with~$\omega$. This motivates us to consider the family of {\em linearized threshold policies}: following the trajectory of~$\Tau^k(\omega)$ until some {\em linear} function $r(\cdot):~\mathbb{R}^K\rightarrow\mathbb{R}$ (\eg the projection $r(\omega)=\omega B'$) maps $\Tau^k(\omega)$ to a value greater than a given one, we activate the arm and reset the value function to one of $V_{\beta,m}(p_0),V_{\beta,m}(p_1),\cdots,V_{\beta,m}(p_{K-1})$ (See~\eqref{eq: long time value function with active action}). Linearized threshold policies are suboptimal in general, especially when~$\beta$ is large. However, they provide an efficient way in solving the approximated value functions and leads to a computable Whittle index function in low complexity and near-optimal performance even when~$\beta$ is close to~$1$, as elaborated in Sec.~\ref{sec:relaxedIdx} and Sec.~\ref{sec:numeric}.
\end{itemize}

Next, we show that the optimal single-arm policy has a general threshold structure. Let $\mathbb{X}$ denote the belief state space as a $(K-1)$-simplex. It is a $(K-1)$-dimensional space of probability measures. For convenience, we still use the $K$-dimensional vector $(\omega_1,\omega_2,\ldots,\omega_K)$ to denote a point in $\mathbb{X}$ by keeping in mind that $\sum_{i=1}^K\omega_k=1$. Now consider an extreme point~$\omega=[0,0,\cdots,1,\cdots,0]$ of the belief state space where it is known that the arm's internal state is~$k$ for some~$k\in\{0,1,\cdots,K-1\}$. In this case, the next belief state is deterministically~$p_k$, independent of the current action, \ie
\begin{eqnarray}
\arg\max\{V_{\beta,m}(\omega;u=1);V_{\beta,m}(\omega;u=0)\}&=& \arg\max\left\{B_k, ~ m\right\}.
\end{eqnarray}
From the above, each extreme point successively joins the passive set as~$m$ increases from~$0$ to~$B_{K-1}$. Consider an~$m\in(0,B_{K-1})$ such that~$0=B_0\le\cdots\le m<B_j\le\cdots\le B_{K-1}$. The first~$j$ states are in the passive set while states~$j,\cdots,K-1$ are in the active set defined as
\begin{eqnarray}
A(m)\defeq \mathbb{X}-P(m)= \{\omega:~V_{\beta,m}(\omega;u=1)> V_{\beta,m}(\omega;u=0)\}.\label{def:active}
\end{eqnarray}

The following lemma shows that the active set $A(m)$ is an open convex region in $\mathbb{X}$ with a {\em decision boundary} $C(m)$ shared by the passive set $P(m)$:
\begin{eqnarray}
C(m)\defeq \{\omega:~V_{\beta,m}(\omega;u=1)=V_{\beta,m}(\omega;u=0)\}.\label{def:boundary}
\end{eqnarray}

\begin{lemma}\label{lm:boundary}
The active set $A(m)$ is an open convex $(K-1)$-dimensional subspace of $\mathbb{X}$. The decision boundary $C(m)$ is a compact and simply connected $(K-2)$-dimensional subspace of $\mathbb{X}$ that partitions $\mathbb{X}$ into two disjoint and connected subspaces: $A(m)$ and $P(m)$ with $C(m)\subset P(m)$.
\end{lemma}

The proof will be given in the e-companion to this paper.

According to Lemma~\ref{lm:boundary}, we can treat $C(m)$ as a $(K-2)$-dimensional threshold {\em without} any holes or discontinuities for the optimal decision making process. One can visualize it as a {\it curve} for $K=3$ or a {\it surface} for $K=4$. Higher dimensions are analogous to compact $(K-2)$-manifolds. If indexability holds, the boundary~$C(m)$ should (continuously) move in a direction such that~$A(m)$ shrinks as~$m$ increases. For each~$\omega$, there exists an~$m$ such that~$C(m)$ reaches ~$\omega$ for the first time and this~$m$ is the Whittle index~$W(\omega)$ of~$\omega$:
\begin{eqnarray}
W(\omega)\defeq \inf\{m:V_{\beta,m}(\omega;u=1)\le V_{\beta,m}(\omega;u=0)\}=\min\{m:\omega\in C(m)\}.\label{def:WI}
\end{eqnarray}
In the above, we have used the minimization operator instead of the infimum by observing that the closure of the nontrivial region~$(0,B_{K-1})$ for the subsidy~$m$ is compact. A {\em sufficient and necessary} condition of indexability for our model with an infinite state space is given in the next subsection.

\subsection{An Equivalent Condition for Indexability}\label{sec:eqvIdx}

In this subsection, we establish a sufficient and necessary condition for indexability by requiring the decision boundary $C(m)$ to satisfy certain properties. Furthermore, we verify this equivalent condition to prove the indexability of our problem when $\beta\le0.5$.

\begin{theorem}\label{prop:indexability}
Let~$\Pi^*_{sa}(m)$ denote the set of all optimal single-arm policies achieving~$V_{\beta,m}(\omega)$ with initial belief state~$\omega$. Define the passive time
\begin{eqnarray}
D_{\beta,m}(\omega)\defeq \max_{\pi^*_{sa}(m)\in\Pi^*_{sa}(m)}\mathbb{E}_{\pi^*_{sa}(m)}[\sum_{t=1}^{\infty}\beta^{t-1}\mathbbm{1}(u(t)=0)|\omega(1)=\omega].\label{def:passiveTime}
\end{eqnarray}
The right derivative of the value function~$V_{\beta,m}(\omega)$ with~$m$, denoted by~$\frac{dV_{\beta,m}(\omega)}{(dm)^+}$, exists at every value of~$m$ and
\begin{eqnarray}
\left.\frac{dV_{\beta,m}(\omega)}{(dm)^+}\right|_{m=m_0}=D_{\beta,m_0}(\omega).\label{eq:rightDiff}
\end{eqnarray}
Furthermore, the single-armed bandit is indexable if and only if for all values of~$\omega$ and~$m_{\omega}$ such that~$\omega\in C(m_{\omega})$, we have
\begin{eqnarray}
\left.\frac{dV_{\beta,m}(\omega;u=0)}{(dm)^+}\right|_{m=m_{\omega}} \ge \left.\frac{dV_{\beta,m}(\omega;u=1)}{(dm)^+}\right|_{m=m_{\omega}},\label{eq:diffIdx}
\end{eqnarray}
and for any~$\omega\in C(m_{\omega})$ with the equality true in~\eqref{eq:diffIdx}, there exists a~$\Delta m(\omega)>0$ such that
\begin{eqnarray}
V_{\beta,m}(\omega;u=0)\ge V_{\beta,m}(\omega;u=1),~\quad\forall~m\in (m_{\omega}+\Delta m(\omega)).\label{eq:diffIdx1}
\end{eqnarray}
\end{theorem}

The proof will be given in the e-companion to this paper.

\noindent\textbf{Remark.} Theorem~\ref{prop:indexability} establishes a crucial relation between the value function $V_{\beta,m}(\omega)$ and the passive time $D_{\beta,m}(\omega)$ as its right derivative. The convexity established in Lemma~\ref{lm:vf} then implies the monotonic property of $D_{\beta,m}(\omega)$ as~$m$ increases. However, the increase of $D_{\beta,m}(\omega)$ needs not be continuous. In the proof of Theorem~\ref{prop:indexability}, we have shown the right continuity of $D_{\beta,m}(\omega)$ but not the left one. These jumping points are essentially caused by the case where the points in the belief state space may not join the passive set~$P(m)$ in a continuous sense. Specifically, if we fix the initial belief state~$\omega$, the state of the arm will move in a countable set as a discrete process. Under the optimal policy that achieves the passive time defined in~\eqref{def:passiveTime} and indexability, it is possible that when~$m$ increases by a sufficiently small amount, the policy remains unchanged, \ie the partition of active and passive sets for the countable state space is the same. Consequently, the passive time $D_{\beta,m}(\omega)$ remains a constant during this increasing period of the subsidy. However, as~$m$ keeps increasing, new states would join the passive set and cause a jump in $D_{\beta,m}(\omega)$. The discontinuity of $D_{\beta,m}(\omega)$ poses a question: how should one make the continuation of $D_{\beta,m}(\omega)$ such that constraint~\eqref{ys:dsctRelax} must be satisfied for the relaxed version of the multi-armed bandit problem? The technique is to use nondeterministic optimal policies: for believe states in the decision boundary~$C(m)$ that causes discontinuities in $D_{\beta,m}(\omega)$, we activate the arm with certain probability $\rho\in[0,1]$ and make it passive with probability~$1-\rho$. As~$\rho$ decreases from~$1$ to~$0$, the corresponding policies provide a continuation of $D_{\beta,m}(\omega)$. For a detailed exposition of this randomization technique in solving the original {\em multi-armed} bandit problem under the relaxed constraint, see \citealt{L2020} that considers a more general model of infinite arm state spaces.

\begin{theorem}\label{thm:indexableBeta}
The restless bandit is indexable if~$\beta\le 0.5$.
\end{theorem}

The proof will be given in the e-companion to this paper.

However, it is difficult to verify~$\eqref{eq:diffIdx}$ and~$\eqref{eq:diffIdx1}$ when~$\beta>0.5$. This requires further analysis on the passive time~$D_{\beta,m}(\omega)$ as well as the value function~$V_{\beta,m}(\omega)$. If we can characterize the boundary function~$C(m)$ associated with subsidy~$m$, then for each~$\omega$, we may obtain the first crossing time~$L(\omega,C(m))$ when it enters the active set under consecutive passive actions:
\begin{eqnarray}
L(\omega,C(m))\defeq \min_{0\le k<\infty}\{k:~~\Tau^k(\omega)\in A(m)\}.\label{def:firstCrossing}
\end{eqnarray}
Define~$\Tau^0(\omega)\defeq\omega$ and if~$\Tau^k(\omega)\notin A(m)$ for all~$k\ge0$, we set~$L(\omega,C(m))=+\infty$. It is clear that for any~$\omega\in C(m)$, we have
\begin{eqnarray}
    &V_{\beta, m}(\omega)& = \omega B' + \beta\omega(V_{\beta, m}(p_0),\cdots,V_{\beta, m}(p_{K-1}))'\nn\\
                         && = \frac{1 - \beta^{L(\omega, C(m))}}{1 - \beta}m +\beta^{
		L(\omega, C(m))}V_{\beta, m}\left(\Tau^{L(\omega, C(m))}(\omega); u = 1\right),\label{eq: linearE1}\\
    &V_{\beta, m}\left(\Tau^{L(\omega, C(m))}(\omega); u = 1\right)& = \Tau^{L(\omega, C(m))}(\omega)B'+\beta \Tau^{L(\omega, C(m))}(\omega)(V_{\beta, m}(p_0),\cdots,V_{\beta, m}(p_{K-1}))',\label{eq: linearE2}\\
    &V_{\beta, m}(p_k)& = \frac{1 - \beta^{L(p_k, C(m))}}{1 - \beta}m +\beta^{
		L(p_k, C(m))}V_{\beta, m}\left(\Tau^{L(p_k, C(m))}(p_k); u = 1\right),\nn\\&&\quad\forall~k\in\{0,\cdots,K-1\},\label{eq: linearE3}\\
    &V_{\beta, m}\left(\Tau^{L(p_k, C(m))}(p_k); u = 1\right)& = \Tau^{L(p_k, C(m))}(p_k)B'+\beta \Tau^{L(p_k, C(m))}(\omega)(V_{\beta, m}(p_0),\cdots,V_{\beta, m}(p_{K-1}))',\nn\\&&\quad\forall~k\in\{0,\cdots,K-1\}.\label{eq: linearE4}
\end{eqnarray}
With both~$\omega$ and $L(\cdot,C(m))$ {\em fixed and known}, the above equation set is linear and has~$2K+3$ equations with~$2K+3$ unknowns (value functions and~$m$), so an exact solution for the value function $V_{\beta, m}(\omega)$ is possible to obtain, as well as the passive time $D_{\beta, m}(\omega)$, and the subsidy~$m$, in terms of~$\omega$ and $L(\cdot,C(m))$. However, even if $L(\cdot,C(m))$ is solved for, such a way in checking indexability and solving for Whittle index is complex since $L(\cdot,C(m))$ itself is a {\em nonlinear} function and appears as an exponent in the expressions of the value functions. Furthermore, the function~$L(\cdot,C(m))$ may be solved only if~$C(m)$ is sufficiently analyzed which involves dynamic programming on an uncountable state space. To circumvent these difficulties, we consider a family of threshold policies that simplifies the analysis of the value functions and establish an approximation of Whittle index under a relaxed requirement for indexability, as detailed in Sec.~\ref{sec:relaxedIdx}.

Before concluding this subsection, let us have a brief review on the case of $K=2$ as considered in \cite{LZ2010}. In this case, our decision boundary $C(m)$ is reduced to a single point ($0$-dimensional as also proven in Lemma~2 in \citealt{LZ2010})! Given {\em any} $\omega^*$ in the $1$-dimensional belief space (\ie homeomorphic to the interval $[0,1]$), the solution to $L(\omega,\omega^*)$ is straightforwardly obtained in closed-form (Equations (16) and (17) in \citealt{LZ2010}). Consequently, the indexability and closed-form Whittle index can be established with the closed-form solutions of the value functions. But these proofs are highly nontrivial even with the closed-form solutions of the value functions (\citealt{LZ2010}).

For the general case of $K>2$, we do not know $C(m)$ solely by fixing an $\omega^*\in C(m)$ since $C(m)$ contains other (uncountable) points with {\em unknown} locations. Henceforth the first-crossing time $L(\cdot,C(m))$ is not solved. Nevertheless, we are free to impose {\em some relations} among the points in $C(m)$ to approximate $L(\cdot,C(m))$. This is the key motivation for generalizing the classical indexability to a relaxed one, as elaborated in the next section.

\section{Threshold Policies and Relaxed Indexability}\label{sec:relaxedIdx}

As mentioned in Sec.~\ref{sec:eqvIdx}, the key to analyze the indexability is to solve for the first crossing time $L(\omega,C(m))$ and subsequently the set of linear equations for {\em every} fixed~$\omega$ and $L(\omega,C(m))$~(\ref{def:firstCrossing}--\ref{eq: linearE4}). Note that the linear equation set has a coefficient matrix {\em nonlinearly dependent} on $L(\omega,C(m))$ and thus on the belief state. Therefore each belief state requires a {\em different} set of linear equations to solve for its Whittle index, \ie the system of these equations is {\em nonlinear} over the belief space in general. Nevertheless, the first step is still to solve for $L(\omega,C(m))$. In this subsection, we approximate $C(m)$ by a family of {\em linearized} threshold policies to solve for the approximate Whittle index under the framework of {\em relaxed indexability}.

We first give a general definition of threshold policies.  A threshold policy is defined by a mapping~$r(\cdot)$ from the belief state space $\mathbb{X}$ to a space $\mathbb{Y}$ with the {\em order} topology such that the binary action (activate or make passive) at any state~$\omega\in\mathbb{X}$ depends only on the order between~$r(\omega)$ and~$r(\omega^*_{\beta}(m))$ for a {\em fixed} $\omega^*_{\beta}(m)\in\mathbb{X}$. This~$\omega^*_{\beta}(m)$ is called to be {\em in} the threshold of the belief state space with respect to~$r(\cdot)$. Furthermore, we require that the set $\{\omega\in\mathbb{X}:r(\omega)=r(\omega^*_{\beta}(m))\}$ to be {\it simply connected}, \ie it is without any holes or disconnectedness but as a whole solid piece. Then we call this set as the threshold and the mapping $r(\cdot)$ the threshold function that specifies the threshold policy. In our problem, the optimal single-arm policy with a fixed subsidy~$m$ is a threshold policy with its threshold function given by $V_{\beta,m}(\omega;u=1)-V_{\beta,m}(\omega;u=0)$, mapping the belief simplex space $\mathbb{X}$ to the real line $\mathbb{R}$. Recall that the original decision boundary~$C(m)$ defined in~\eqref{def:boundary} is thus a nonlinear~$(K-2)$-dimensional threshold as it is simply connected (Lemma~\ref{lm:boundary}).

To approximate the decision boundary $C(m)$, we consider the following {\em linear} threshold function~$r(\cdot)$ defined as
\begin{eqnarray}
r(\omega)\defeq \omega B'.\label{def:immediateR}
\end{eqnarray}
We immediately observe that the above threshold function {\em linearizes} the original decision boundary $C(m)$: fixing any $\omega^*_{\beta}(m)\in\mathbb{X}$ as a point in the threshold, the entire threshold is specified by the $(K-2)$-dimensional hyperplane
\begin{eqnarray}
\{\omega\in\mathbb{X}:~~\omega B'=\omega^*_{\beta}(m) B'\}.\label{eqn:hyperplane}
\end{eqnarray}
It is easy to see the hyperplane is reduced to a point when $K=2$ which is just $\omega^*_{\beta}(m)$ itself. In this case, the linearized threshold policy is equivalent to the optimal policy \citep{LZ2010}.

The threshold function~$r(\cdot)$ defined in~\eqref{def:immediateR} is a simple and intuitive definition which is identical to the immediate reward function by activating the arm. As the belief state~$\omega$ varies in the~$(K-1)$-dimensional probability space, we measure the attractiveness of activating the arm by the expected reward that can be immediately obtained under activation. For the original problem with multiple arms, if they share the same parameters (homogeneous arms), \ie the transition matrix~$P$ and reward vector~$B$ are arm-independent, this threshold policy on a single arm with subsidy corresponds to the myopic policy: at each time we activate the~$M$ arms that will yield the highest expected reward. However, for inhomogeneous arms, the myopic policy may yield a significant performance loss (see Sec.~\ref{sec:numeric}). It is thus important to precisely characterize the attractiveness of a state as a function of the arm parameters. Our attempt is to solving for the subsidy~$m$ that makes a belief state~$\omega$ as a point in the threshold and define this~$m$ (if exists) as its approximate Whittle index~$W(\omega)$. Given the initial belief state~$\omega$ as the threshold, the action~$u_{\beta,m}^*(\omega(t))$ to take at~$t\ge1$ is given by:
\begin{equation}\label{def:threshold}
	u_{\beta,m}^*(\omega(t)) =
	\begin{cases}
		1~\text{(active)}, \quad \text{if }r(\omega(t)) > r(\omega^*_{\beta}(m))\\
		0~\text{(passive)}, \quad \text{if }r(\omega(t)) \le r(\omega^*_{\beta}(m))
	\end{cases}.
\end{equation}
It is important to observe that when the current arm state is equal to the threshold, \eg at~$t=1$, we always make the arm passive (for now). This is because activating the arm does {\em not} necessarily yield the same performance when {\em confined} in the family of linearized threshold policies. Nevertheless, the suboptimality of this threshold policy is alleviated if the belief update has a sharp slope projected into the real line by~\eqref{def:immediateR} and the discount factor~$\beta$ is small, in which case the comparison in~\eqref{eq: long time value function} is dominated by the order between the expected immediate reward~$\omega B'$ and the subsidy~$m$.

\subsection{The Value Function}
Consider the linearized threshold policy~$\pi_{\beta,m}$ with~$\omega^*_{\beta}(m)$ fixed as a point in the threshold and the belief points $p_0,p_1,\cdots p_{K-1}$. Define
\begin{eqnarray}
L(\omega_1,\omega^*_{\beta}(m))\defeq \min_{0\le k<\infty}\{k:~~\Tau^k(\omega_1)B'> \omega^*_{\beta}(m) B'\}.\label{def:firstCrossingT}
\end{eqnarray}
The above is just a simplified expression by using $\omega^*_{\beta}(m)$ to delegate the entire decision boundary (compare~\eqref{def:firstCrossing}) because the threshold function is specified now. If~$\Tau^k(\omega_1)B'\le \omega_2 B'$ for all~$k\ge0$, we set~$L(\omega_1,\omega_2)=+\infty$. Under~$\pi_{\beta,m}$, we have
\begin{eqnarray}
    &\hat{V}_{\beta, m}(p_k)& = \frac{1 - \beta^{L(p_k, \omega)}}{1 - \beta}m +\beta^{
		L(p_k, \omega)}\hat{V}_{\beta, m}\left(\Tau^{L(p_k, \omega)}(p_k); u = 1\right),\nn\\&&\quad\forall~k\in\{0,\cdots,K-1\},\label{eq: linearET1}\\
    &\hat{V}_{\beta, m}\left(\Tau^{L(p_k, \omega)}(p_k); u = 1\right)& = \Tau^{L(p_k, \omega)}(p_k)B'+\beta \Tau^{L(p_k, \omega)}(\omega)(\hat{V}_{\beta, m}(p_0),\cdots,\hat{V}_{\beta, m}(p_{K-1}))',\nn\\&&\quad\forall~k\in\{0,\cdots,K-1\},\label{eq: linearET2}
\end{eqnarray}
where the value function $\hat{V}_{\beta, m}(\omega_1)$ denotes the expected total discounted reward under~$\pi_{\beta,m}$, starting from a belief state~$\omega_1$. If $L(\cdot,\omega)$ is solved, the above equation set is linear and has~$2K$ unknowns with~$2K$ equations. In this case, we show~\eqref{eq: linearET1} and~\eqref{eq: linearET2} yield a unique solution consisting of $\{\hat{V}_{\beta, m}(p_k)\}_{k=0}^{K-1}$.
\begin{lemma}\label{lemma:linearEq1}
Given the first crossing time function $L(\cdot,\omega)$ with $\omega^*_{\beta}(m)=\omega$ fixed as a point in the threshold, the linear equation set~\eqref{eq: linearET1} and~\eqref{eq: linearET2} has a unique solution consisting of the value functions $\{\hat{V}_{\beta, m}(p_k)\}_{k=0}^{K-1}$ in terms of~$\omega$ and~$m$.
\end{lemma}

The proof will be given in the e-companion to this paper.

\noindent\textbf{Remark}
\begin{itemize}
  \item The proof of Lemma~\ref{lemma:linearEq1} does not require any particular form of~$L(\cdot,\cdot)$ so for any optimal single-arm policy, equations~\eqref{eq: linearE3} and~\eqref{eq: linearE4} have a unique solution consisting of $\{V_{\beta, m}(p_k)\}_{k=0}^{K-1}$ in terms of~$m$ and~$\omega$. For equation~\eqref{eq: linearE1} that solves for the Whittle index~$m$ for a given belief state~$\omega$, the existence of the decision boundary~$C(m)$ under the optimal policy leads to its validity. For the linearized threshold policy with $\omega^*_{\beta}(m)=\omega$ fixed as a point in the threshold, if we add to~\eqref{eq: linearET1} and~\eqref{eq: linearET2} the following additional constraint similar to~\eqref{eq: linearE1}, then the subsidy~$m$ may be solved as an {\em approximated Whittle index} under the threshold policy.
    \begin{eqnarray}
    &\hat{V}_{\beta, m}(\omega)& = \omega B' + \beta\omega(\hat{V}_{\beta, m}(p_0),\cdots,\hat{V}_{\beta, m}(p_{K-1}))' = m + \beta\hat{V}_{\beta, m}(\Tau^1(\omega)).\label{eq:equal at threshold}
    \end{eqnarray}
    Equation~\eqref{eq:equal at threshold} essentially requires that there exists a subsidy~$m$ such that taking active and passive actions at the threshold~$\omega$ achieve the same performance by {\em following the threshold policy}. In general, this might not hold and we need to redefine the approximated Whittle index as detailed in Sec.~\ref{sec:relaxedIdx}.
  \item The value function $\hat{V}_{\beta, m}(\omega)$ with the fixed threshold $\omega^*_{\beta}(m)=\omega$ is a linear function of~$m$ as
  \begin{eqnarray}
  &\hat{V}_{\beta, m}(\omega)& = \frac{1 - \beta^{L(\omega, \omega)}}{1 - \beta}m +\beta^{
		L(\omega, \omega)}(\Tau^{L(\omega, \omega)}(\omega)B'+\beta \Tau^{L(\omega, \omega)}(\omega)(\hat{V}_{\beta, m}(p_0),\cdots,\hat{V}_{\beta, m}(p_{K-1}))')
  \end{eqnarray}
  and $\{\hat{V}_{\beta, m}(p_k)\}_{k=0}^{K-1}$ are linear in~$m$ as well, since $L(\cdot, \omega)$ is independent of~$m$ with~$\omega$ fixed as a point in the linearized threshold. Furthermore, the coefficient of~$m$ in the linear expression of $\hat{V}_{\beta, m}(\omega)$ is equal to the expected total discounted passive time starting from~$\omega$ under the threshold policy. With the indifference at~$\omega$ as given in~\eqref{eq:equal at threshold}, the passive time must be unique and equal to the derivative of $\hat{V}_{\beta, m}(\omega)$ with~$m$ under the threshold policy with $\omega^*_{\beta}(m)=\omega$ fixed as a point in the threshold.
\end{itemize}

\subsection{Relaxed Indexability and Approximate Whittle Index}\label{sec:relaxedIdx}

To have a well-defined value function $\hat{V}_{\beta, m}(\omega)$ to calculate the approximated Whittle index of a belief state~$\omega$, we need one more equation (that is,~\eqref{eq:equal at threshold}) that makes active and passive actions indistinguishable at~$\omega$ fixed as a point in the threshold (Lemma~\ref{lemma:linearEq1}). Now we introduce the following definition of {\em relaxed indexability}:

\begin{definition}
A restless multi-armed bandit satisfies the relaxed indexability {\em with respect to} a threshold policy if for each arm state there exists a unique subsidy~$m$ such that, when making this state as a point in the threshold, then taking the passive and active actions at this state followed by the threshold policy achieve the same performance.
\end{definition}

For our model, relaxed indexability is equivalent to the {\em unique} solution of~\eqref{eq: linearET1},~\eqref{eq: linearET2} and~\eqref{eq:equal at threshold}. Recall the linearity of the value function $\hat{V}_{\beta, m}(\omega_1)$ with~$m$ given a belief state $\omega_1$ fixed as a point in the threshold, we take the derivatives of the value functions in~\eqref{eq:equal at threshold} with $m$ and arrive at the following equivalent condition for relaxed indexability.
\begin{theorem}\label{prop:relaxIdx}
Define the passive time~$\hat{D}_{\beta,m}(\omega_1)$ under a threshold policy~$\pi_{\beta,m}$ as
\begin{eqnarray}
\hat{D}_{\beta,m}(\omega_1)\defeq \mathbb{E}_{\pi_{\beta,m}}[\sum_{t=1}^\infty\beta^{t-1}\mathbbm{1}(u(t)=0)|\omega(1)=\omega_1]=\frac{d\hat{V}_{\beta, m}(\omega_1)}{dm}.\label{eq: passiveTimeTh}
\end{eqnarray}
By fixing $\omega^*_{\beta}(m)=\omega$ as a point in the threshold and an $m$-independent function~$r(\cdot)$ in~\eqref{def:immediateR}, the passive time starting from any initial belief state~$\omega_1$ is independent of~$m$ and denoted by~$\hat{D}_{\beta}(\omega_1)$. The restless bandit of POMDP satisfies the relaxed indexability if and only if for any arm, the corresponding single-armed bandit problem with subsidy~$m$ and with any belief state~$\omega$ fixed as a point in the threshold on the arm state space, we have
\begin{eqnarray}
\beta\omega(\hat{D}_{\beta}(p_0),\cdots,\hat{D}_{\beta}(p_{K-1}))' \neq 1+\beta\hat{D}_{\beta}(\Tau^1(\omega)).\label{eq: relaxedIdxIff}
\end{eqnarray}
Under the relaxed indexability, the approximated Whittle index~$\hat{W}(\omega)$ for a belief state~$\omega$ is given by
\begin{equation}\label{eq: Approximated Whittle}
	\hat{W}(\omega) = \frac{\omega B' - \beta g(\omega\textbf{P})\biggl[
		\textbf{I}_K + \beta H(\textbf{P})G(\textbf{P})
		\biggr]B' + \beta\omega H(\textbf{P})G(\textbf{P})B'}
	{1 + \beta f(\omega\textbf{P}) + \beta\biggl[\beta g(\omega\textbf{P}) - \omega\biggr]
		H(\textbf{P})F(\textbf{P})},
\end{equation}
where $L(\cdot) := L(\cdot,\omega),\ f(\cdot) :=
\frac{1 - \beta^{L(\cdot)}}{1 - \beta},\ g(\omega) := \beta^{L(\omega)}\Tau^{L(\omega)}
(\omega) = \beta^{L(\omega)}\omega\textbf{P}^{L(\omega)}$, and
$$F(\textbf{P}) := \left(
\begin{matrix}
	f(p_0)\\f(p_1)\\\cdots\\f(p_{K-1})
\end{matrix}
\right),\quad G(\textbf{P}) := \left(
\begin{matrix}
	g(p_0)\\g(p_1)\\\cdots\\g(p_{K-1})
\end{matrix}
\right),\quad H(\textbf{P})=\biggl(\textbf{I}_K - \beta G(\textbf{P})\biggr)^{-1}.$$
\end{theorem}

The proof will be given in the e-companion to this paper.

\noindent\textbf{Remark}
\begin{itemize}
  \item Fixing any~$\omega$ as in the threshold and an initial~$\omega_1$, Lemma~3 allows us to solve for $\hat{V}_{\beta, m}(\omega_1)$ as a linear function of~$m$ and its coefficient (derivative) is just $\hat{D}_{\beta,m}(\omega_1)$. Thus the relaxed indexability condition~\eqref{eq: relaxedIdxIff} is {\em equivalent} to requiring that the denominator of the Whittle index~\eqref{eq: Approximated Whittle} obtained from~\eqref{eq:equal at threshold} is nonzero. If there exists a belief state~$\omega$ such that~\eqref{eq: relaxedIdxIff} does not hold, we can simply use~$\omega B'$ as a substitute for its Whittle index to measure the attractiveness of activating the arm. For the original multi-armed bandit problem~\eqref{max:dsctStrict} and~\eqref{ys:dsctStrict}, starting from the initial belief states for all arms~$\Omega(1)$, we only need to calculate the (approximate) Whittle index for the state of each arm and select the~$M$ arms with the highest Whittle indices at each time. Since solving for the Whittle index as well as verifying condition~\eqref{eq: relaxedIdxIff} for each arm has a complexity determined by the process of solving the set of linear equations~\eqref{eq: linearET1},~\eqref{eq: linearET2} and~\eqref{eq:equal at threshold}, we have an efficient {\em online} algorithm for arm selections with a polynomial running time of the arm internal state size~$K$ and a linear running time of the number of arms~$N$ at each time~$t$. Specifically, our algorithm has complexity~$O(K^3NT)$, given that the first-crossing function~$L(\cdot,\omega)$ is solved for any~$\omega$ fixed as a point in the threshold.
  \item Recall the definition of~$L(\cdot,\omega)$ in~\eqref{def:firstCrossingT}. Since~$\Tau^k(\omega_1)=\omega_1\textbf{P}^k=\omega_1QJ^kQ^{-1}$ with~$\textbf{P}=QJQ^{-1}$ in its Jordan canonical form, it is possible to have an analytical solution for~$L(\cdot,\omega)$. In Sec~\ref{sec:k=3}, we focus on the case of~$K=3$ and obtain detailed forms of~$L(\cdot,\omega)$ in various scenarios. In general, one could use the exhaustion method to search for the first crossing time with an upper bound on the number of steps. If the number of search steps exceeds the upper bound, we simply set~$L(\cdot,\omega)=\infty$. As~$\beta^k$ decreases geometrically with~$k$ and~$\omega_1\textbf{P}^k$ with any belief state~$\omega_1$ also converges geometrically with~$k$ for regular transition matrices, such a numerical exhaustion has its practical convenience.
  \item The relaxed indexability is reduced to the classical one given the threshold function $V_{\beta,m}(\omega;u=1)-V_{\beta,m}(\omega;u=0)$. First, the classical indexability implies the unique existence of the minimum~$m$ for any belief state~$\omega$ at which the active and passive actions are indifferent\footnotemark\footnotetext{In case of an interval, we choose the minimum~$m$.}. On the other hand, the unique existence of~$m$ under the relaxed indexability implies the monotonic property of the passive set, for otherwise there will be two different values of $m$ that pushes some~$\omega$ to the passive set as~$m$ increases from $-\infty$ to $+\infty$ due to the continuity of the value functions. For $K=2$, the relaxed indexability relative to the {\em linearized} threshold function is equivalent to the classical indexability as the threshold consists of only one point.
\end{itemize}

Now we are ready to present the general algorithm based on the approximate Whittle index under the relaxed indexability relative to the linearized threshold function, as detailed in Algorithm~\ref{alg: whittle}.

\begin{algorithm}
	\caption{Whittle Index Policy}
	\hspace*{0.02in} {\bf Input:}
	time period $T$, arm number $N$, active arm number $M$\\
	\hspace*{0.02in} {\bf Input:}
	 initial belief state $\omega_n(1)$, transition matrix $\textbf{P}^{(n)}$, reward vector $B_n,~n=1,...,N$\\
    \hspace*{0.02in} {\bf Input:}
    discount factor $\beta$, first-crossing search maximum $l_{max}$
	\begin{algorithmic}[1]
	\For{$t=1,2,...,T$}
	       \For{$n=1,..., N$}
	       \State For each $i$, do a binary search for $L\left(\omega_n(t)\textbf{P}^{(n)}, \omega_n(t)\right)$ and $L\left(p^{(n)}_i, \omega_n(t)\right)$ in $\{1,\ldots,l_{max}\}$
           \State For any value of $L(\cdot,\cdot)$ not found in $\{1,\ldots,l_{max}\}$, set it to $+\infty$
		   \State Compute $f(p^{(n)}_i) = \frac{1 - \beta^{L\left(p^{(n)}_i, \omega_n(t)\right)}}{1 - \beta},\quad i=0,1,2,\cdots,K-1$
		   \State Set $F(\textbf{P}^{(n)}) = [f(p^{(n)}_0), f(p^{(n)}_1), f(p^{(n)}_2), \cdots, f(p^{(n)}_{K-1})]'$
		   \State Compute $g(p^{(n)}_i) = \beta^{L\left(p^{(n)}_i,\omega_n(t)\right)}\omega_n(t)(\textbf{P}^{(n)})^{L\left(p^{(n)}_i, \omega_n(t)\right)},\quad i=0,1,2,\cdots,K-1$
		   \State Set $G(\textbf{P}^{(n)}) = [g(p^{(n)}_0), g(p^{(n)}_1), g(p^{(n)}_2), \cdots, g(p^{(n)}_{K-1})]'$ and $H(\textbf{P}^{(n)}) = \left(\textbf{I}_K - \beta G(\textbf{P}^{(n)})\right)^{-1}$
		   \State Compute $f(\omega_n(t)\textbf{P}^{(n)}) = \frac{1   - \beta^{L\left(\omega_n(t)\textbf{P}^{(n)}, \omega_n(t)\right)}}{1 - \beta}$
           \State Compute $g(\omega_n(t)\textbf{P}^{(n)}) = \beta^{L\left(\omega_n(t)\textbf{P}^{(n)}, \omega_n(t)\right)}\omega_n(t)(\textbf{P}^{(n)})^{L\left(\omega_n(t)\textbf{P}^{(n)}, \omega_n(t)\right)+1}$
           \State Define $A\defeq1 + \beta f(\omega_n(t)\textbf{P}^{(n)}) + \beta\biggl[\beta g(\omega_n(t)\textbf{P}^{(n)}) - \omega_n(t)\biggr]H(\textbf{P}^{(n)})F(\textbf{P}^{(n)})$
           \State Set $W(\omega_n(t))=\omega_n(t) B_n'$ and skip Step 13 if $A=0$
		   \State Compute $W(\omega_n(t)) = \frac{\omega_n(t) B_n' - \beta g(\omega_n(t)\textbf{P}^{(n)})\biggl[\textbf{I}_K + \beta H(\textbf{P}^{(n)})G(\textbf{P}^{(n)})\biggr]B_n' + \beta\omega_n(t) H(\textbf{P}^{(n)})G(\textbf{P}^{(n)})B_n'}{1 + \beta f(\omega_n(t)\textbf{P}^{(n)}) + \beta\biggl[\beta g(\omega_n(t)\textbf{P}^{(n)}) - \omega_n(t)\biggr]H(\textbf{P}^{(n)})F(\textbf{P}^{(n)})}$
   \EndFor
   \State Choose the top $M$ arms with the largest Whittle indices $W(\omega_n(t))$
   \State Observe the states $S_n(t)$ of the selected $M$ arms and accrue the reward
   \For{$n=1,...,N$}
        \If{arm $n$ is active}
             \State $\omega_n(t+1) = \left[p^{(n)}_{S_n(t),0},p^{(n)}_{S_n(t),1},p^{(n)}_{S_n(t),2},\cdots,p^{(n)}_{S_n(t),K-1}\right]$
        \Else
             \State $\omega_n(t+1) = \omega_n(t)\textbf{P}^{(n)}$
        \EndIf
   \EndFor
\EndFor
\end{algorithmic}\label{alg: whittle}
\end{algorithm}

\section{Numerical Analysis}\label{sec:numeric}

In this section, we demonstrate the near-optimality of the (approximate) Whittle index policy by setting $M=1$ (each time one arm is played) and $K=3$ (the first-crossing time is more efficiently computed as in Sec.~\ref{sec:k=3}). 

\subsection{Comparison to Optimal and Myopic Policies}
Through extensive numerical examples, we compute the performance of the optimal policy by dynamic programming and simulate the low-complexity Whittle index policy and the myopic policy by Monte-Carlo runs. The performance of Whittle index policy has been observed to be quite close to the optimal one from all these numerical trials. Here we list a few examples in~\Cref{fig: result-1,fig: result-2,fig: result-3,fig: result-4,fig: result-5,fig: result-6} with their system parameters shown in~\Cref{table: experiment setting 1,table: experiment setting 2-1,table: experiment setting 2-2}. The comparison between the Whittle index policy and the myopic policy in Figures~\ref{fig: result-5} and~\ref{fig: result-6} demonstrates the superiority of the former.

To better illustrate the efficiency of the Whittle index policy, we further plot its performance versus the myopic policy for large number of arms and long time horizon. The optimal policy was not plotted due to the curse of dimensionality for large systems. As observed in~\Cref{fig: result-7,fig: result-8,fig: result-9,fig: result-10,fig: result-11,fig: result-12,fig: result-13,fig: result-14,fig: result-15,fig: result-16} with time horizon $T=100$ and 1000 monte-carlo runs for smoothing each curve, the Whittle index policy clearly shows a stronger performance for arm number $N=20,~30,~40,50~\mbox{and}~60$, respectively (two figures for each case). For larger systems or longer time horizons in consideration, the Whittle index policy becomes more significant since it has only a linear complexity with the number of arms (as well as with the length of time horizon), while solving for the optimal policy has an exponential complexity as the joint-state space grows geometrically with the number of arms. As times goes, the (approximate) Whittle index policy maintains a better balance between exploitation and exploration than the myopic one which only maximizes the immediate reward.

\subsection{More Performance Insights}
\begin{figure}
	\begin{minipage}[t]{0.5\linewidth}
		\centering
		\includegraphics[height=5cm,width=7cm]{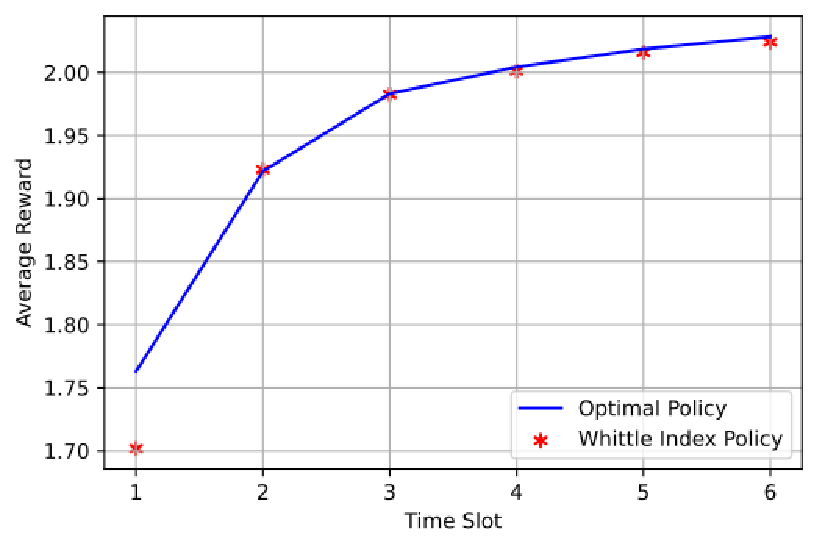}
		\caption{Experiment 1: machine 1}
		\label{fig: result-1}
	\end{minipage}%
	\begin{minipage}[t]{0.5\linewidth}
		\centering
		\includegraphics[height=5cm,width=7cm]{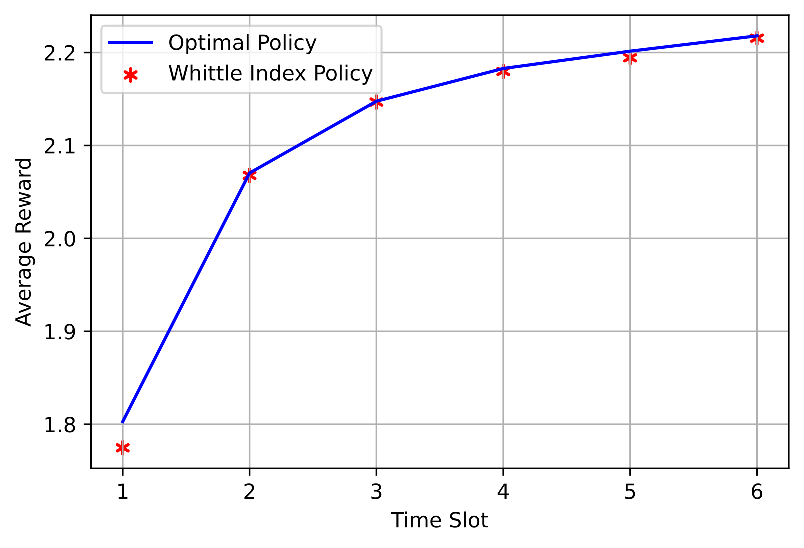}
		\caption{Experiment 1: machine 2}
		\label{fig: result-2}
	\end{minipage}
\end{figure}

\begin{figure}
	\begin{minipage}[t]{0.5\linewidth}
		\centering
		\includegraphics[height=5cm,width=7cm]{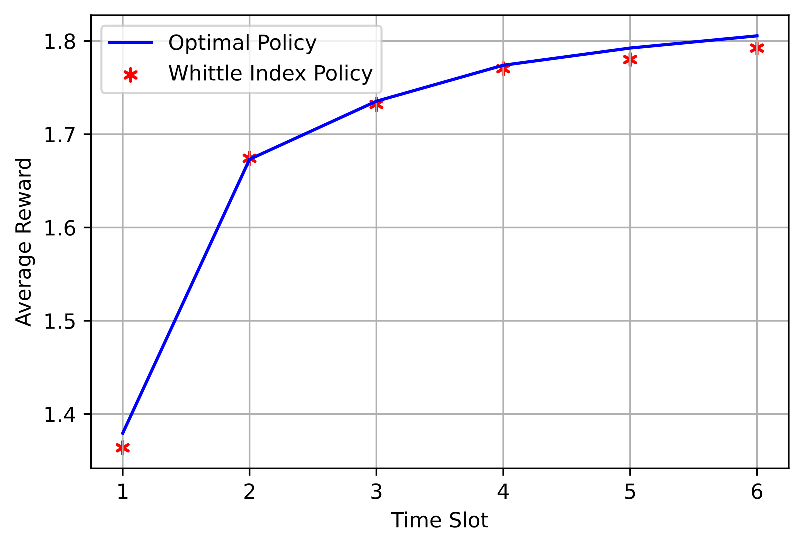}
		\caption{Experiment 2: machine 1}
		\label{fig: result-3}
	\end{minipage}%
	\begin{minipage}[t]{0.5\linewidth}
		\centering
		\includegraphics[height=5cm,width=7cm]{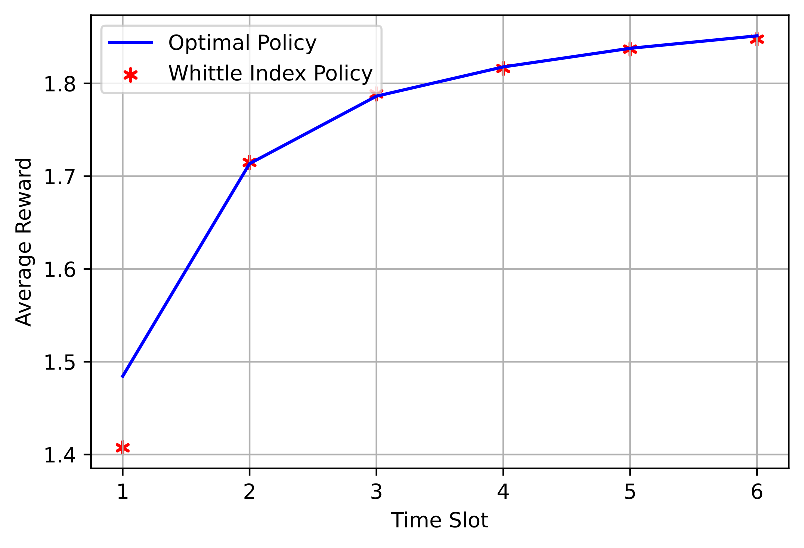}
		\caption{Experiment 2: machine 2}
		\label{fig: result-4}
	\end{minipage}
\end{figure}

\begin{figure}
	\begin{minipage}[t]{0.5\linewidth}
		\centering
		\includegraphics[height=5cm,width=7cm]{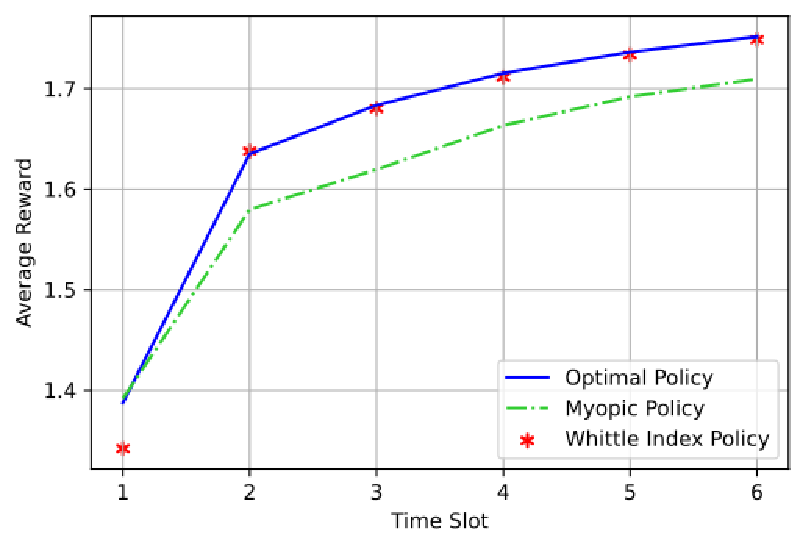}
		\caption{Experiment 2: machine 3}
		\label{fig: result-5}
	\end{minipage}%
    \begin{minipage}[t]{0.5\linewidth}
		\centering
		\includegraphics[height=5cm,width=7cm]{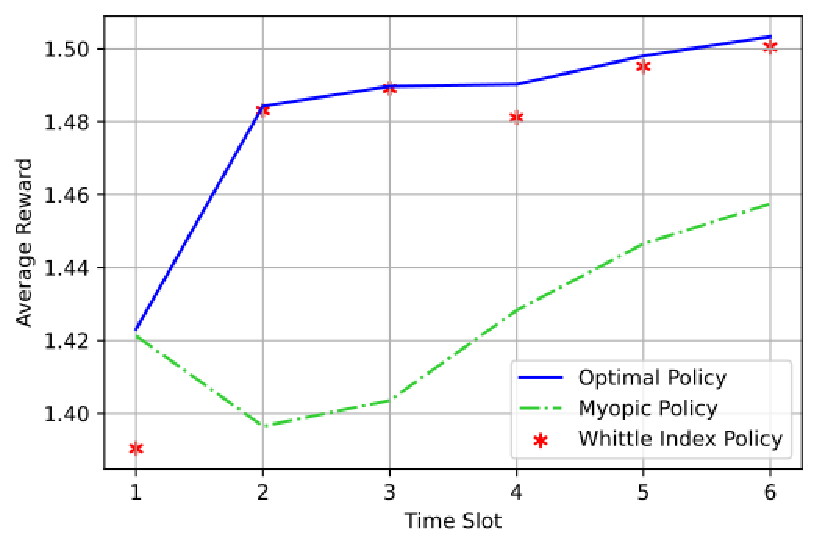}
		\caption{Experiment 2: machine 4}
		\label{fig: result-6}
	\end{minipage}
\end{figure}

If we take a closer look at \Cref{fig: result-1,fig: result-2,fig: result-3,fig: result-4,fig: result-5,fig: result-6}, there is an interesting phenomenon that the Whittle index policy can be a little worse than the optimal policy in the middle but eventually catches up as time goes. This illustrates the asymptotic optimality of Whittle index studied by \cite{WW1990}. In other words, Whittle's relaxation from \eqref{ys:dsctStrict} to \eqref{ys:dsctRelax} should not fundamentally change the optimal long-term system state and action paths from the perspective of large deviation theory \citep{WW1990}. Therefore, applying the Whittle index derived from the relaxed problem to the original one (exactly $M$ arms are played at each time) is effective, especially when time goes to infinity. We emphasize on the observation that the near-optimality achieved here is from the (approximate) Whittle index derived from the relaxed indexability. This further demonstrates that the relaxed indexability does not suffer significant loss in terms of utilizing the arm state dynamics compared to the original indexability, but can be much more efficiently implemented for the high-dimensional case ($K>2$).

\begin{figure}
	\begin{minipage}[t]{0.5\linewidth}
		\centering
		\includegraphics[height=5cm,width=7cm]{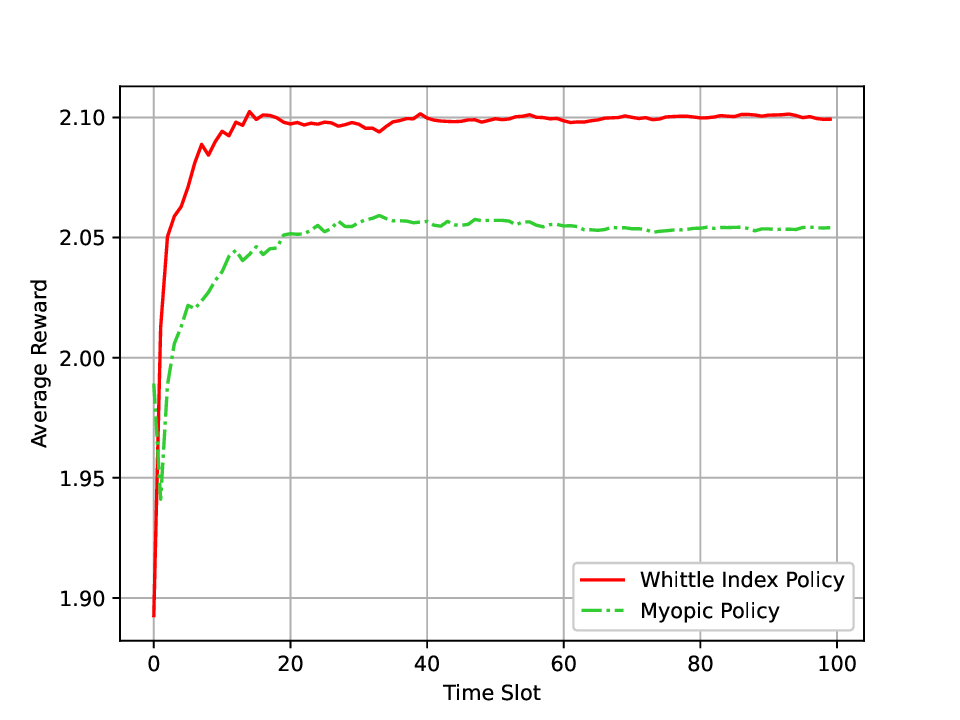}
		\caption{Large System 1: $N=20$}
		\label{fig: result-7}
	\end{minipage}%
    \begin{minipage}[t]{0.5\linewidth}
		\centering
		\includegraphics[height=5cm,width=7cm]{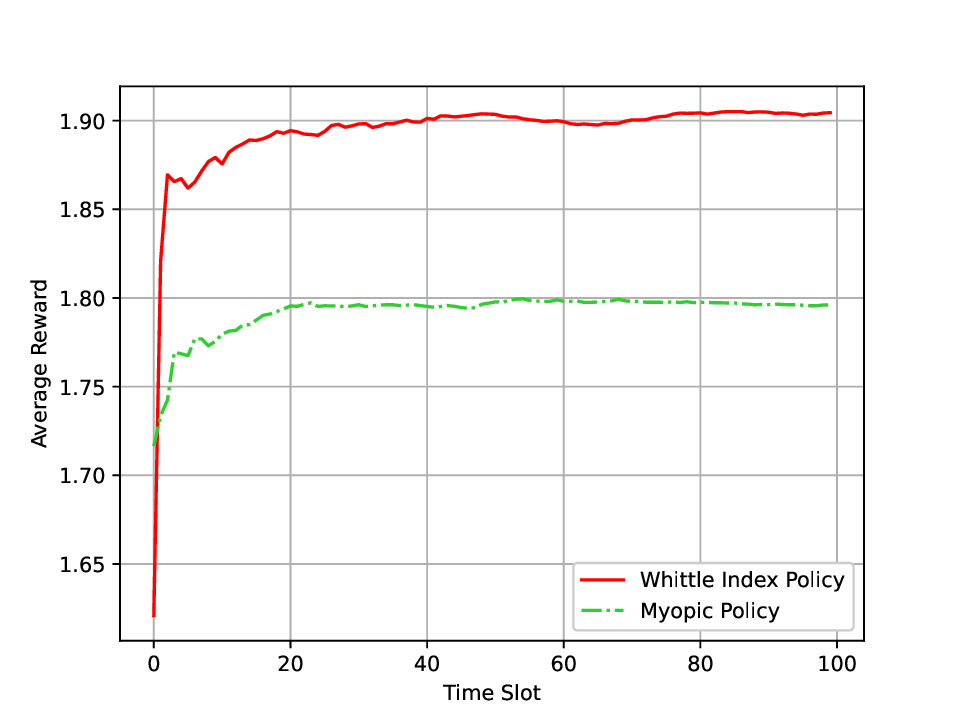}
		\caption{Large System 2: $N=20$}
		\label{fig: result-8}
	\end{minipage}
\end{figure}

\begin{figure}
	\begin{minipage}[t]{0.5\linewidth}
		\centering
		\includegraphics[height=5cm,width=7cm]{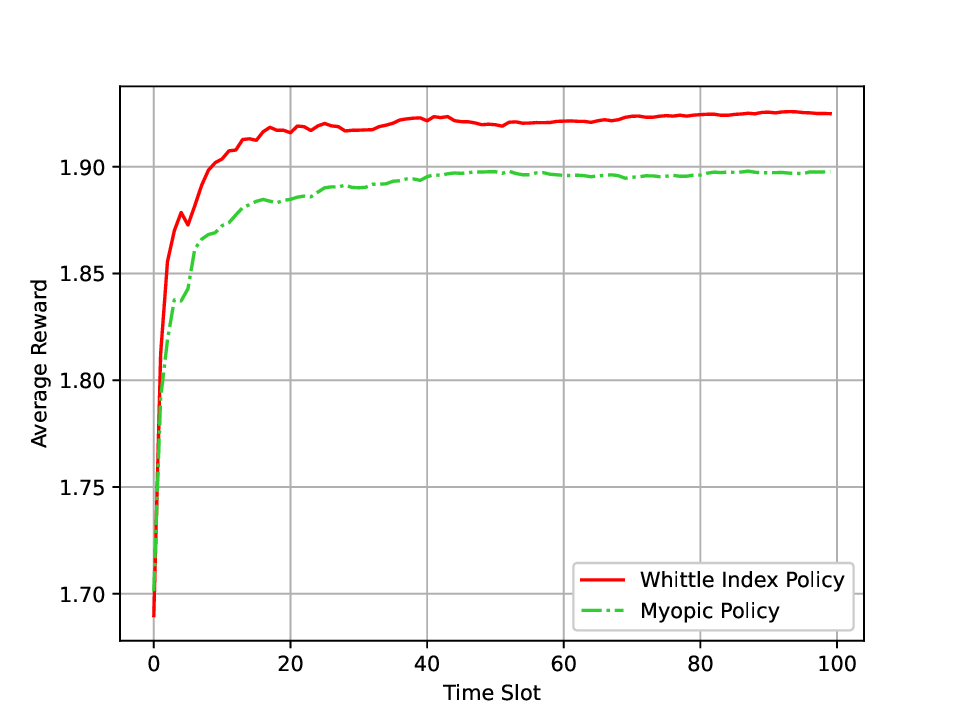}
		\caption{Large System 3: $N=30$}
		\label{fig: result-9}
	\end{minipage}%
    \begin{minipage}[t]{0.5\linewidth}
		\centering
		\includegraphics[height=5cm,width=7cm]{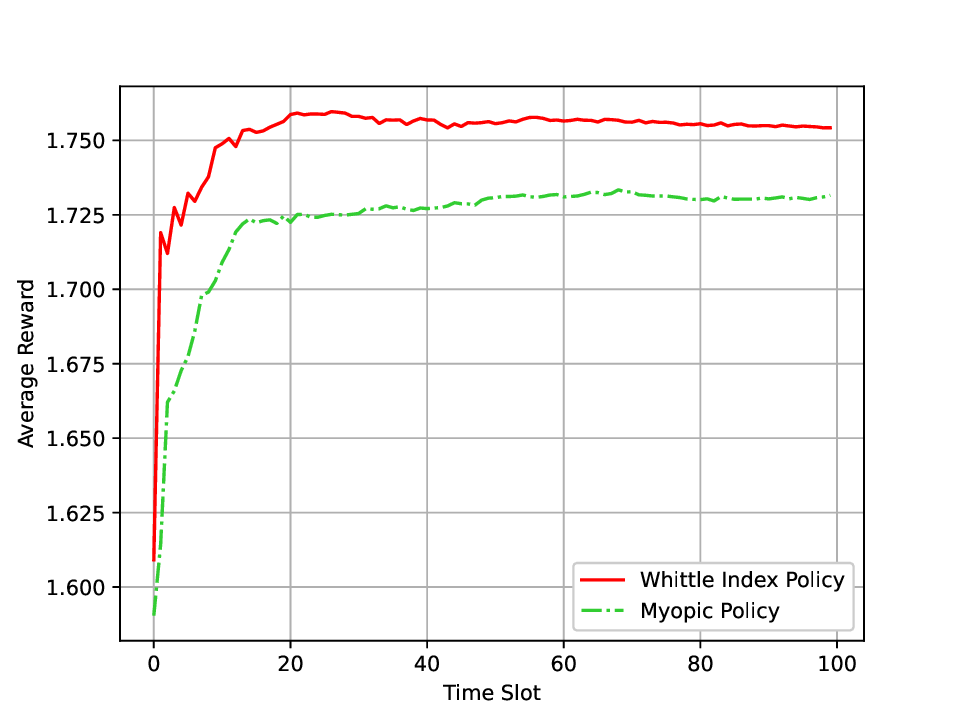}
		\caption{Large System 4: $N=30$}
		\label{fig: result-10}
	\end{minipage}
\end{figure}

\begin{figure}
	\begin{minipage}[t]{0.5\linewidth}
		\centering
		\includegraphics[height=5cm,width=7cm]{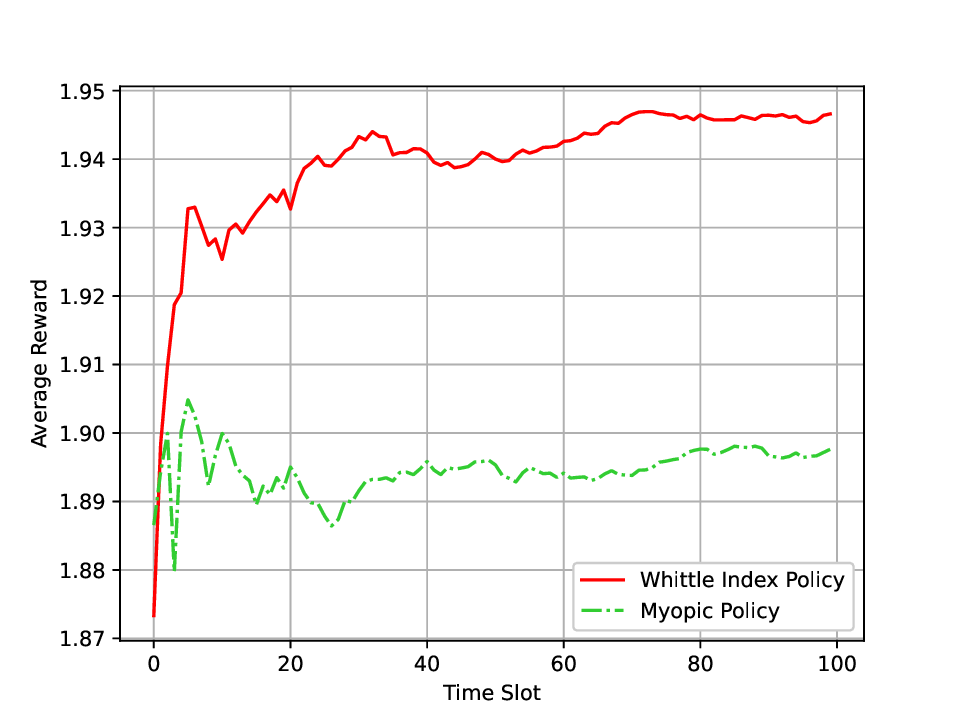}
		\caption{Large System 5: $N=40$}
		\label{fig: result-11}
	\end{minipage}%
    \begin{minipage}[t]{0.5\linewidth}
		\centering
		\includegraphics[height=5cm,width=7cm]{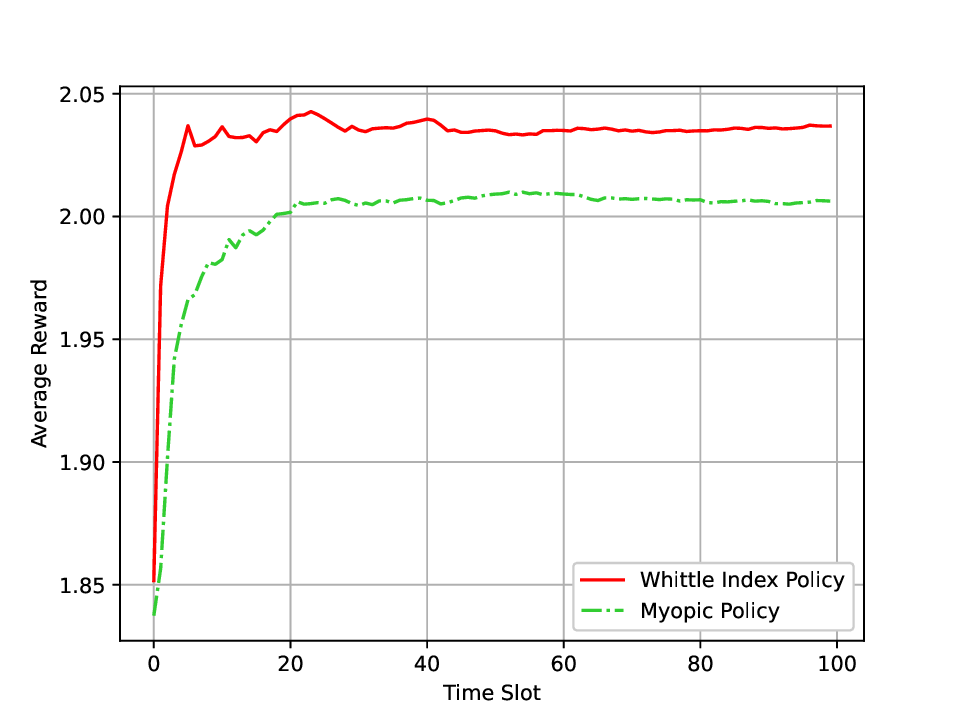}
		\caption{Large System 6: $N=40$}
		\label{fig: result-12}
	\end{minipage}
\end{figure}

\begin{figure}
	\begin{minipage}[t]{0.5\linewidth}
		\centering
		\includegraphics[height=5cm,width=7cm]{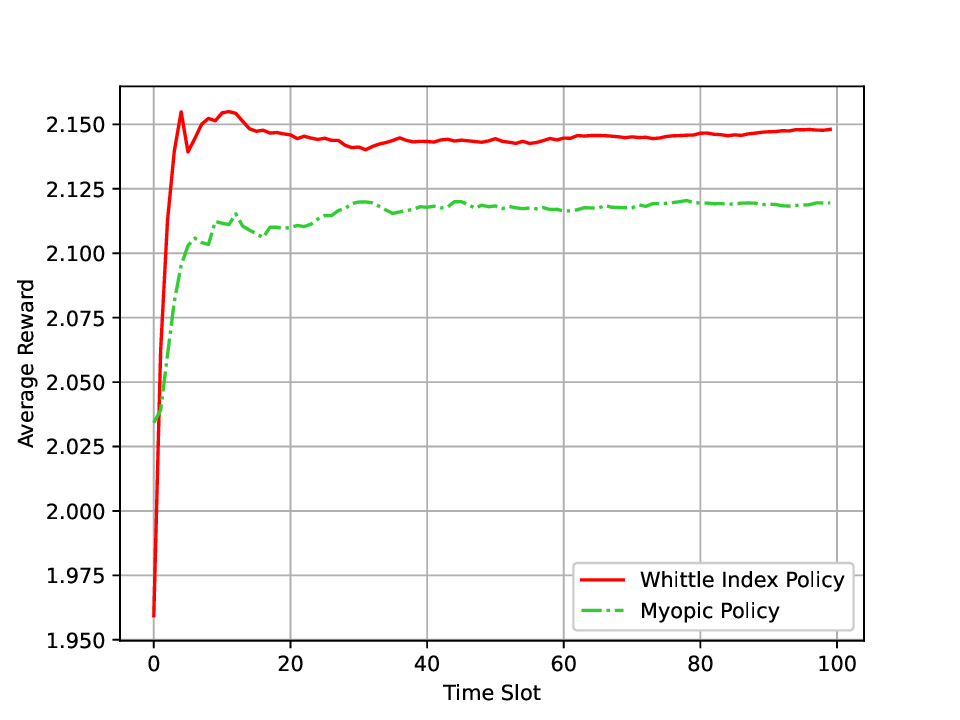}
		\caption{Large System 7: $N=50$}
		\label{fig: result-13}
	\end{minipage}%
    \begin{minipage}[t]{0.5\linewidth}
		\centering
		\includegraphics[height=5cm,width=7cm]{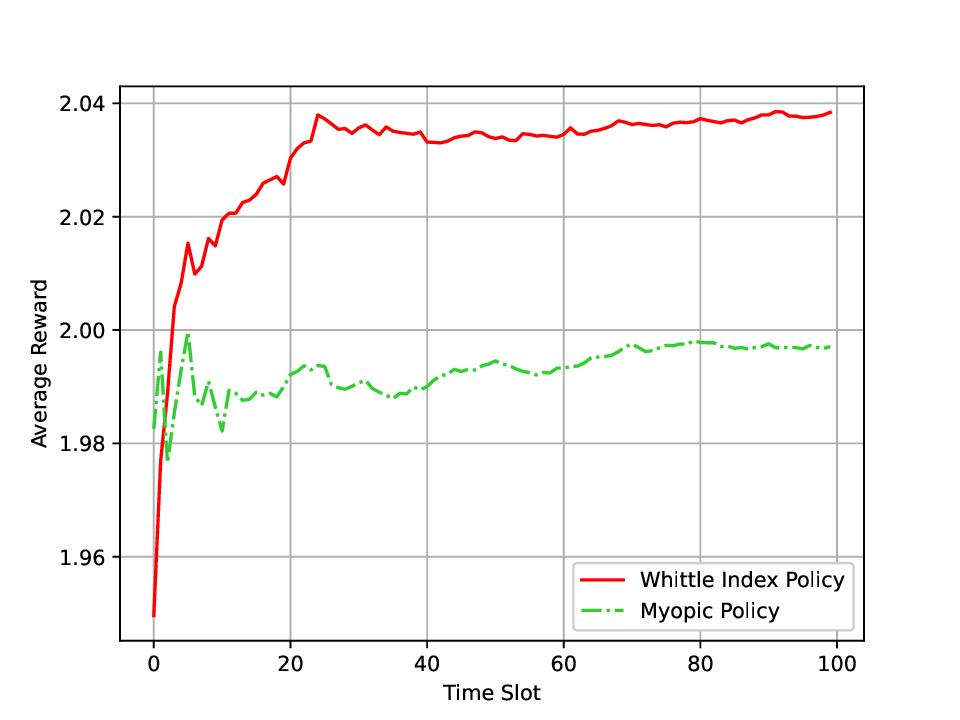}
		\caption{Large System 8: $N=50$}
		\label{fig: result-14}
	\end{minipage}
\end{figure}

\begin{figure}
	\begin{minipage}[t]{0.5\linewidth}
		\centering
		\includegraphics[height=5cm,width=7cm]{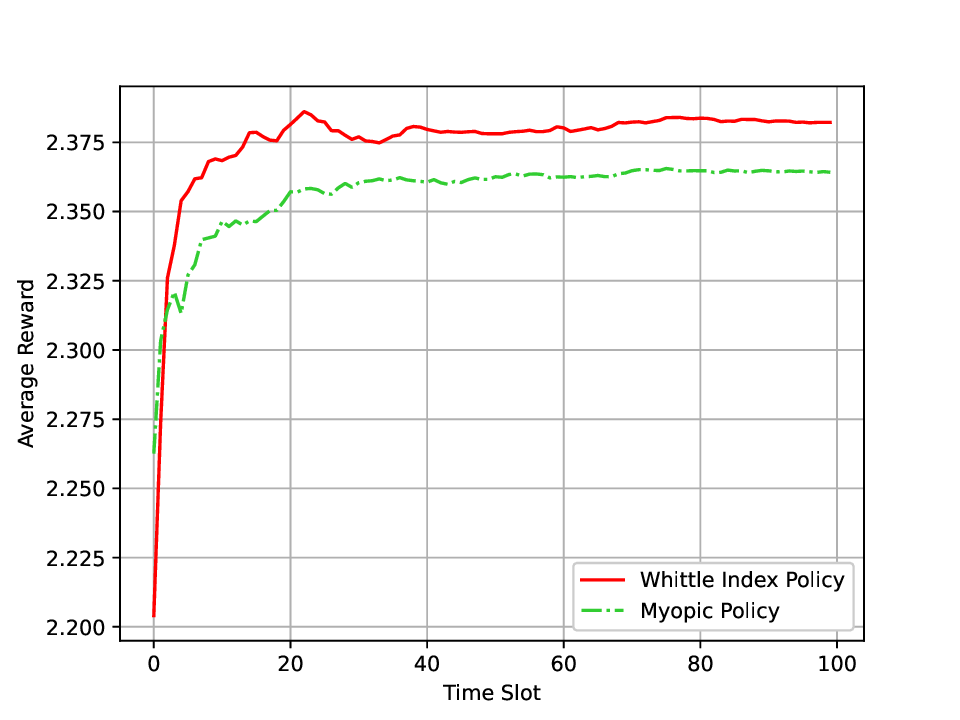}
		\caption{Large System 7: $N=60$}
		\label{fig: result-15}
	\end{minipage}%
    \begin{minipage}[t]{0.5\linewidth}
		\centering
		\includegraphics[height=5cm,width=7cm]{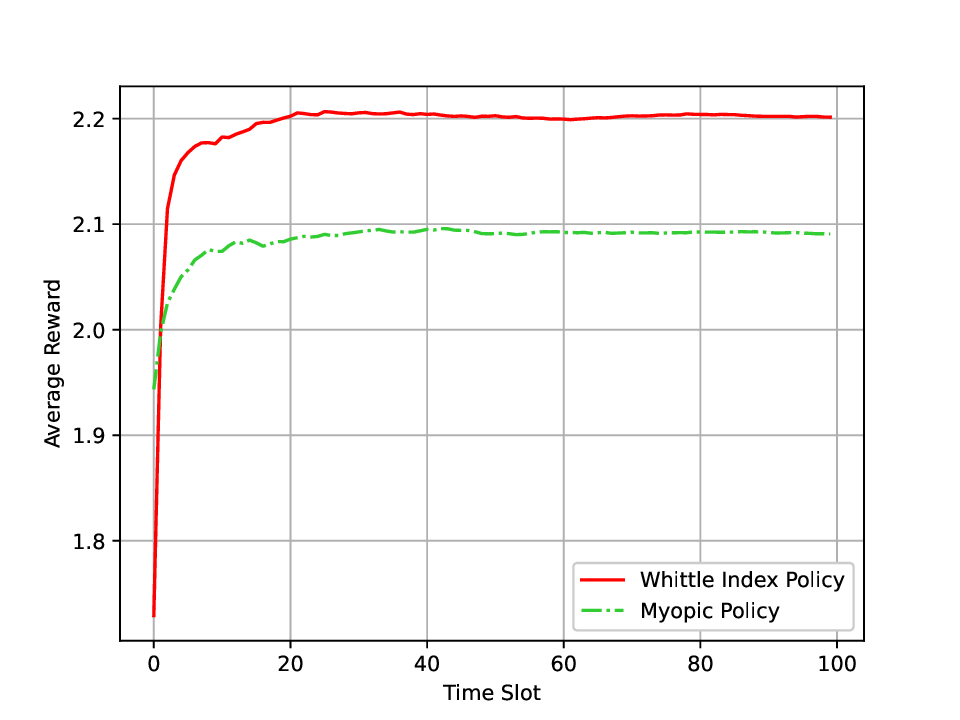}
		\caption{Large System 8: $N=60$}
		\label{fig: result-16}
	\end{minipage}
\end{figure}

There are also some interesting observations on the performance of Whittle index under the relaxed indexability when compared to the myopic policy (\Cref{fig: result-5,fig: result-6,fig: result-7,fig: result-8,fig: result-9,fig: result-10,fig: result-11,fig: result-12,fig: result-13,fig: result-14,fig: result-15,fig: result-16}). First of all, the curve of the myopic policy always converges to a point lower than the limit of the curve of the Whittle index policy as time goes to infinity: the myopic policy can never catch up! This is in accordance with the nature of the myopic policy as it does not consider the long-term benefit. In contrast, the Whittle index policy leads the system state and action paths to follow the optimal ones by evaluating the priorities of different arms based on their long-term performance (value functions). Second, the myopic policy may get the player trapped in a set of arms and not able to sufficiently explore other arms which may yield higher reward in the future. This is because some arms yield better immediate rewards and keep being played for the majority of time but other arms scarcely get a chance to be played and their states never transit to more advantageous positions.

\section{Conclusion and Future Work}
In this paper, we proposed an efficient algorithm to achieve strong performance for a class of restless multi-armed bandits arising within in the general POMDP framework. By formulating the problem with a $(K-1)$-dimensional belief state space, we extended the Whittle index policy previously studied for the case of $K=2$ to $K>2$ by introducing the concept of relaxed indexability. An interesting finding is that through the online computation process for the first crossing time, all our numerical studies have shown that the relaxed indexability relative to the linearized threshold function was satisfied. Future work includes the extensions of Whittle index to more general POMDP models, \eg those with observation errors or different state transition dynamics. Furthermore, the approximation of the decision boundary can be readily implemented by the classical $t$-step lookahead approach in dynamic programming with $t$ chosen to control the the tradeoff between the approximation accuracy and the time complexity \citep{Bert1987}. Specifically, the $t$-step threshold function is defined as
\begin{eqnarray}
r_t(\omega)\defeq V_{\beta,m,[1,\cdots,t]}(\omega;u=1)-V_{\beta,m,[1,\cdots,t]}(\omega;u=0),
\end{eqnarray}
where $V_{\beta,m,[1,\cdots,t]}(\cdot)$ denotes the maximum expected reward obtained over~$t$ steps. Obviously our linearized threshold function is equivalent to the case of $t=1$. The larger $t$ is, the better approximation to the optimal threshold function since
\begin{eqnarray}
r_t(\cdot)\rightarrow V_{\beta,m}(\omega;u=1)-V_{\beta,m}(\omega;u=0), \ \mbox{as}~t\rightarrow\infty.
\end{eqnarray}
So the approximation error to the original decision boundary $C(m)$ converges to zero as $t\rightarrow\infty$. Nevertheless, the algorithmic complexity definitely increases with~$t$.

Future work will include the theoretical verification of relaxed indexability in relation to the linearized threshold function (\ie $t=1$), systematic categorization of threshold functions to develop efficient algorithms for general partially observable RMAB problems, and analysis of the duality gap introduced by the relaxation for arm-decoupling (specifically, the relaxation concerning the expected number of arms to activate).

\clearpage

\begin{table}[h]\notsotiny{
	\caption{Experiment setting 1 ($\beta = 0.9999, B_i=[0, 2, 3], i=1,...,7$)}
    \label{table: experiment setting 1}
    \centering
	\begin{tabular}{|c|c|c|}
		\hline
		\diagbox{arm}{machine}&1&2\\
		\hline
		1 & $\textbf{P}^{(1)}=\left(
		\begin{matrix}
			0.514&0.321&0.165\\0.123&0.159&0.718\\0.420&0.195&0.385
		\end{matrix}
		\right)$ &
		$\textbf{P}^{(1)}=\left(
		\begin{matrix}
			0.519&0.445&0.036\\0.188&0.292&0.520\\0.043&0.292&0.665
		\end{matrix}
		\right)$\\
		& $\omega_1(1)=$(0.279, 0.618, 0.103)&
		$\omega_1(1)=$(0.354, 0.164, 0.482)\\
		\hline
		2 & $\textbf{P}^{(2)} = \left(
		\begin{matrix}
			0.372&0.543&0.085\\0.103&0.633&0.264\\0.417&0.301&0.282
		\end{matrix}
		\right)$ &
		$\textbf{P}^{(2)}=\left(
		\begin{matrix}
			0.193&0.534&0.273\\0.275&0.485&0.240\\0.234&0.694&0.072
		\end{matrix}
		\right)$ \\
		& $\omega_2(1)=$(0.688, 0.024, 0.288)&
		$\omega_2(1)=$(0.426, 0.188, 0.386)\\
		\hline
		3& $\textbf{P}^{(3)} = \left(
		\begin{matrix}
			0.405&0.129&0.466\\0.413&0.328&0.259\\0.327&0.502&0.171
		\end{matrix}
		\right)$ &
		$\textbf{P}^{(3)}=\left(
		\begin{matrix}
			0.250&0.274&0.476\\0.600&0.242&0.158\\0.271&0.199&0.530
		\end{matrix}
		\right)$  \\
		& $\omega_3(1)=$(0.489, 0.408, 0.103)&
		$\omega_3(1)=$(0.333, 0.498, 0.169)\\
		\hline
		4& $\textbf{P}^{(4)} = \left(
		\begin{matrix}
			0.461&0.272&0.267\\0.555&0.431&0.014\\0.058&0.689&0.253
		\end{matrix}
		\right)$&
		$\textbf{P}^{(4)}=\left(
		\begin{matrix}
			0.721&0.203&0.076\\0.201&0.621&0.178\\0.444&0.319&0.237
		\end{matrix}
		\right)$ \\
		& $\omega_4(1)=$(0.554, 0.061 , 0.385)&
		$\omega_4(1)=$(0.455, 0.285, 0.260) \\
		\hline
		5& $\textbf{P}^{(5)} = \left(
		\begin{matrix}
			0.339&0.427&0.234\\0.161&0.469&0.370\\0.265&0.296&0.439
		\end{matrix}
		\right)$&
		$\textbf{P}^{(5)}=\left(
		\begin{matrix}
			0.161&0.445&0.394\\0.249&0.394&0.357\\0.164&0.363&0.473
		\end{matrix}
		\right)$ \\
		& $\omega_5(1)=$(0.313, 0.297, 0.390)&
		$\omega_5(1)=$(0.352, 0.424, 0.224)\\
		\hline
		6& $\textbf{P}^{(6)} = \left(
		\begin{matrix}
			0.071&0.556&0.373\\0.158&0.308&0.534\\0.412&0.459&0.129
		\end{matrix}
		\right)$&
		$\textbf{P}^{(6)}=\left(
		\begin{matrix}
			0.080&0.279&0.641\\0.027&0.780&0.193\\0.418&0.265&0.317
		\end{matrix}
		\right)$ \\
		& $\omega_6(1)=$(0.332, 0.305, 0.363)&
		$\omega_6(1)=$(0.102, 0.893, 0.005)\\
		\hline
		7& $\textbf{P}^{(7)} = \left(
		\begin{matrix}
			0.427&0.324&0.249\\0.478&0.356&0.166\\0.326&0.490&0.184
		\end{matrix}
		\right)$&
		$\textbf{P}^{(7)}=\left(
		\begin{matrix}
			0.130&0.536&0.334\\0.377&.253&0.370\\0.334&0.120&0.546
		\end{matrix}
		\right)$ \\
		& $\omega_7(1)=$(0.234, 0.722, 0.044)&
		$\omega_7(1)=$(0.367, 0.276, 0.357)\\
		\hline
	\end{tabular}}
\end{table}

\begin{table}\notsotiny{
	\caption{ Experiment setting 2-1 ($\beta = 0.9999$)}
    \label{table: experiment setting 2-1}
	\begin{tabular}{|c|c|c|}
		\hline
		\diagbox{arm}{machine}&1&2\\
		\hline
		1 & $\textbf{P}^{(1)}=\left(
		\begin{matrix}
			0.036&0.607&0.357\\0.053&0.126&0.821\\0.579&0.359&0.062
		\end{matrix}
		\right)$ &
		$\textbf{P}^{(1)}=\left(
		\begin{matrix}
			0.538&0.305&0.157\\0.575&0.097&0.328\\0.346&0.168&0.486
		\end{matrix}
		\right)$\\
		& $\omega_1(1)=$(0.284, 0.404, 0.312)&$\omega_1(1)=$(0.462, 0.418, 0.120)\\
		&$B_1 = $(0, 1.004, 2.186)&$B_1 = $(0, 2.422, 2.698)\\
		\hline
		2 & $\textbf{P}^{(2)} = \left(
		\begin{matrix}
			0.447&0.021&0.532\\0.485&0.164&0.351\\0.461&0.409&0.130
		\end{matrix}
		\right)$ &
		$\textbf{P}^{(2)}=\left(
		\begin{matrix}
			0.367&0.114&0.519\\0.798&0.089&0.113\\0.367&0.354&0.279
		\end{matrix}
		\right)$ \\
		& $\omega_2(1)=$(0.297, 0.361, 0.342)&$\omega_2(1)=$(0.459, 0.528, 0.013)\\
		&$B_2 = $(0, 1.155, 2.761)&$B_2 = $(0, 2.745, 2.754)\\
		\hline
		3& $\textbf{P}^{(3)} = \left(
		\begin{matrix}
			0.407&0.200&0.393\\0.435&0.180&0.385\\0.245&0.465&0.290
		\end{matrix}
		\right)$ &
		$\textbf{P}^{(3)}=\left(
		\begin{matrix}
			0.218&.015&0.767\\0.849&0.129&0.022\\0.405&0.151& 0.444
		\end{matrix}
		\right)$ \\
		& $\omega_3(1)=$(0.043, 0.421, 0.536)&$\omega_3(1)=$(0.519, 0.413, 0.068)\\
		&$B_3 = $(0, 0.437, 0.7826)&$B_3 = $(0, 2.917, 2.916)\\
		\hline
		4& $\textbf{P}^{(4)} = \left(
		\begin{matrix}
			0.087&0.454&0.459\\0.181&0.672&0.147\\0.462&0.492&0.046
		\end{matrix}
		\right)$&
		$\textbf{P}^{(4)}=\left(
		\begin{matrix}
			0.428&0.294&0.278\\0.431&0.022&0.547\\0.011&0.511& 0.478
		\end{matrix}
		\right)$  \\
		& $\omega_4(1)=$(0.642, 0.026, 0.332)& $\omega_4(1)=$(0.113, 0.499, 0.388) \\
		&$B_4 = $(0, 0.568, 0.619)&$B_4 = $(0, 0.051, 0.503)\\
		\hline
		5& $\textbf{P}^{(5)} = \left(
		\begin{matrix}
			0.331&0.181&0.488\\0.347&0.117&0.536\\0.245&0.197&0.558
		\end{matrix}
		\right)$&
		$\textbf{P}^{(5)}=\left(
		\begin{matrix}
			0.317&0.413&0.270\\0.376&0.333&0.291\\0.351&0.203& 0.446
		\end{matrix}
		\right)$  \\
		& $\omega_5(1)=$(0.606, 0.017, 0.377)& $\omega_5(1)=$(0.555, 0.400, 0.045)\\
		&$B_5 = $(0, 2.448, 2.63 )&$B_5 = $(0, 1.51 , 2.688)\\
		\hline
		6& $\textbf{P}^{(6)} = \left(
		\begin{matrix}
			0.449&0.008&0.543\\0.782&0.198&0.020\\0.338&0.614&0.048
		\end{matrix}
		\right)$&
		$\textbf{P}^{(6)}=\left(
		\begin{matrix}
			0.320&0.649&0.031\\0.112&0.037&0.851\\0.377&0.364& 0.259
		\end{matrix}
		\right)$\\
		& $\omega_6(1)=$(0.362, 0.560, 0.078)&$\omega_6(1)=$(0.348, 0.212, 0.440) \\
		&$B_6 = $(0, 1.327, 1.945)&$B_6 = $(0, 1.623, 1.777)\\
		\hline
		7& $\textbf{P}^{(7)} = \left(
		\begin{matrix}
			0.067&0.435&0.498\\0.334&0.290&0.376\\0.258&0.483&0.259
		\end{matrix}
		\right)$&
		$\textbf{P}^{(7)}=\left(
		\begin{matrix}
			0.046&0.213&0.741\\0.031&0.430&0.539\\0.238&0.238& 0.524
		\end{matrix}
		\right)$ \\
		& $\omega_7(1)=$(0.296, 0.298, 0.406)& $\omega_7(1)=$(0.483, 0.050, 0.467)\\
		&$B_7 = $(0, 1.858, 2.033)&$B_7 = $(0, 0.897, 2.443)\\
		\hline
	\end{tabular}}
\end{table}

\begin{table}\notsotiny{
	\caption{Experiment setting 2-2 ($\beta = 0.9999$)}
    \label{table: experiment setting 2-2}
	\begin{tabular}{|c|c|c|}
		\hline
		\diagbox{arm}{machine}&3&4\\
		\hline
		1 & $\textbf{P}^{(1)}=\left(
		\begin{matrix}
			0.488&0.258&0.254\\0.012&0.790&0.198\\0.681&0.208&0.111
		\end{matrix}
		\right)$ &
		$\textbf{P}^{(1)}=\left(
		\begin{matrix}
			0.413&0.329&0.258\\0.089&0.511&0.400\\0.086&0.309&0.605
		\end{matrix}
		\right)$\\
		&$\omega_1(1)=$(0.405,0.415,0.180)&$\omega_1(1)=$(0.486, 0.028, 0.486)\\
		&$B_1 = $(0, 2.146, 2.491)&$B_1 = $(0, 0.233, 2.853)\\
		\hline
		2 & $\textbf{P}^{(2)} = \left(
		\begin{matrix}
			0.354&0.311&0.335\\0.278&0.027&0.695\\0.502&0.341& 0.157
		\end{matrix}
		\right)$ &
		$\textbf{P}^{(2)}=\left(
		\begin{matrix}
			0.031&0.171&0.798\\0.678&0.134&0.188\\0.597&0.358&0.045
		\end{matrix}
		\right)$ \\
		&$\omega_2(1)=$(0.551,0.328,0.121)&$\omega_2(1)=$(0.408, 0.496, 0.096)\\
		&$B_2 = $(0, 1.579, 2.444)&$B_2 = $(0, 2.358, 2.632)\\
		\hline
		3& $\textbf{P}^{(3)} = \left(
		\begin{matrix}
			0.342&0.036&0.622\\0.451&0.219&0.330\\0.471&0.073& 0.456
		\end{matrix}
		\right)$ &
		$\textbf{P}^{(3)}=\left(
		\begin{matrix}
			0.358&0.263&0.379\\0.264&0.249&0.487\\0.400&0.364&0.236
		\end{matrix}
		\right)$  \\
		&$\omega_3(1)=$(0.555,0.315,0.130)&$\omega_3(1)=$(0.014, 0.247, 0.739)\\
		&$B_3 = $(0, 0.286, 0.644)&$B_3 = $(0, 0.378, 1.241)\\
		\hline
		4& $\textbf{P}^{(4)} = \left(
		\begin{matrix}
			0.304&0.639&0.057\\0.457&0.380&0.163\\0.262&0.357& 0.381
		\end{matrix}
		\right)$&
		$\textbf{P}^{(4)}=\left(
		\begin{matrix}
			0.598&0.028&0.374\\0.762&0.109&0.129\\0.313&0.391&0.296
		\end{matrix}
		\right)$\\
		&$\omega_4(1)=$(0.495,0.117,0.388)& $\omega_4(1)=$(0.490, 0.256, 0.254)\\
		&$B_4 = $(0, 2.391, 2.852)&$B_4 = $(0, 2.002, 2.374)\\
		\hline
		5& $\textbf{P}^{(5)} = \left(
		\begin{matrix}
			0.404&0.282&0.314\\0.621&0.106&0.273\\0.204&0.657& 0.14
		\end{matrix}
		\right)$&
		$\textbf{P}^{(5)}=\left(
		\begin{matrix}
			0.323&0.177&0.500\\0.174&0.138&0.688\\0.416&0.310&0.274
			
		\end{matrix}
		\right)$\\
		&$\omega_5(1)=$(0.474,0.239,0.287)&$\omega_5(1)=$(0.358, 0.501, 0.141)\\
		&$B_5 = $(0, 0.111, 1.420)&$B_5 = $(0, 1.502, 2.258)\\
		\hline
		6& $\textbf{P}^{(6)} = \left(
		\begin{matrix}
			0.586&0.024&0.390\\0.455&0.027&0.518\\0.365&0.464&0.171
		\end{matrix}
		\right)$&
		$\textbf{P}^{(6)}=\left(
		\begin{matrix}
			0.424&0.442&0.134\\0.301&0.182&0.517\\0.164&0.360&0.476
		\end{matrix}
		\right)$ \\
		&$\omega_6(1)=$(0.413,0.388,0.199)&$\omega_6(1)=$(0.263, 0.502, 0.235)\\
		&$B_6 = $(0, 0.324, 0.755)&$B_6 = $(0, 0.715, 1.022)\\
		\hline
		7& $\textbf{P}^{(7)} = \left(
		\begin{matrix}
			0.612&0.335&0.053\\0.333&0.486&0.181\\0.483&0.513& 0.004
		\end{matrix}
		\right)$&
		$\textbf{P}^{(7)}=\left(
		\begin{matrix}
			0.613&0.136&0.251\\0.454&0.383&0.163\\0.287&0.693&0.020
		\end{matrix}
		\right)$ \\
		&$\omega_7(1)=$(0.369,0.262,0.369)&$\omega_7(1)=$(0.707, 0.226, 0.067)\\
		&$B_7 = $(0, 0.491, 0.797)&$B_7 = $(0, 2.013, 2.436)\\
		\hline
	\end{tabular}}
\end{table}

\clearpage

\ACKNOWLEDGMENT{The author gratefully acknowledges the help from his students, Jiale Zha and Chengzhong Zhang, for the numerical analysis and figures. The author's colleague Prof. Ting Wu and the anonymous reviewers provided very helpful comments for improving this paper.}

%
%
%






\ECSwitch

\noindent {\bf Title:} Relaxed Indexability and Index Policy for Partially Observable Restless Bandits

\noindent {\bf Author:} Keqin Liu

\vspace{1em}


\ECHead{Proofs of Lemmas and Theorems}

\section{Proof of Lemma~\ref{lm:vf}.}\label{proof:lmaCC}
\proof{} Consider a horizon of~$T~(T\ge1)$ time slots and define~$V_{\beta,m,1}(\omega)$ as the maximum expected total discounted reward over~$T$ slots that can be obtained starting from initial state~$\omega$ at~$t=1$:
\begin{eqnarray}
V_{\beta,m,1}(\omega)=\max_{\pi\in\Pi_{sa}(T)}&\mathbb{E}_{\pi}[\sum_{t=1}^{T}\beta^{t-1}\sum_{n=1}^N\{B_{n,S_n(t)}\mathbbm{1}(u(t)=1)+m\mathbbm{1}(u(t)=0)\}],\label{max:finiteT}
\end{eqnarray}
where~$\Pi_{sa}(T)$ is the set of single-arm policies that map the belief state~$\omega(t)$ to the action $u(t)\in\{1~\mbox{(active)},~0~\mbox{(passive)}\}$ for~$t=1,2,\cdots,T$. Note that~$\omega(1)=\omega$ and an optimal policy~$\pi^*_{sa}(T)$ achieving~$V_{\beta,m,T}(\omega)$ is generally non-stationary, \ie the mapping from~$\omega(t)$ to~$u(t)$ is dependent on~$t$. Especially when~$t=T$, we have only one more step to go and the myopic policy that maximizes the immediate reward is obviously optimal:
\begin{eqnarray}
u^*(T) = \arg\max_{u\in\{0,1\}}\{u\cdot\omega(T) B'+(1-u)\cdot m\}.
\end{eqnarray}
Let~$V_{\beta,m,t}(\cdot)$ denote the maximum expected total discounted reward accumulated from slot~$t$ to~$T$ under~$\pi^*_{sa}(T)$. We have the following dynamic equations:
\begin{eqnarray}
V_{\beta,m,t}(\omega(t)) &=& \max\left\{\omega B' + \beta \omega
	\left(
	\begin{matrix}
		V_{\beta,m,t+1}(p_0)\\
        V_{\beta,m,t+1}(p_1)\\
        \vdots\\
        V_{\beta,m,t+1}(p_{K-1})
	\end{matrix}
	\right), ~ m+V_{\beta,m,t+1}(\Tau^1(\omega(t)))\right\},\ t\le T,\\
V_{\beta,m,T+1}(\cdot) &\equiv& 0.
\end{eqnarray}
We first prove the properties of~$V_{\beta,m}(\omega)$ regarding to~$\omega$ with~$m$ fixed. Our approach is based on backward induction on~$t$ with~$T$ fixed and then taking the limit~$T\rightarrow\infty$. When~$t=T$, it is clear that~$V_{\beta,m,T}(\omega)$ is the maximum of a linear function of~$\omega$ and a constant function~($m$), and is thus continuous, convex and piecewise linear. By the induction hypothesis that~$V_{\beta,m,t+1}(\omega)$ is continuous, convex and piecewise linear, we again have that~$V_{\beta,m,t}(\omega)$ is the maximum of two continuous, convex and piecewise linear functions and is thus continuous, convex and piecewise linear. Therefore~$V_{\beta,m,t}(\omega)$ is continuous, convex and piecewise linear in~$\omega$ for all~$t\in\{1,2,\cdots, T\}$. Using~$\|\cdot\|\defeq\|\cdot\|_1$ norm on~$\mathbb{R}^K$ and consider two states~$\omega_1,\omega_2$ such that~$\|\omega_1-\omega_2\|>0$. At~$t=T$, we have
\begin{eqnarray}
|V_{\beta,m,T}(\omega_1)-V_{\beta,m,T}(\omega_2)| &=& \max\{\omega_1 B', m\} - \max\{\omega_2 B', m\}.\label{eq:bd1}
\end{eqnarray}
Without loss of generality, assume~$\omega_1 B'\le\omega_2 B'$. We consider the following 3 cases:\\
\noindent i) if~$m<\omega_1 B'$, then
\[|V_{\beta,m,T}(\omega_1)-V_{\beta,m,T}(\omega_2)|=|\omega_1B'-\omega_2B'|\le B_K\|\omega_1-\omega_2\|;\]
\noindent ii) if~$\omega_1 B'\le m\le \omega_2 B'$, then
\[|V_{\beta,m,T}(\omega_1)-V_{\beta,m,T}(\omega_2)|\le B_K\|\omega_1-\omega_2\|;\]
\noindent iii) if~$m> \omega_2 B'$, then
\[|V_{\beta,m,T}(\omega_1)-V_{\beta,m,T}(\omega_2)|=0.\]
From the above, we have that
\begin{eqnarray}
|V_{\beta,m,T}(\omega_1)-V_{\beta,m,T}(\omega_2)| &\le& B_K\|\omega_1-\omega_2\|.\label{eq:bd2}
\end{eqnarray}
At time~$t+1$, we make the following induction hypothesis that
\[|V_{\beta,m,t+1}(\omega_1)-V_{\beta,m,t+1}(\omega_2)| \le \frac{1-\beta^{T-t}}{1-\beta}B_K\|\omega_1-\omega_2\|.\]
At time~$t$, we have, by a similar case analysis as above, that
\begin{eqnarray}
|V_{\beta,m,t}(\omega_1)-V_{\beta,m,t}(\omega_2)| \le \frac{1-\beta^{T-t+1}}{1-\beta}B_K\|\omega_1-\omega_2\|.\label{eq:bd3}
\end{eqnarray}
Note that we have used the fact that~$\|\Tau^1(\omega_1)-\Tau^1(\omega_2)\|\le\|\omega_1-\omega_2\|$. Therefore, we have that
\begin{eqnarray}
|V_{\beta,m,1}(\omega_1)-V_{\beta,m,1}(\omega_2)| \le \frac{1-\beta^{T}}{1-\beta}B_K\|\omega_1-\omega_2\|\le\frac{1}{1-\beta}B_K\|\omega_1-\omega_2\|.\label{eq:bd4}
\end{eqnarray}
Furthermore, for all~$t\in\{1,2,\cdots,T\}$, we have that
\begin{eqnarray}
|V_{\beta,m,t}(\omega_1)-V_{\beta,m,t}(\omega_2)| \le \frac{1}{1-\beta}B_K\|\omega_1-\omega_2\|.\label{eq:bd5}
\end{eqnarray}
This proves that the finite-horizon value function $V_{\beta,m,t}(\omega)$ is Lipschitz continuous in~$\omega$ with constant~$\frac{1}{1-\beta}B_K$, independent of horizon length~$T$ and starting point~$t$. Fix~$t=1$, if we can show as~$T$ goes to infinity~$V_{\beta,m,1}(\cdot)$ converges to~$V_{\beta,m}(\cdot)$ pointwise, then~$V_{\beta,m}(\cdot)$ must be Lipschitz continuous with the same constant. This is because that given any two states~$\omega_1,\omega_2$ and any~$\epsilon>0$, there exists a positive integer~$T_0$ such that
\[|V_{\beta,m}(\omega_1)-V_{\beta,m}(\omega_2)| \le 2\epsilon+|V_{\beta,m,1}(\omega_1)-V_{\beta,m,1}(\omega_2)|\le 2\epsilon+\frac{1}{1-\beta}B_K\|\omega_1-\omega_2\|.\]
Since~$\epsilon>0$ is arbitrary, the Lipschitz continuity of~$V_{\beta,m}(\cdot)$ follows. To prove the convergence of~$V_{\beta,m,1}(\cdot)$ with~$T$, we first apply the optimal policy~$\pi_{sa}^*(T)$ to the first~$T$ slots followed by an (stationary) optimal policy~$\pi^*_{sa}$ for the infinite-horizon problem in subsequent time slots~$t>T$, then
\begin{eqnarray}
V_{\beta,m}(\omega)\ge V_{\beta,m,1}(\omega)+\beta^{T}\mathbb{E}[V_{\beta,m}(\omega(T+1))],\label{eq:bd6}
\end{eqnarray}
where the expectation is taken with respect to~$\omega(T+1)$ which is determined by the past observations and actions in the first~$T$ slots. It is clear that~$V_{\beta,m}(\cdot)$ is bounded:
\begin{eqnarray}
0\le V_{\beta,m}(\cdot)\le\frac{\max\{B_K,m\}}{1-\beta}.\label{eq:bd10}
\end{eqnarray}
From~\eqref{eq:bd6} and~\eqref{eq:bd10}, we know that
\begin{eqnarray}
V_{\beta,m,1}(\omega)-V_{\beta,m}(\omega)\le0.\label{eq:bd7}
\end{eqnarray}
Now we apply~$\pi^*_{sa}$ to the finite-horizon problem with length~$T$ and compare the reward accumulated in the~$T$ slots:
\begin{eqnarray}
V_{\beta,m,1}(\omega)\ge V_{\beta,m}(\omega)-\beta^{T}\mathbb{E}[V_{\beta,m}(\omega(T+1))].\label{eq:bd8}
\end{eqnarray}
From~\eqref{eq:bd10},~\eqref{eq:bd7} and~\eqref{eq:bd8}, we have, for any initial value of~$\omega$ at~$t=1$,
\begin{eqnarray}
-\beta^{T}\frac{\max\{B_K,m\}}{1-\beta}\le V_{\beta,m,1}(\omega)-V_{\beta,m}(\omega)\le0.\label{eq:bd9}
\end{eqnarray}
Taking the limit~$T\rightarrow\infty$, we proved the (uniform) convergence of~$V_{\beta,m,1}(\cdot)$ to~$V_{\beta,m}(\cdot)$. Consequently~$V_{\beta,m}(\cdot)$ is Lipschitz continuous. Its convexity is clear as a limiting function of convex functions.

Now we consider the properties of~$V_{\beta,m}(\omega)$ regarding to~$m$ with~$\omega$ fixed. By a similar argument as above, we have that~$V_{\beta,m,1}(\omega)$ is convex, continuous and piecewise linear in~$m$. Furthermore, it is Lipschitz continuous in~$m$ with constant~$\frac{1}{1-\beta}$, \ie, for any~$m_1$ and~$m_2$,
\[|V_{\beta,m_1,1}(\omega)-V_{\beta,m_2,1}(\omega)|\le \frac{1}{1-\beta}|m_1-m_2|.\]
It remains to show the pointwise convergence of~$V_{\beta,m,1}(\omega)$ to~$V_{\beta,m}(\omega)$ for every fixed~$m$ as~$T\rightarrow\infty$. However, it is a direct result of~\eqref{eq:bd9}.
\Halmos
\endproof

\section{Proof of Lemma~\ref{lm:boundary}.}\label{proof:lmBoundary}

\proof{} We prove the lemma step-by-step.

{\it Step 1.} We first show that $A(m)$ is convex. From~\eqref{eq: long time value function},~\eqref{eq: long time value function with active action},~\eqref{eq: long time value function with passive action} and~Lemma~\ref{lm:vf}, given any~$\omega_1,\omega_2\in A(m)$ and~$\lambda\in(0,1)$, we have
\begin{eqnarray}
V_{\beta,m}(\lambda\omega_1+(1-\lambda)\omega_2;u=1) &=& \lambda V_{\beta,m}(\omega_1;u=1) + (1-\lambda) V_{\beta,m}(\omega_2;u=1)\\
&>& \lambda V_{\beta,m}(\omega_1;u=0) + (1-\lambda) V_{\beta,m}(\omega_2;u=0)\\
&\ge& V_{\beta,m}(\lambda\omega_1+(1-\lambda)\omega_2;u=0).
\end{eqnarray}
The first equality in the above is due to the linearity of~$V_{\beta,m}(\cdot;u=1)$ by~\eqref{eq: long time value function with active action}, the second last inequality is by Definition~\eqref{def:active}, and the last inequality is due to the convexity of~$V_{\beta,m}(\cdot)$ established in Lemma~\ref{lm:vf} and the linearity of~$\Tau^1(\cdot)$ by~\eqref{eq: state dist update}. Therefore, the active set~$A(m)$ is convex.

{\it Step 2.} The openness of $A(m)$ is due to the strict inequality in \eqref{def:active} and the continuity of the value function established in Lemma~\ref{lm:vf}. Since every point in $A(m)$ has an $\epsilon$-neighborhood as a $(K-1)$-dimensional ball of $\mathbb{R}^{K-1}$ in $A(m)$, the dimension of $A(m)$ is $K-1$.

{\it Step 3.} It is obvious that the closure $\overline{A(m)}$ of the open convex set $A(m)$ in $\mathbb{X}$ is formed by the linear boundaries of the simplex space $X$ and $C(m)$. Therefore $C(m)$ is closed and bounded and thus compact as a subspace of $\mathbb{R}^{K-1}$. The proof that $C(m)$ is a simply connected $(K-2)$-dimensional subspace of $X$ requires familiarity to the theory of algebraic topology and is concisely sketched as follows. Since $\overline{A(m)}$ is also convex, it is {\it homeomorphic} to the open $(K-1)$-dimensional unit ball in $\mathbb{R}^{K-1}$ \citep{Ges2012}. Let $h$ denote this homeomorphism. Choose the extreme point $x_0=(0,0,\ldots,1)\in \overline{A(m)}-C(m)$. Note that the boundary of $A(m)$ consists of parts of the $(K-2)$-dimensional linear boundaries (hyperplanes) of $X$ and $C(m)$. Let $P$ be any of such hyperplanes that intersect with $C(m)$. Then there is a deformation retract from $P-x_0$ to the intersection $P\cap C(m)$ (\ie a continuous map from $P-x_0$ to $P\cap C(m)$ that is homotopic to the identity map on $P-x_0$ with $P\cap C(m)$ fixed during the homotopy). Followed by $h$, this induces a deformation retract from the punctured $(K-2)$-sphere $\mathbb{S}^{K-2}-h(x_0)$ to $h(C(m))$. Therefore the fundamental group of $C(m)$ is isomorphic to that of $\mathbb{S}^{K-2}-h(x_0)$ by the deformation retract (Theorem 58.3 on \citealt{Munkres2003}). Since the punctured sphere $\mathbb{S}^{K-2}-h(x_0)$ is homeomorphic to $\mathbb{R}^{K-2}$ (by the stereographic projection) which is simply connected, we proved that $C(m)$ is also simply connected. Finally, the homeomorphism $h$ shows that $C(m)$ is a $(K-2)$-dimensional subspace of $\mathbb{X}$ with $h(C(m))$ as a subset with a nonempty interior in $\mathbb{S}^{K-2}$.

{\it Step 4.} Now it should be clear that $C(m)$ and the linear boundaries of the simplex space $X$ containing the extreme points in $A(m)$ form the closure of the active set. While $C(m)$ and the linear boundaries of $X$ containing the rest of extreme points in $P(m)$ form the passive set. Therefore $C(m)$ partitions $X$ into {\it two disjoint and connected} subspaces: the active set and the passive set.

{\it Step 5.} There is a small point in the above argument worth a little more discussion. Since we include $C(m)$ in the passive set $P(m)$ by definition, is it possible that $C(m)$ is a thick boundary as a $(K-1)$-dimensional space (bulged in the direction to the interior of the passive set)? The answer is no. Because in this case the convex value function conditional on $u=0$ will always lie above the linear value function conditional on $u=1$, then the problem is reduced to the trivial scenario.

\Halmos
\endproof

\section{Proof of Theorem~\ref{prop:indexability}.}\label{proof:thmIdx}

\proof{} The existence of the right (or left) derivative follows directly from the convexity of~$V_{\beta,m}(\omega)$. Fix an~$m_0$ and apply a change~$\Delta m$ to the single-armed bandit, we have
\begin{eqnarray}
V_{\beta,m_0+\Delta m}(\omega)\ge V_{\beta,m_0}(\omega) + D_{\beta,m}(\omega){\Delta m}.\label{eq:compare6}
\end{eqnarray}
Now if we apply an optimal policy for the arm with subsidy~$m=m_0+\Delta m$ to the case of~$m=m_0$, we have
\begin{eqnarray}
V_{\beta,m_0}(\omega)\ge V_{\beta,m_0+\Delta m}(\omega) - D_{\beta,m+\Delta m}(\omega){\Delta m}.\label{eq:compare7}
\end{eqnarray}
From~\eqref{eq:compare6} and~\eqref{eq:compare7}, it is clear that
\begin{eqnarray}
D_{\beta,m}(\omega)\le\frac{V_{\beta,m_0+\Delta m}(\omega)-V_{\beta,m_0}(\omega)}{\Delta m}
\le D_{\beta,m+\Delta m}(\omega),~\quad\forall~\Delta m>0.\label{eq:compare8}
\end{eqnarray}
Note that the above implies the monotonically nondecreasing property of~$D_{\beta,m}(\omega)$ with~$m$. To prove~\eqref{eq:rightDiff}, we only need to show that~$D_{\beta,m}(\omega)$ is right continuous in~$m$. Assume this is {\em not} true so there exists a decreasing sequence~$\{m_k\}$ converging to~$m_0$ and an~$\epsilon>0$ such that
\begin{eqnarray}
D_{\beta,m_k}(\omega)-D_{\beta,m_0}(\omega)>\epsilon.\label{eq:compare9}
\end{eqnarray}
Since~$D_{\beta,m_k}(\omega)$ has a value ranging in the compact set~$[0,\frac{1}{1-\beta}]$, we can find a convergent subsequence~$\{m_{k_i}\}$ of~$\{m_k\}$ such that
\begin{eqnarray}
\lim_{i\rightarrow\infty}D_{\beta,m_{k_i}}(\omega)=D>D_{\beta,m_0}(\omega),\label{eq:converge1}
\end{eqnarray}
where~$D\in(0,\frac{1}{1-\beta}]$ is the limit of the passive time as~$m_{k_i}\rightarrow m_0$. If we can show that~$D$ can be achieved by a policy~$\pi^*\in\Pi^*_{sa}(m_0)$, then we have a contradiction to~\eqref{def:passiveTime}.

To construct~$\pi^*$ with passive time~$D$ and achieving~$V_{\beta,m_0}(\omega)$, we look at a finite horizon~$T$. Starting from the initial belief state~$\omega$, the possible belief states within~$T$ must be finite, leading to a finite set of possible policies. Specifically, if the number of possible states to observe is~$h(T)$, the number of policies up to time~$T$ is at most~$2^{h(T)}$ as each state is applied with either~$u=0$ or~$u=1$. We can thus choose a subsequence~$\{m_j(T)\}$ of~$\{m_{k_i}\}$ such that the optimal policy achieving $D_{\beta,m_{m_j(T)}}$ under~$m_j(T)$ is the same for all~$j$ within the first~$T$ slots. Repeat the process for slots~$T+1$ up to~$2T$ and keep doubling the time horizon, we arrive at a policy for all states that may happen at any time. For any time horizon~$T'$, this policy coincides with the optimal policies for a subsequence~$\{m_j(T'')\}$ of~$\{m_{k_i}\}$ for some~$T''>T'$ and by taking~$T'$ large enough, this policy achieves a passive time at least~$D-\epsilon_1$ and a total reward~$V_{\beta,m_0}(\omega)-\epsilon_1$ for any arbitrarily small~$\epsilon_1>0$ due to~\eqref{eq:converge1} and the continuity of~$V_{\beta,\cdot}(\omega)$. This policy is thus optimal for the infinite-horizon single-armed bandit problem with subsidy~$m_0$ with passive time~$D$, as desired for~$\pi^*$.

To prove the sufficiency of~\eqref{eq:diffIdx} and~\eqref{eq:diffIdx1}, we assume that the arm is not indexable, \ie there exists~$m_0$ and~$\omega\in C(m_0)\subset P(m_0)$ such that for any~$\epsilon>0$, we can find an~$m_1~(m_0<m_1<m_0+\epsilon)$ with~$\omega\in A(m_1)$. This means that as the boundary~$C(m)$ moves (continuously) as~$m$ increases, some belief state moves from the passive set to the active set. Under this scenario, we have
\begin{eqnarray}
V_{\beta,m_0}(\omega;u=1)= V_{\beta,m_0}(\omega;u=0). \label{eq:compare10}\\
V_{\beta,m_1}(\omega;u=1) > V_{\beta,m_1}(\omega;u=0).\label{eq:compare11}
\end{eqnarray}
According to~\eqref{eq: long time value function with active action} and~\eqref{eq: long time value function with passive action}, both~$V_{\beta,m}(\omega;u=1)$ and~$V_{\beta,m}(\omega;u=0)$ are right differentiable with~$m$ for any belief state~$\omega$, so is their difference. Therefore, by~\eqref{eq:compare10} and~\eqref{eq:compare11},
\begin{eqnarray}
\left.\frac{dV_{\beta,m}(\omega;u=1)}{(dm)^+}\right|_{m=m_0}&=&\lim_{m_1\rightarrow m_0}\frac{V_{\beta,m_1}(\omega;u=1)- V_{\beta,m_0}(\omega;u=1)}{m_1-m_0}\label{eq:compare12}\\
&\ge& \lim_{m_1\rightarrow m_0}\frac{V_{\beta,m_1}(\omega;u=0)-V_{\beta,m_0}(\omega;u=0)}{m_1-m_0}=\left.\frac{dV_{\beta,m}(\omega;u=0)}{(dm)^+}\right|_{m=m_0}.\label{eq:compare13}
\end{eqnarray}
This would contradict~\eqref{eq:diffIdx} unless the equality in~\eqref{eq:compare13} holds, which would contradict~\eqref{eq:diffIdx1} given~\eqref{eq:compare11} and that~$\epsilon$ can be chosen arbitrarily small.

To prove the necessity of~\eqref{eq:diffIdx} and~\eqref{eq:diffIdx1}, assume there exists an~$\omega\in C(m_{\omega})$ such that
\begin{eqnarray}
\left.\frac{dV_{\beta,m}(\omega;u=0)}{(dm)^+}\right|_{m=m_{\omega}} < \left.\frac{dV_{\beta,m}(\omega;u=1)}{(dm)^+}\right|_{m=m_{\omega}},\label{eq:contra1}
\end{eqnarray}
and when
\begin{eqnarray}
\left.\frac{dV_{\beta,m}(\omega;u=0)}{(dm)^+}\right|_{m=m_{\omega}} = \left.\frac{dV_{\beta,m}(\omega;u=1)}{(dm)^+}\right|_{m=m_{\omega}},\label{eq:contra2}
\end{eqnarray}
for any~$\epsilon_1>0$, there exists an~$m_2~(m_{\omega}<m_2<m_{\omega}+\epsilon_1)$ such that
\begin{eqnarray}
V_{\beta,m_2}(\omega;u=0)< V_{\beta,m_2}(\omega;u=1).\label{eq:contra3}
\end{eqnarray}
By~\eqref{eq:contra1}, there exists~$\Delta m>0$ such that
\begin{eqnarray}
V_{\beta,m(\omega)+\Delta m}(\omega;u=0)- V_{\beta,m(\omega)}(\omega;u=0)< V_{\beta,m(\omega)+\Delta m}(\omega;u=1)- V_{\beta,m(\omega)}(\omega;u=1).
\end{eqnarray}
Together with the fact that
\begin{eqnarray}
V_{\beta,m(\omega)}(\omega;u=0)= V_{\beta,m(\omega)}(\omega;u=1),
\end{eqnarray}
we have that~$\omega\in A(m(\omega)+\Delta m)$ and obtained a contradiction to indexability as~$\omega\in C(\omega)\subset P(m(\omega))$. Furthermore, when~\eqref{eq:contra2} holds, it is straightforward that~\eqref{eq:contra3} contradicts~\eqref{eq:diffIdx1}.
\Halmos
\endproof

\section{Proof of Theorem~\ref{thm:indexableBeta}.}\label{proof:thmIdxBeta}

\proof{} Note that~\eqref{eq:diffIdx} is equivalent to
  \begin{eqnarray}
    \beta\sum_{k=0}^{K-1}\omega_kD_{\beta,m}(p_k)\le 1+\beta D_{\beta,m}(\Tau^1(\omega)),~\quad\forall~\omega\in C(m).\label{eq:compare14}
  \end{eqnarray}
  The above clearly holds if~$\beta\le 0.5$ as~$D_{\beta,m}(\cdot)$ is lower and upper bounded by $0$ and~$1-\beta$, respectively. The strict inequality in~\eqref{eq:compare14} holds if~$\beta<0.5$, satisfying~$\eqref{eq:diffIdx}$. When~$\beta=0.5$, the equality in~\eqref{eq:compare14} holds if~$\beta\sum_{k=0}^{K-1}\omega_kD_{\beta,m}(p_k)=\frac{\beta}{1-\beta}=1$ and~$D_{\beta,m}(\Tau^1(\omega))=0$. In this case, as~$m$ keeps increasing, the left-hand side of~\eqref{eq:compare14} can not increase while the right-hand side cannot decrease. Apply any~$\Delta m>0$ to~$\omega\in C(m)$, we have
  \begin{eqnarray}
    V_{\beta,m+\Delta m}(\omega;u=0)&\ge& \Delta m (1+\beta D_{\beta,m}(\Tau^1(\omega))) + V_{\beta,m}(\omega;u=0) \\
    &=& \Delta m (1+\beta D_{\beta,m}(\Tau^1(\omega))) + V_{\beta,m}(\omega;u=1) \\
    &=& \Delta m (1+\beta D_{\beta,m}(\Tau^1(\omega))) + \omega B'+\beta\frac{m}{1-\beta} \\
    &=& \omega B'+\beta\frac{m+\Delta m}{1-\beta} \\
    &=& V_{\beta,m+\Delta m}(\omega;u=1),
  \end{eqnarray}
  where the second last equality is due to that $D_{\beta,m}(\Tau^1(\omega))=0$ and $\beta=0.5$. The last equality is due to the fact that any future state~$p_k$ after activating at~$\omega$ must remain in the passive set as~$m$ increases due to the monotonic nondecreasing property of $D_{\beta,m}(p_k)$ with~$m$ and that~$D_{\beta,m}(p_k)$ is already equal to the upper bound~$\frac{1}{1-\beta}$. Therefore~$\eqref{eq:diffIdx1}$ is satisfied as well.
\Halmos
\endproof

\section{Proof of Lemma~\ref{lemma:linearEq1}.}\label{proof:lmaLE}

\proof{} It is helpful to rewrite~\eqref{eq: linearET1} and~\eqref{eq: linearET2} in the following matrix form~$AX=b$:
{\small\begin{align}\label{eq: value function in matrix form}
	\left[
	\begin{matrix}
		\textbf{I}_K - \beta\left(
		\begin{matrix}
			\beta^{L(p_0,\omega)}p_0\textbf{P}^{L(p_0,\omega)}\\
			\beta^{L(p_1,\omega)}p_1\textbf{P}^{L(p_1,\omega)}\\
            \vdots\\
			\beta^{L(p_{K-1},\omega)}p_K\textbf{P}^{L(p_{K-1},\omega)}\\
		\end{matrix}
		\right)
	\end{matrix}
	\right]\times
    \left(
	\begin{matrix}
		V_{\beta,m}(p_{0})\\V_{\beta,m}(p_{1})\\\vdots\\V_{\beta,m}(p_{K-1})
	\end{matrix}
	\right) = \left[
	\left(
	\begin{matrix}
		\frac{1 - \beta^{L(p_0,\omega)}}{1 - \beta}\\
		\frac{1 - \beta^{L(p_1,\omega)}}{1 - \beta}\\
        \vdots\\
		\frac{1 - \beta^{L(p_{K-1},\omega)}}{1 - \beta}\\
	\end{matrix}
	\right)m + \left(
	\begin{matrix}
		\beta^{L(p_0,\omega)}p_0\textbf{P}^{L(p_0,\omega)}\\
		\beta^{L(p_1,\omega)}p_1\textbf{P}^{L(p_1,\omega)}\\
        \vdots\\
		\beta^{L(p_{K-1},\omega)}p_{K-1}\textbf{P}^{L(p_{K-1},\omega)}\\
	\end{matrix}
	\right)B'
	\right]
\end{align}}
To prove the claim, we only need to show the coefficient matrix is invertible. By the Perron-Frobenius theorem, the eigenvalues~$\{\lambda_i\}_{i=1}^h$ of the following matrix satisfy~$|\lambda_i|\le 1$ for all~$i\in\{1,\cdots,h\}$ since it is a transition matrix with nonnegative elements and the sum of each row is equal to~$1$:
\begin{align}
\left(
		\begin{matrix}
			p_0\textbf{P}^{L(p_0,\omega)}\\
			p_1\textbf{P}^{L(p_1,\omega)}\\
            \vdots\\
			p_{K-1}\textbf{P}^{L(p_{K-1},\omega)}\\
		\end{matrix}
		\right) = Q\left(
        \begin{matrix}
            \lambda_1&\\
                     &~\ddots\\
                     &~&~~~\lambda_i&~~1\\
                     &~&~&~~~~\lambda_i\\
                     &~&~&~&~&~~\ddots\\
                     &~&~&~&~&~&~~~~\lambda_h
        \end{matrix}
\right)Q^{-1},
\end{align}
where the above equation shows the Jordan canonical form of the matrix and the square matrix~$Q$ has full rank~$K$. Therefore we can rewrite the coefficient matrix~$A$ as
\begin{align}
	A&=\left[
	\begin{matrix}
		\textbf{I}_K - \beta\left(
		\begin{matrix}
			\beta^{L(p_0,\omega)}p_0\textbf{P}^{L(p_0,\omega)}\\
			\beta^{L(p_1,\omega)}p_1\textbf{P}^{L(p_1,\omega)}\\
            \vdots\\
			\beta^{L(p_{K-1},\omega)}p_K\textbf{P}^{L(p_{K-1},\omega)}\\
		\end{matrix}
		\right)
	\end{matrix}
	\right]\nn\\&=Q\left(
        \begin{matrix}
            1-\beta^{L(p_0,\omega)+1}\lambda_1&\\
                     &\ddots\\
                     &~&1-\beta^{L(p_j,\omega)+1}\lambda_i&-\beta^{L(p_i,\omega)+1}\\
                     &~&~&1-\beta^{L(p_j,\omega)+1}\lambda_i\\
                     &~&~&~&~&\ddots\\
                     &~&~&~&~&~&1-\beta^{L(p_{K-1},\omega)+1}\lambda_h
        \end{matrix}
\right)Q^{-1}\nn\\
&=QJQ^{-1}.
\end{align}
It is easy to see that no eigenvalue of~$J$ can be zero so it has a full rank, leading to the full rank of~$A$ as it is similar to~$J$.
\Halmos
\endproof

\section{Proof of Theorem~\ref{prop:relaxIdx}.}\label{proof:thmRI}
\proof{} According to the remark following Lemma~\ref{lemma:linearEq1}, the value function~$\hat{V}_{\beta,m}(\omega_1)$ is linear in~$m$ for any~$\omega_1$ because the threshold policy is independent of~$m$ if~$\omega^*_{\beta}(m)$ is fixed. Since the subsidy~$m$ is paid if and only if the arm is made passive, the linear coefficient of~$m$ in~$\hat{V}_{\beta,m}(\omega_1)$ is simply~$\hat{D}_{\beta}(\omega_1)$. The passive time $\hat{D}_{\beta}(\omega_1)$ is clearly independent of~$m$ conditional on the fixed threshold. Since~\eqref{eq:equal at threshold} has a unique solution for~$m$ if and only if its left and right hand sides have different coefficients of~$m$, we proved the equivalence of~\eqref{eq: relaxedIdxIff} to the relaxed indexability. The expression~\eqref{eq: Approximated Whittle} of the approximated Whittle index follows directly from the unique solution of~\eqref{eq: linearET1},~\eqref{eq: linearET2} and~\eqref{eq:equal at threshold} under the relaxed indexability.
\Halmos
\endproof

\section{The Case of $K=3$}\label{sec:k=3}

In this section, we consider the case that an arm has a $2$-dimensional belief state space, \ie the internal Markov chain has 3 states. For simplicity of presentation, we assume that the Markov chain is irreducible and aperiodic, thus having a unique stationary (limiting) distribution.

\subsection{The Jordan Canonical Form}

To compute~$L(\cdot,\omega)$, it is crucial to analyze the form of~$\textbf{P}^k$ with~$k$. A general approach is to use the Jordan canonical form of the stochastic matrix when computing its power. It is well known that any~$K\times K$ square matrix~$\textbf{P}$ can be written in its Jordan form as
\begin{align}
\textbf{P} = QJQ^{-1}= Q\left(
        \begin{matrix}
            J_0&\\
                ~     &\ddots\\
                 ~    &~&J_{V-1}
        \end{matrix}
\right)Q^{-1},\quad J_v = \left(
        \begin{matrix}
            \lambda_v& ~1 &\\
                    & ~\lambda_v&~\ddots &\\
                     &&~~\ddots&~1\\
                     &&&~~\lambda_v
        \end{matrix}
\right),
\end{align}
where~$Q$ is a square matrix of full rank~$K$ and the upper diagonal~$K_v\times K_v$ matrix~$J_v$ is the~$v$-th {\em Jordan block} of size~$K_v~(1\le K_v\le K)$ with the eigenvalue~$\lambda_v$ and~$\sum_{v=0}^{V-1} = K$. Note that if~$K_v=1$ then the Jordan block is simply a scalar~$(\lambda_v)$ and different blocks can have the same eigenvalue, \ie there may exist~$0\le v_1\neq v_2\le V-1$ such that~$\lambda_{v_1}=\lambda_{v_2}$. The~$k$th power of~$\textbf{P}$ can thus be computed as
\begin{align}
\textbf{P}^k = Q\left(
        \begin{matrix}
            J_0^k&\\
                ~     &\ddots\\
                 ~    &~&J_{V-1}^k
        \end{matrix}
\right)Q^{-1},\quad J_v^k = \left(
        \begin{matrix}
            \lambda_v^k& \binom{k}{1}\lambda_v^{k-1} & \binom{k}{2}\lambda_v^{k-2}&\cdots&\binom{k}{K_v-1}\lambda_v^{k-K_v+1}\\
            0&\lambda_v^k &\binom{k}{1}\lambda_v^{k-1}&\cdots&\binom{k}{K_v-2}\lambda_v^{k-K_v+2}\\
            \vdots&\vdots&\ddots&\ddots& \vdots\\
            0&0&\cdots&\lambda_v^k& \binom{k}{1}\lambda_v^{k-1}\\
            0&0&\cdots&0& \lambda_v^k\\
        \end{matrix}
\right). \label{eq:jordanPower}
\end{align}
For a finite irreducible and aperiodic Markov chain, the transition matrix~$\textbf{P}$ is {\em regular}, \ie there exists an integer~$k\ge1$ such that~$\textbf{P}^k>0$ (element-wise). Therefore its Perron-Frobenius eigenvalue~$\lambda_{pf}=1$ has algebraic multiplicity~$1$, \ie the Jordan block associated with the eigenvalue~$1$ is unique and of size~$1$. Furthermore, any other eigenvalue~$\lambda_v\neq 1~(1\le v\le K)$ of~$\textbf{P}$ satisfies~$|\lambda_v|<1$. In this case, the Markov chain has a unique stationary distribution to which~$\omega_1\textbf{P}^k$ converges at a geometric rate as~$k\rightarrow\infty$ for any belief sate~$\omega_1$.

In the case of~$K=3$, the Jordan canonical form of the transition matrix~$\textbf{P}$ takes one of the following two forms (assuming irreducible and aperiodic Markov chains):
\begin{itemize}

  \item[1.] $\textbf{P}$ has~$3$ linearly independent eigenvectors: there exist~$\lambda_1,\lambda_2\in\mathbb{R}$ or~$\lambda_1=\overline{\lambda_2}\in\mathbb{C}$ with~$|\lambda_1|, |\lambda_2|\in[0,1)$,
  \begin{equation}\label{def: 1st Jordan form}
    J_{(1)} = \left(
    \begin{matrix}
        1&0&0\\0&\lambda_1&0\\0&0&\lambda_2
    \end{matrix}
    \right);
  \end{equation}
  \item[2.] $\textbf{P}$ has~$2$ linearly independent eigenvectors: there exists~$\lambda_1\in\mathbb{R},~|\lambda_1|\in[0,1)$,
  \begin{equation}\label{def: 2nd Jordan form}
    J_{(2)} = \left(
        \begin{matrix}
            1&0&0\\0&\lambda_1&1\\0&0&\lambda_1
        \end{matrix}
    \right).
  \end{equation}
\end{itemize}
Since the eigenvalue~$1$ corresponds to a single Jordan block of size~$1$ under our assumption, the matrix~$\textbf{P}$ has at least~$2$ linearly independent eigenvectors.

\subsection{The $k$-Step Reward Function}

Fix a belief state~$\omega$. Define the~$k$-step reward function as
\begin{eqnarray}
h(k) = \Tau^k(\omega)B' = \omega\textbf{P}^kB',\quad k\ge0.
\end{eqnarray}
To analyze~$L(\omega,\omega^*)$ for any threshold~$\omega^*$, we only need to find out the maximum of~$h(k)$ and if it exceeds~$\omega^*B'$, the first~$k$ that makes~$h(k)>\omega^*B'$. In the following lemma, we show that~$h(k)$ can only take three forms and then establish a detailed form of~$L(\omega,\omega^*)$ for the three cases respectively in Sec.~\ref{sec:computeL}.

\begin{lemma}\label{lemma:kStepForm}
The~$k$-step reward function~$h(k)$ takes one of the following three forms:
\begin{itemize}

  \item[1.] $\textbf{P}$ has only real eigenvalues and~$3$ linearly independent eigenvectors:
  \[h(k) = a_1b_1^k + a_2b_2^k + c,\quad a_1, b_1, a_2, b_2, c\in\mathbb{R},\ |b_1|,|b_2| < 1;\]
  \item[2.] $\textbf{P}$ has only real eigenvalues and~$2$ linearly independent eigenvectors:
  \[h(k) = ab^k + ckb^{k-1} + d,\quad a,b,c,d\in\mathbb{R},\  |b|<1;\]
  \item[3.] $\textbf{P}$ has a pair of conjugate complex eigenvalues:
  \[h(k) = a'A^k\sin(k\theta + b') + c',\quad a',A,b',c',\theta\in\mathbb{R},\ A\in(0,1),\ a' \ge 0,\ \theta\in(0,2\pi),\ b'\in[0,2\pi).\]
\end{itemize}
\end{lemma}
\proof{Proof.} Case~$1$ and~$2$ follow directly from the power of Jordan matrices with~$b_1=\lambda_1$, $b_2=\lambda_2$, or~$b=\lambda_1=\lambda_2$:
$$J_{(1)}^k = \left(
    \begin{matrix}
        1&0&0\\0&\lambda_1^k&0\\0&0&\lambda_2^k
    \end{matrix}
\right),\quad
J_{(2)}^k = \left(
    \begin{matrix}
        1&0&0\\0&\lambda_1^k&k\lambda_1^{k-1}\\0&0&\lambda_1^k
    \end{matrix}
\right).$$
For Case~$3$, write~$\textbf{P}=QJ_{(1)}Q^{-1}$. We have that~$\lambda_2 = \overline{\lambda_1}$. Let $Q = \{q_{ij}\}_{i,j=0,1,2},\ Q^{-1} =
\{\tilde{q}_{ij}\}_{i,j=0,1,2}$ and $Q_i = (q_{0i}, q_{1i}, q_{2i})',\ \tilde{Q}_i = (\tilde{q}_{i0},
\tilde{q}_{i1}, \tilde{q}_{i2}),\ i=0,1,2$. Then~$Q_2 = \overline{Q}_1,~\tilde{Q}_2 = \overline{\tilde{Q}_1}$
and~$Q_0,\tilde{Q}'_0\in\mathbb{R}^3$ (since they are respectively the right and left eigenvectors of~$\textbf{P}$ corresponding to the
eigenvalue~$1$):
\begin{align}
    h(k)
    &= \omega\left(Q_0\tilde{Q}_0 + \lambda_1^kQ_1\tilde{Q}_1 + \overline{\lambda_1^k}Q_2\tilde{Q}_2\right)B'
    \nonumber\\
    &= \omega Q_0\tilde{Q}_0B' + 2Re(\lambda_1^k\omega Q_1\tilde{Q}_1B')\quad (\mbox{Let}~r+si=\omega Q_1\tilde{Q}_1B',~\lambda_1=Ae^{i\theta})\nonumber\\
    &= \omega Q_0\tilde{Q}_0B' + 2A^k(r\cos k\theta-s\sin k\theta)\quad (\mbox{Let}~a'\sin(k\theta+b')=2(r\cos k\theta-s\sin k\theta),~c'=\omega Q_0\tilde{Q}_0B')\nonumber\\
    &= a'A^k\sin(k\theta + b')+ c',\nonumber
\end{align}
where~$\lambda_1=Ae^{i\theta},\ A\in(0,1),\ \theta\in(0,2\pi)$. Without loss of generality, we choose~$a' \ge 0,\ b'\in[0,2\pi)$.\Halmos
\endproof

\subsection{The Computation of~$L(\omega,\omega^*)$}\label{sec:computeL}

In the following theorem, we give the forms of the first crossing time~$L(\omega,\omega^*)$ for various cases mentioned in Lemma~\ref{lemma:kStepForm}.

\begin{theorem}\label{thm: firstCross}
Fix~$\omega$ and $\omega^*$. Let~$r^*\defeq\omega^*B'$. The first crossing time~$L(\omega,\omega^*)$ takes following forms:
\begin{eqnarray}
L(\omega,\omega^*) =
\begin{cases}
0,\quad \text{if}~h(0)>r^*\\
1,\quad \text{if}~h(1)>r^*\ge h(0)\\
2,\quad \text{if}~h(2)>r^*\ge \max\{h(0),h(1)\}
\end{cases},\label{eq:baseCase}
\end{eqnarray}
where~$h(k)$ is the~$k$-step reward function that depends on~$\omega$. The other cases are summarized below.
\begin{itemize}
  \item[1.] $\textbf{P}$ has only real eigenvalues and~$3$ linearly independent eigenvectors: $h(k) = a_1b_1^k + a_2b_2^k + c$.
  \begin{itemize}
  \item[1.1] $b_1=b_2\neq0~\&\&~b_1>0~\&\&~a_1+a_2<0$:
  \begin{eqnarray}
  L(\omega,\omega^*) =
  \begin{cases}
  \lfloor\log_{b_1}^{\frac{c-r^*}{-(a_1+a_2)}}\rfloor+1,\quad\text{if}~h(0)\le r^*<c\\
  \infty, \quad\text{if}~r^*\ge c
  \end{cases};
  \end{eqnarray}
  \item[1.2] $b_1=b_2\neq0~\&\&~(b_1<0~||~a_1+a_2\ge0)$: $L(\omega,\omega^*)=\infty$ if~$r^*\ge \max\{h(0),h(1)\}$;
  \item[1.3] $a_1b_1=0~\&\&~b_2>0~\&\&~a_2<0$:
  \begin{eqnarray}
  L(\omega,\omega^*) =
  \begin{cases}
  \lfloor\log_{b_2}^{\frac{c-r^*}{-a_2}}\rfloor+2,\quad\text{if}~\max\{h(0),h(1)\}\le r^*<c\\
  \infty, \quad\text{if}~r^*\ge \max\{h(0),c\}
  \end{cases};
  \end{eqnarray}
  \item[1.4] $a_1b_1=0~\&\&~(b_2\le0~||~a_2\ge0)$: $L(\omega,\omega^*)=\infty$ if~$r^*\ge \max\{h(0),h(1),h(2)\}$;
  \item[1.5] $a_2b_2=0~\&\&~b_1>0~\&\&~a_1<0$:
  \begin{eqnarray}
  L(\omega,\omega^*) =
  \begin{cases}
  \lfloor\log_{b_1}^{\frac{c-r^*}{-a_1}}\rfloor+2,\quad\text{if}~\max\{h(0),h(1)\}\le r^*<c\\
  \infty, \quad\text{if}~r^*\ge \max\{h(0),c\}
  \end{cases};
  \end{eqnarray}
  \item[1.6] $a_2b_2=0~\&\&~(b_1\le0~||~a_1\ge0)$: $L(\omega,\omega^*)=\infty$ if~$r^*\ge \max\{h(0),h(1),h(2)\}$;
  \item[1.7] $a_1,a_2,b_1,b_2>0$: $L(\omega,\omega^*)=\infty$ if~$r^*\ge h(0),h(1)$;
  \item[1.8] $a_1<0,a_2>0,b_1>b_2>0$:
  \begin{eqnarray}
  L(\omega,\omega^*) =
  \begin{cases}
  \min\{k:~k>\max\{\lceil\log_{\frac{b_1}{b_2}}^{-\frac{a_2(1-b_2)}{a_1(1-b_1)}}\rceil,0\},~h(k)>r^*\},\quad\text{if}~h(0)\le r^*<c\\
  \infty, \quad\text{if}~r^*\ge \max\{h(0),c\}
  \end{cases};
  \end{eqnarray}
  \item[1.9] $a_1<0,a_2>0,b_2>b_1>0$:
  \begin{eqnarray}
  L(\omega,\omega^*) =
  \begin{cases}
  \min\{k:~0<k\le\lfloor\log_{\frac{b_1}{b_2}}^{z_0}\rfloor+1,~h(k)>r^*\},\quad\text{if}~z_0<1
  ~\&\&~h(\lfloor\log_{\frac{b_1}{b_2}}^{z_0}\rfloor+1)>r^*\ge h(0)\\
  \infty, \quad\text{if}~(z_0\ge1\&\&r^*\ge h(0))~||~(z_0<1~\&\&~h(\lfloor\log_{\frac{b_1}{b_2}}^{z_0}\rfloor+1)\le r^*)
  \end{cases}
  \end{eqnarray}
  where~$z_0=-\frac{a_2(1-b_2)}{a_1(1-b_1)}$;
  \item[1.10] $b_1<0,a_1,a_2,b_2>0$: $L(\omega,\omega^*)=\infty$ if~$r^*\ge h(0)$;
  \item[1.11] $a_1,b_1<0,a_2,b_2>0$: $L(\omega,\omega^*)=\infty$ if~$r^*\ge \max\{h(0),h(1)\}$;
  \item[1.12] $a_2,b_1<0,a_1,b_2>0,|b_1|>b_2$:
  \begin{eqnarray}
  L(\omega,\omega^*) =
  \begin{cases}
  \min\{k:~0<k\le\lfloor\log_{\frac{b_1}{b_2}}^{z_1}\rfloor+2,~h(k)>r^*\},\quad\text{if}~z_1\ge1
  ~\&\&~h(\lfloor\log_{\frac{b_1}{b_2}}^{z_1}\rfloor+2)>r^*\ge h(0)\\
  \infty, \quad\text{if}~(z_1<1\&\&r^*\ge h(0))~||~(z_1\ge1~\&\&~h(\lfloor\log_{\frac{b_1}{b_2}}^{z_1}\rfloor+2)\le r^*)
  \end{cases}
  \end{eqnarray}
  where~$z_1=-\frac{a_2(1-b_2^2)}{a_1(1-b_1^2)}$;
  \item[1.13] $a_2,b_1<0,a_1,b_2>0,|b_1|<b_2$:
  \begin{eqnarray}
  L(\omega,\omega^*) =
  \begin{cases}
  \min\{k:~k\ge\max\{\lceil\log_{\frac{b_1}{b_2}}^{-\frac{a_2(1-b_2^2)}{a_1(1-b_1^2)}}\rceil,1\},~h(k)>r^*\},\quad\text{if}
  ~h(0)\le r^*< c\\
  \infty, \quad\text{if}~\max\{h(0),c\}\le r^*
  \end{cases};
  \end{eqnarray}
  \item[1.14] $a_2,b_1<0,a_1,b_2>0,|b_1|=b_2$:
  \begin{eqnarray}
  L(\omega,\omega^*) =
  \begin{cases}
  \min\{k:~k\ge1,~h(k)>r^*\},\quad\text{if}~h(0)\le r^*< c\\
  \infty, \quad\text{if}~\max\{h(0),c\}\le r^*
  \end{cases};
  \end{eqnarray}
  \item[1.15] $b_1,b_2<0,a_1,a_2>0$: $L(\omega,\omega^*)=\infty$ if~$r^*\ge h(0)$;
  \item[1.16] $b_1,b_2>0,a_1,a_2<0$:
  \begin{eqnarray}
  L(\omega,\omega^*) =
  \begin{cases}
  \min\{k:~k\ge0,~h(k)>r^*\},\quad\text{if}~r^*<c\\
  \infty, \quad\text{if}~c\le r^*
  \end{cases};
  \end{eqnarray}
  \item[1.17] $a_1,a_2,b_1<0,b_2>0,|b_1|>b_2$:
  \begin{eqnarray}
  L(\omega,\omega^*) =
  \begin{cases}
  \min\{k:~0<k\le\lfloor\log_{\frac{-b_1}{b_2}}^{z_2}\rfloor+2,~h(k)>r^*\},\quad\text{if}~z_2>\frac{-b_1}{b_2}
  ~\&\&~h(\lfloor\log_{\frac{-b_1}{b_2}}^{z_2}\rfloor+2)>r^*\ge h(0)\\
  \infty, \quad\text{if}~(z_2\le\frac{-b_1}{b_2}\&\&r^*\ge h(1))~||~(z_2>\frac{-b_1}{b_2}~\&\&~h(\lfloor\log_{\frac{-b_1}{b_2}}^{z_2}\rfloor+2)\le r^*)
  \end{cases}
  \end{eqnarray}
  where~$z_2=\frac{a_2(b_2^2-1)}{a_1(b_1^2-1)}$;
  \item[1.18] $a_1,a_2,b_1<0,b_2>0,|b_1|<b_2$:
  \begin{eqnarray}
  L(\omega,\omega^*) =
  \begin{cases}
  \min\{k:~k\ge\max\{\lceil\log_{\frac{-b_1}{b_2}}^{z_2}\rceil,1\},~h(k)>r^*\},\quad\text{if}~h(1)\le r^*< c\\
  \infty, \quad\text{if}~\max\{h(1),c\}\le r^*
  \end{cases}
  \end{eqnarray}
  where~$z_2=\frac{a_2(b_2^2-1)}{a_1(b_1^2-1)}$;
  \item[1.19] $a_1,a_2,b_1<0,b_2>0,|b_1|=b_2$:
  \begin{eqnarray}
  L(\omega,\omega^*) =
  \begin{cases}
  \min\{k:~k>1,~h(k)>\omega^*\},\quad\text{if}~h(1)\le r^*< c\\
  \infty, \quad\text{if}~\max\{h(1),c\}\le r^*
  \end{cases};
  \end{eqnarray}
  \item[1.20] $a_2,b_1,b_2<0,a_1>0,|b_1|>|b_2|$:
  \begin{eqnarray}
  L(\omega,\omega^*) =
  \begin{cases}
  \min\{k:~1\le k\le\lfloor\log_{\frac{b_1}{b_2}}^{z_1}\rfloor_e+2 ,~h(k)>r^*\},\quad\text{if}~z_1>1~\&\&~h(0)\le r^*< h(\lfloor\log_{\frac{b_1}{b_2}}^{z_1}\rfloor_e+2)\\
  \infty, \quad\text{if}~(z_1>1~\&\&~\max\{h(0),h(1),h(\lfloor\log_{\frac{b_1}{b_2}}^{z_1}\rfloor_e+2)\}\le r^*)~||~(z_1\le1~\&\&~h(0)\le r^*)
  \end{cases}
  \end{eqnarray}
  where~$z_1=-\frac{a_2(1-b_2^2)}{a_1(1-b_1^2)}$ and~$\lfloor n\rfloor_e$ denotes the maximum even integer not exceeding~$n$;
  \item[1.21] $a_2,b_1,b_2<0,a_1>0,|b_1|<|b_2|$:
  \begin{eqnarray}
  L(\omega,\omega^*) =
  \begin{cases}
  \min\{k:~1\le k\le\lfloor\log_{\frac{b_1}{b_2}}^{z_1}\rfloor_o+2 ,~h(k)>r^*\},\quad\text{if}~z_1<\frac{b_1}{b_2}~\&\&~h(0)\le r^*< h(\lfloor\log_{\frac{b_1}{b_2}}^{z_1}\rfloor_o+2)\\
  \infty, \quad\text{if}~(z_1\ge\frac{b_1}{b_2}~\&\&~\max\{h(0),h(1)\}\le r^*)~||~(z_1<\frac{b_1}{b_2}~\&\&~\{h(0),h(1),h(\lfloor\log_{\frac{b_1}{b_2}}^{z_1}\rfloor_o+2)\}\le r^*)
  \end{cases}
  \end{eqnarray}
  where~$z_1=-\frac{a_2(1-b_2^2)}{a_1(1-b_1^2)}$ and~$\lfloor n\rfloor_o$ denotes the maximum odd integer not exceeding~$n$;
  \item[1.22] $a_2,b_1,b_2<0,a_1<0$: $L(\omega,\omega^*)=\infty$ if~$\omega^*\ge h(1)$;
  \item[1.23] any other case, which is symmetric to one of the above and omitted here.
  \end{itemize}
  \item[2.] $\textbf{P}$ has only real eigenvalues and~$2$ linearly independent eigenvectors: $h(k) = ab^k + ckb^{k-1} + d$.
  \begin{itemize}
  \item[2.1] $b,c>0$:
  \begin{eqnarray}
  L(\omega,\omega^*) =
  \begin{cases}
  \min\{k:~0\le k<\lceil z_3\rceil+1 ,~h(k)>r^*\},\quad\text{if}~z_3>0~\&\&~h(0)\le r^*< h(z_3+1)\\
  \infty, \quad\text{if}~(z_3>0~\&\&~h(z_3+1)\le r^*)~||~(z_3\le0~\&\&~h(0)\le r^*)
  \end{cases}
  \end{eqnarray}
  where~$z_3=\frac{ab-ab^2-cb}{c(b-1)}$;
  \item[2.2] $b>0,c<0$:
  \begin{eqnarray}
  L(\omega,\omega^*) =
  \begin{cases}
  \min\{k:~k\ge\max\{\lceil z_3\rceil+1,0\},~h(k)>r^*\},\quad\text{if}~h(0)\le r^*< d\\
  \infty, \quad\text{if}~\max\{h(0),d\}\le r^*
  \end{cases}
  \end{eqnarray}
  where~$z_3=\frac{ab-ab^2-cb}{c(b-1)}$;
  \item[2.3] $b<0,c<0$:
  \begin{eqnarray}
  L(\omega,\omega^*) =
  \begin{cases}
  \min\{k:~1\le k\le\lfloor z_4\rfloor_e+2 ,~h(k)>r^*\},\quad\text{if}~z_4>0~\&\&~h(0)\le r^*< h(\lfloor z_4\rfloor_e+2)\\
  \infty, \quad\text{if}~(z_4>0~\&\&~\max\{h(0),h(1),h(\lfloor z_4\rfloor_e+2)\}\le r^*)~||~(z_4\le0~\&\&~h(0)\le r^*)
  \end{cases}
  \end{eqnarray}
  where~$z_4=\frac{ab-ab^3-2cb^2}{c(b^2-1)}$ and~$\lfloor n\rfloor_e$ denotes the maximum even integer not exceeding~$n$;
  \item[2.4] $b<0,c>0$:
  \begin{eqnarray}
  L(\omega,\omega^*) =
  \begin{cases}
  \min\{k:~1\le k\le\lfloor z_4\rfloor_o+2 ,~h(k)>r^*\},\quad\text{if}~z_4>1~\&\&~h(0)\le r^*< h(\lfloor z_4\rfloor_o+2)\\
  \infty, \quad\text{if}~(z_4\le1~\&\&~\max\{h(0),h(1)\}\le r^*)~||~(z_4>1~\&\&~\{h(0),h(1),h(\lfloor z_4\rfloor_o+2)\}\le r^*)
  \end{cases}
  \end{eqnarray}
  where~$z_4=\frac{ab-ab^3-2cb^2}{c(b^2-1)}$ and~$\lfloor n\rfloor_o$ denotes the maximum odd integer not exceeding~$n$;
  \item[2.5] $bc=0$: $h(k)$ is reduced to forms similar to those in (1.1)-(1.6) so details are omitted .
  \end{itemize}
  \item[3.] $\textbf{P}$ has a pair of conjugate complex eigenvalues: $h(k) = a'A^k\sin(k\theta + b') + c'$.
  \begin{itemize}
  \item[3.1] $d'=\frac{r^*-c'}{a'}>0$:
  \begin{eqnarray}
  L(\omega,\omega^*) =
  \begin{cases}
  \min\{k:~0\le k<\lceil\log_A^{d'}\rceil ,~h(k)>r^*\},\quad\text{if}~\log_A^{d'}>0~\&\&~\max\{h(k):~0\le k<\lceil\log_A^{d'}\rceil\}>r^*\\
  \infty, \quad\text{if}~\log_A^{d'}\le0~||~(\log_A^{d'}>0~\&\&~\max\{h(k):~0\le k<\lceil\log_A^{d'}\rceil\}\le r^*)
  \end{cases}
  \end{eqnarray}
  \item[3.2] $d'=\frac{r^*-c'}{a'}<0$: $L(\omega,\omega^*) = \min\{k:~k\ge0,~h(k)>r^*\}$;
  \item[3.3] $d'=\frac{r^*-c'}{a'}=0, b'\in(0,\pi)$: $L(\omega,\omega^*)=0$;
  \item[3.4] $d'=\frac{r^*-c'}{a'}=0, \theta=\pi, b'\in\{0,\pi\}$: $L(\omega,\omega^*)=\infty$;
  \item[3.5] $d'=\frac{r^*-c'}{a'}=0, \theta\neq\pi, b'=0$: $L(\omega,\omega^*)=\lfloor\frac{\pi}{2\pi-\theta}\rfloor+1$;
  \item[3.6] $d'=\frac{r^*-c'}{a'}=0, \theta\neq\pi, b'=\pi$: $L(\omega,\omega^*)=\lfloor\frac{\pi}{\theta}\rfloor+1$;
  \item[3.7] $d'=\frac{r^*-c'}{a'}=0, \theta\in(0,\pi], b'\in(\pi,2\pi)$: $L(\omega,\omega^*)=\lfloor\frac{2\pi-b}{\theta}\rfloor+1$;
  \item[3.8] $d'=\frac{r^*-c'}{a'}=0, \theta\in(\pi,2\pi), b'\in(\pi,2\pi)$: $L(\omega,\omega^*)=\lfloor\frac{b-\pi}{2\pi-\theta}\rfloor+1$.
  \end{itemize}
\end{itemize}
\end{theorem}
\proof{Proof.} The base case~\eqref{eq:baseCase} is clear. We prove the rest case by case in the same order as appeared in the theorem.
\begin{itemize}
  \item[1.] $\textbf{P}$ has only real eigenvalues and~$3$ linearly independent eigenvectors: $h(k) = a_1b_1^k + a_2b_2^k + c$.
  \begin{itemize}
  \item[1.1] $b_1=b_2\neq0~\&\&~b_1>0~\&~a_1+a_2<0$: $h(k)=(a_1+a_2)b_1^k+c$ is monotonically increasing over~$k\ge0$ and the result follows.
  \item[1.2] $b_1=b_2\neq0~\&\&~(b_1<0~||~a_1+a_2\ge0)$: $L(\omega,\omega^*)$ achieves the maximum value at either~$h(0)$ or~$h(1)$ and the result follows.
  \item[1.3] $a_1b_1=0~\&\&~b_2>0~\&\&~a_2<0$:  $h(k)=a_2b_2^k+c$ which is monotonically increasing over~$k\ge1$ and the result follows.
  \item[1.4] $a_1b_1=0~\&\&~(b_2\le0~||~a_2\ge0)$: $L(\omega,\omega^*)$ achieves the maximum value at one of~$\{h(0),h(1),h(2)\}$ and the result follows.
  \item[1.5] $a_2b_2=0~\&\&~b_1>0~\&~a_1<0$: similar to (1.3).
  \item[1.6] $a_2b_2=0~\&\&~(b_1\le0~||~a_1\ge0)$: similar to (1.4).
  \item[1.7] $a_1,a_2,b_1,b_2>0$: $h(k)$ achieves the maximum value at~$h(0)$ and the result follows.
  \item[1.8] $a_1<0,a_2>0,b_1>b_2>0$: observe that
  $$h(k+1) - h(k) > 0 \Leftrightarrow z(k)\defeq\left(\frac{b_1}{b_2}\right)^k > -\frac{a_2(b_2 - 1)}{a_1(b_1 - 1)} (> 0).$$ If there exists a~$k_1\ge0$ satisfying the above, then~$h(k)$ is monotonically decreasing until~$k_1$ after which it increases. So the supremum of~$h(k)$ is achieved at either~$0$ or~$\infty$. If such~$k_1$ does not exist, $h(k)$ is monotonically increasing for all~$k\ge0$ and achieves its supremum at~$\infty$. The result thus follows.
  \item[1.9] $a_1<0,a_2>0,b_2>b_1>0$: contrary to (1.8), if there exists a~$k_1\ge0$ such that~$h(k_1+1) - h(k_1) > 0$, then~$h(k)$ is monotonically increasing until~$k_1+1$ after which it decreases to the stationary reward~$c$ (see Fig.~\ref{fig: Lemma 3-2-1}). So the maximum of~$h(k)$ is achieved at either~$0$ or~$k_1+1$. If such~$k_1$ does not exist, $h(k)$ is monotonically decreasing for all~$k\ge0$ and achieves its maximum at~$0$. The result thus follows.
    \begin{figure}
    	\centering
    	\includegraphics[scale=0.6]{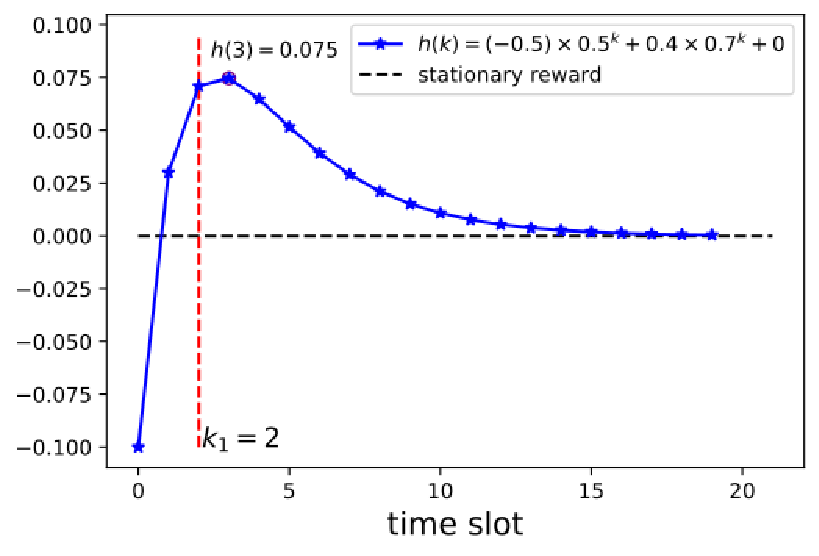}
    	\includegraphics[scale=0.6]{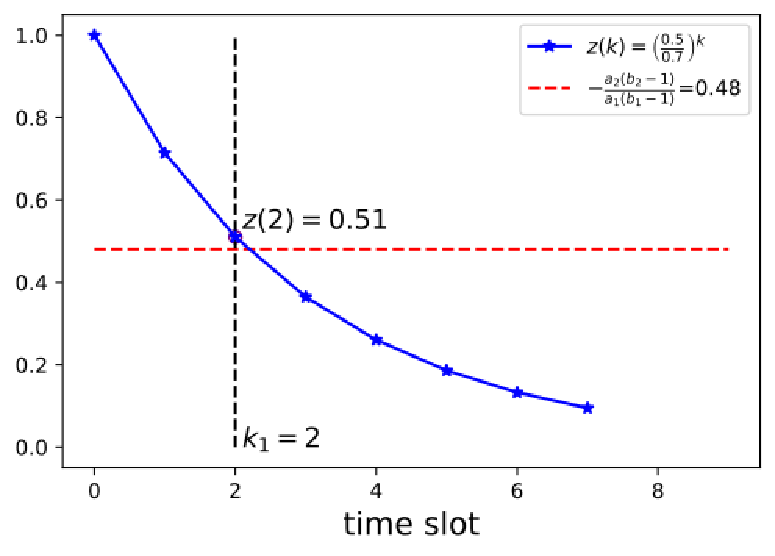}
    	\caption{$h(k) = -0.5\times 0.5^k + 0.4\times 0.7^k$, $z(k) = \left(\frac{0.5}{0.7}\right)^k$}
    	\label{fig: Lemma 3-2-1}
    \end{figure}
  \item[1.10] $b_1<0,a_1,a_2,b_2>0$: since
  $$h(k+1) < h(k),\ h(k+2) < h(k),\quad \forall \text{even number } k\ge0,$$ so $h(k)$ achieves its maximum at~$0$ and the result follows.
  \item[1.11] $a_1,b_1<0,a_2,b_2>0$: observe that
    $$h(k+1) - h(k) > 0 \Leftrightarrow z(k)=\left(\frac{b_1}{b_2}\right)^k > -\frac{a_2(b_2 - 1)}{a_1(b_1 - 1)} (> 0)$$
    $$h(k+2) - h(k) > 0 \Leftrightarrow z(k)=\left(\frac{b_1}{b_2}\right)^k > -\frac{a_2(b_2^2 - 1)}{a_1(b_1^2 - 1)} (> 0)$$
    which directly lead to the following properties:
    $$f(k+1) < f(k),\ f(k+2) < f(k),\quad \forall \text{odd number } k\ge1.$$
    So $h(k)$ achieves its maximum at~$0$ or~$1$ and the result follows. See Fig.~\ref{fig: Lemma 3-4-1} for an example.
    \begin{figure}	
	\centering
	\includegraphics[scale=0.7]{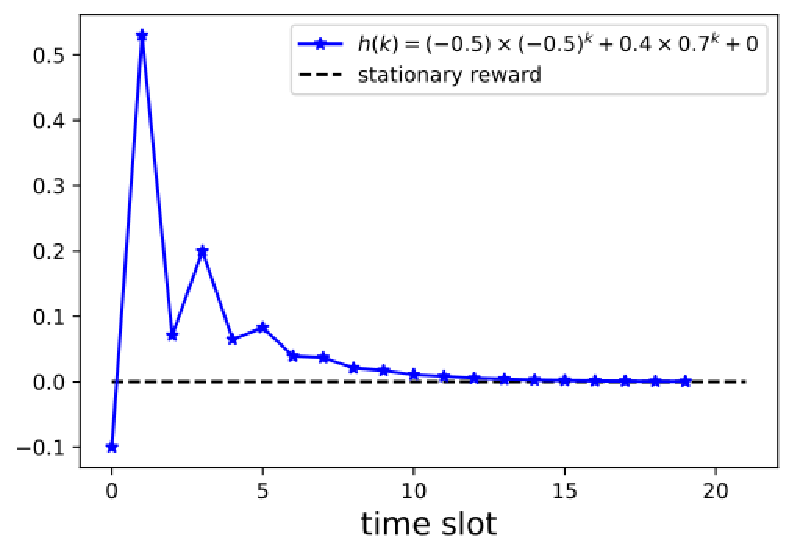}
	\caption{$h(k)=-0.5\times (-0.5)^k + 0.4\times 0.7^k$}
	\label{fig: Lemma 3-4-1}
    \end{figure}
    \item[1.12] $a_2,b_1<0,a_1,b_2>0,|b_1|>b_2$: observe that
    $$h(k+1) - h(k) > 0 \Leftrightarrow z(k)=\left(\frac{b_1}{b_2}\right)^k < -\frac{a_2(b_2 - 1)}{a_1(b_1 - 1)} (> 0),$$
    $$h(k+2) - h(k) > 0 \Leftrightarrow z(k)=\left(\frac{b_1}{b_2}\right)^k < -\frac{a_2(b_2^2 - 1)}{a_1(b_1^2 - 1)} (> 0).$$
    Let~$k_1$ and~$k_2$ be the maximum even integers achieving the above inequalities, respectively. Note that~$k_2\ge k_1$. If both of them are nonnegative, then~$h(k)$ is monotonically increasing until $k_1+2$, then moving up with oscillations until~$k_2+2$ and finally moving downward to converge to~$c$. If~$k_1<0\le k_2$, then~$h(k)$ still achieves its maximum~$k_2+2$. Finally, if~$k_2<0$, $h(k)$ has its maximum at~$0$. The result thus follows. See Figs.~\ref{fig: Lemma 3-5.1-1} and~\ref{fig: Lemma 3-5.1-2} for an example.
    \begin{figure}[hb]
	\centering
	\includegraphics[scale=0.6]{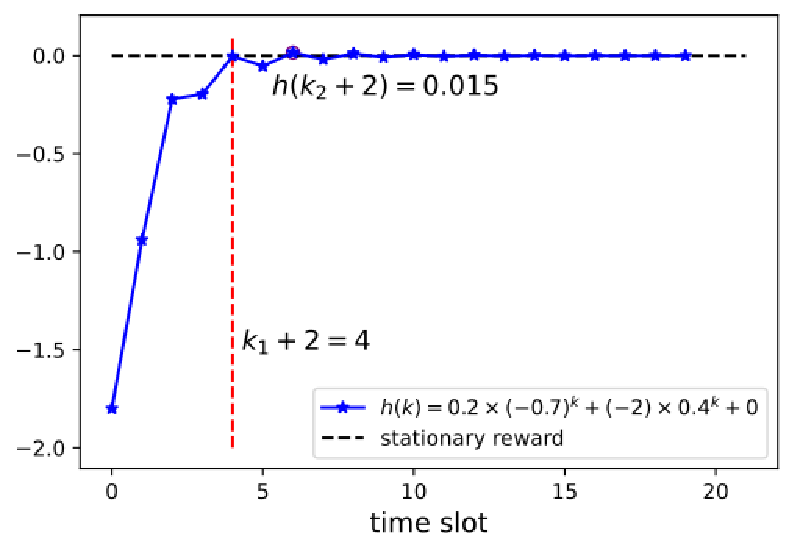}
	\includegraphics[scale=0.6]{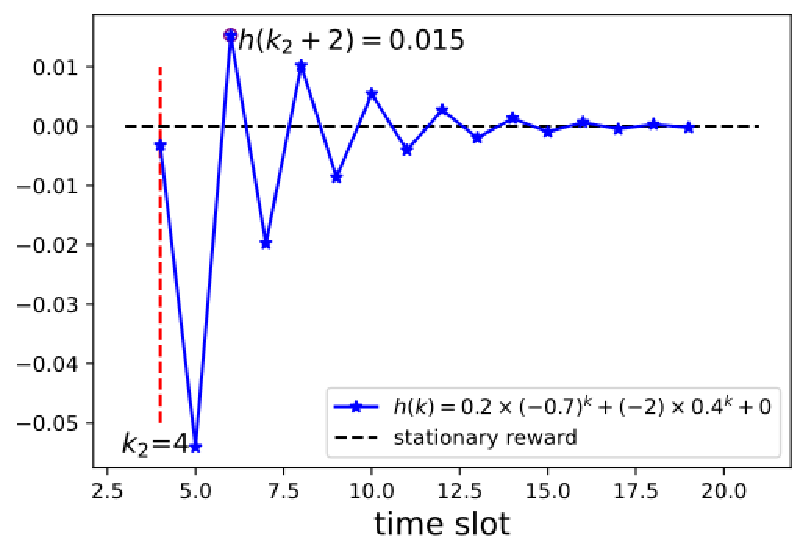}
	\caption{$h(k)=0.2\times (-0.7)^k -2\times 0.4^k$}
    \label{fig: Lemma 3-5.1-1}
    \end{figure}
    \begin{figure}[h]
    	\centering
    	\includegraphics[scale=0.7]{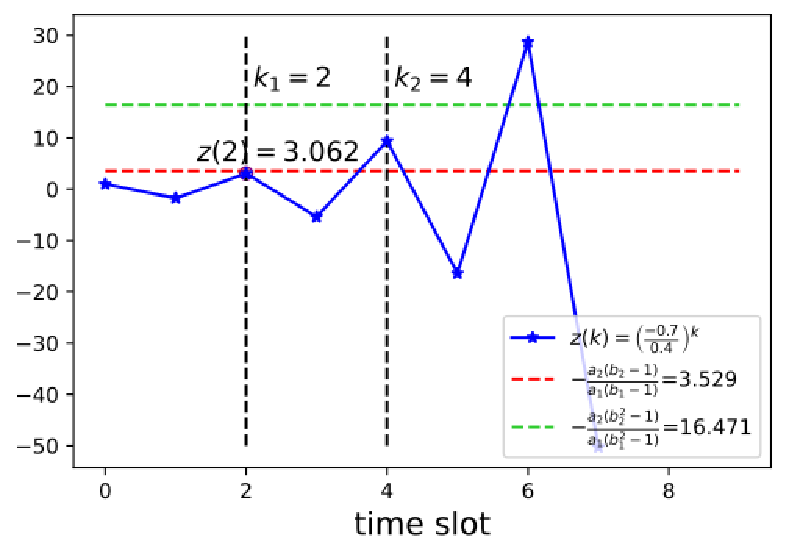}
    	\caption{\ $z(k) = \left(-\frac{0,7}{0.4}\right)^k$}
    	\label{fig: Lemma 3-5.1-2}
    \end{figure}
  \item[1.13] $a_2,b_1<0,a_1,b_2>0,|b_1|<b_2$: this case is sort of the reversed version to (1.12). Let~$k_1\ge0$ and~$k_2\ge0$ be the minimum even integers achieving the two inequalities in (1.12), respectively. In this case, $h(k)$ moves down with oscillations until~$k_2$, then it moves up with oscillations until~$k_1$ and finally increases to the stationary reward~$c$. Therefore~$h(k)$ achieves its supremum at~$0$ or~$\infty$. The result thus follows.
  \item[1.14] $a_2,b_1<0,a_1,b_2>0,|b_1|=b_2$: under this case, the following holds
    $$h(k+1) - h(k) > 0 \Leftrightarrow (-1)^k < -\frac{a_2(b_2 - 1)}{a_1(b_1 - 1)} (> 0),$$
    $$h(k+2) - h(k) > 0 \Leftrightarrow (-1)^k < -\frac{a_2(b_2^2 - 1)}{a_1(b_1^2 - 1)} (> 0).$$
    If $-\frac{a_2(b_2 - 1)}{a_1(b_1 - 1)} \ge 1$, then~$h(k)$ is monotonically increasing to the stationary reward~$c$. If $-\frac{a_2(b_2 - 1)}{a_1(b_1 - 1)} < 1$ and $ -\frac{a_2(b_2^2 - 1)}{a_1(b_1^2 - 1)} \le 1$, then~$h(k)$ oscillates but its maximum value cannot exceed~$h(0)$. If $-\frac{a_2(b_2 - 1)}{a_1(b_1 - 1)} < 1$ and $-\frac{a_2(b_2^2 - 1)}{a_1(b_1^2 - 1)} > 1$, then $h(k)$ moves up with oscillations to its supremum~$c$. The result thus follows.
  \item[1.15] $b_1,b_2<0,a_1,a_2>0$: the maximum of~$h(k)$ clearly happens at~$0$ and the result thus follows.
  \item[1.16] $b_1,b_2>0,a_1,a_2<0$: $h(k)$ monotonically converges to~$c$ from below and the result thus follows.
  \item[1.17] $a_1,a_2,b_1<0,b_2>0,|b_1|>b_2$: under this case, the following holds
    $$h(k+1) - h(k) > 0 \Leftrightarrow z(k)=\left(\frac{b_1}{b_2}\right)^k > -\frac{a_2(b_2 - 1)}{a_1(b_1 - 1)} (< 0),$$
    $$h(k+2) - h(k) > 0 \Leftrightarrow z(k)=\left(\frac{b_1}{b_2}\right)^k > -\frac{a_2(b_2^2 - 1)}{a_1(b_1^2 - 1)} (< 0).$$
    Any even~$k$ clear satisfies the above two inequalities. Let~$k_1\ge1$ and~$k_2\ge1$ be the maximum odd integers achieving the above, respectively. If $k_1, k_2$ exist, then $k_1\le k_2$ and~$h(k)$ monotonically increases until $k_1 + 2$ after which it goes up with oscillations until $k_2 + 2$, and finally it falls with oscillations and converges to $c$. As long as $k_2$ exists, $h(k)$ has its maximum at $k_2 + 2$. When~$k_2$ does not exist, it is clear that~$h(k)$ achieves its maximum at~$1$ and the result follows.
    \item[1.18] $a_1,a_2,b_1<0,b_2>0,|b_1|<b_2$: let~$k_1\ge1$ and~$k_2\ge1$ be the minimum odd integers achieving the two inequalities in (1.17), respectively. Note that~$k_2\le k_1$. Then~$h(k)$ moves down with oscillations until~$k_2$ after which it goes up with oscillations until~$k_1$ and finally~$h(k)$ monotonically increases to~$c$. The result thus follows.
    \item[1.19] $a_1,a_2,b_1<0,b_2>0,|b_1|=b_2$: similar to (1.14), if $-\frac{a_2(b_2 - 1)}{a_1(b_1 - 1)} \le -1$, then~$h(k)$ is monotonically increasing to the stationary reward~$c$. If $-\frac{a_2(b_2 - 1)}{a_1(b_1 - 1)} > -1$ and $ -\frac{a_2(b_2^2 - 1)}{a_1(b_1^2 - 1)} < -1$, then~$h(k)$ moves up with oscillations to~$c$. If $-\frac{a_2(b_2 - 1)}{a_1(b_1 - 1)} > -1$ and $-\frac{a_2(b_2^2 - 1)}{a_1(b_1^2 - 1)} \ge -1$, then~$h(k)$ achieves its maximum at~$1$. The result thus follows.
    \item[1.20] $a_2,b_1,b_2<0,a_1>0,|b_1|>|b_2|$: let
    \begin{eqnarray}
    k_{o1} &=& \min\{k: h(k+1)>h(k), k\text{ is positive and odd}\} - 2,\nn\\
    k_{o2} &=& \min\{k: h(k+2)>h(k), k\text{ is positive and odd}\} - 2,\nn\\
    k_{e1} &=& \max\{k: h(k+1)>h(k),k\text{ is nonnegative and even}\},\nn\\
    k_{e2} &=& \max\{k: h(k+2)>h(k), k\text{ is nonnegative and even}\}.\nn
    \end{eqnarray}
    If~$k_{e1}$ exists, then $|k_{o1} - k_{e1}|=1$ and $|k_{o2}-k_{e2}|=1$. Furthermore, we have that $k_{o1} \le k_{o2}$ and $k_{e1} \le k_{e2}$ (see Fig.~\ref{fig: Lemma 3-9.1-2} for an example). Let $k_1 = \max\{k_{o1}, k_{e1}\}$. From the origin~$0$ to $k_1+1$, we have that $\max_{0\le k\le k_1+1}h(k) = h(1)$. Then from $k_1 + 1$ to $k_{e2} + 2$, it reaches a local maximum $\max_{k>k_1}h(k) = h(k_{e2} + 2)$ after which it moves down to the stationary reward~$c$ (see Fig.~\ref{fig: Lemma 3-9.1-1} for an example). If~$k_{e1}$ does not exist but~$k_{e2}$ does, $h(k)$ attains its maximum value at either~$0$ or~$h(k_{e2} + 2)$. If~$k_{e2}$ does not exist, then~$h(k)$ attains its maximum value at~$0$. The result thus follows.
    \begin{figure}
    	\centering
    	\includegraphics[scale=0.7]{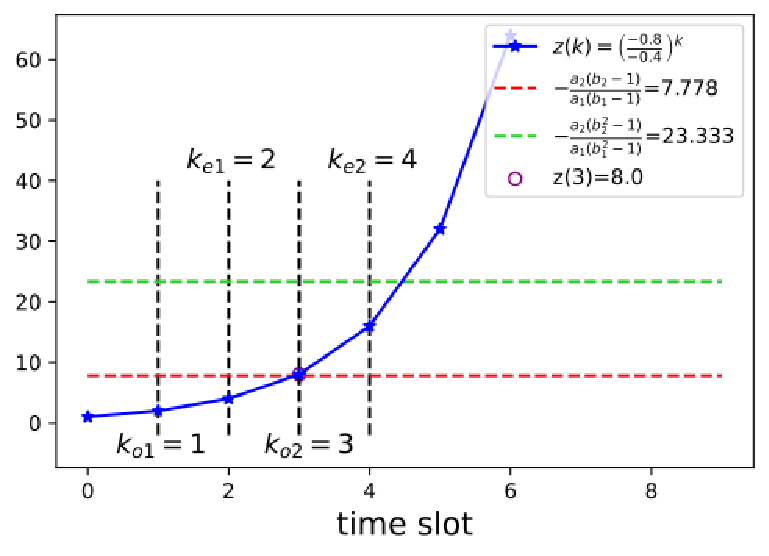}
    	\caption{$z(k) = \left(\frac{-0.8}{-0.4}\right)^k$}
    	\label{fig: Lemma 3-9.1-2}
    \end{figure}
    \begin{figure}
    	\centering
    	\includegraphics[scale=0.6]{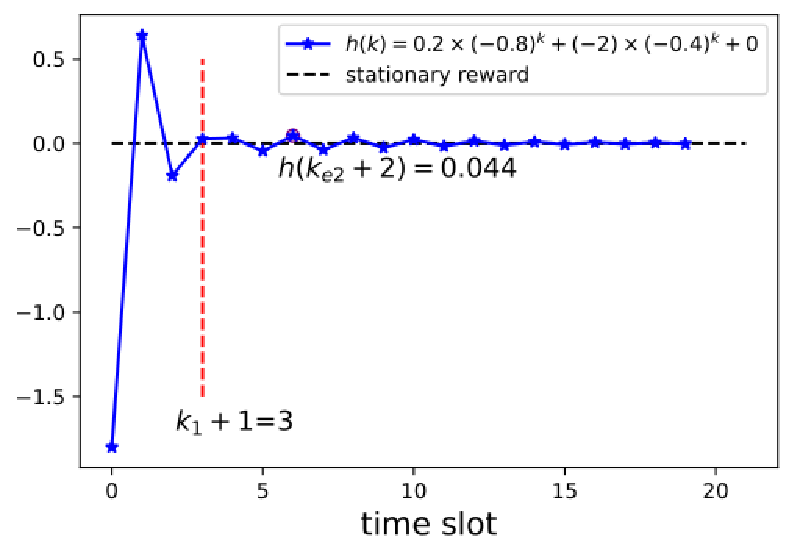}
    	\includegraphics[scale=0.6]{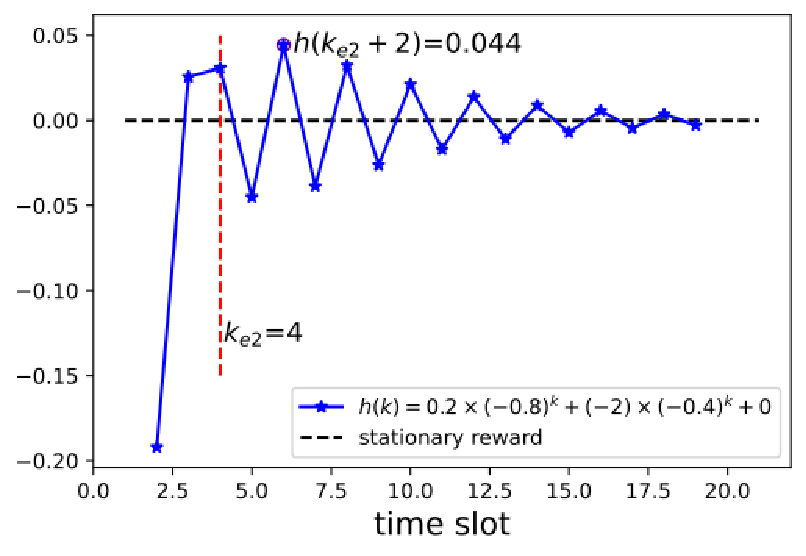}
    	\caption{$h(k) = 0.2\times (-0.8)^k - 2\times (-0.4)^k$}
    	\label{fig: Lemma 3-9.1-1}
    \end{figure}
    \item[1.21] $a_2,b_1,b_2<0,a_1>0,|b_1|<|b_2|$: let
    \begin{eqnarray}
    k_{o1} &=& \max\{k: h(k+1)>h(k), k\text{ is positive and odd}\},\nn\\
    k_{o2} &=& \max\{k: h(k+2)>h(k), k\text{ is positive and odd}\},\nn\\
    k_{e1} &=& \min\{k: h(k+1)>h(k),k\text{ is nonnegative and even}\}-2,\nn\\
    k_{e2} &=& \min\{k: h(k+2)>h(k), k\text{ is nonnegative and even}\}-2.\nn
    \end{eqnarray}
    Similar to (1.20), if~$k_{o1}$ exists, then $|k_{o1} - k_{e1}|=1$, $|k_{o2} - k_{e2}|=1$, $k_{o1} \le k_{o2}$, and $k_{e1} \le k_{e2}$.
    Let $k_1 = \max\{k_{o1}, k_{e1}\}$. We have that $h(0) = \max_{0\le k\le k_1+1}h(k)$ and $\max_{k\ge k_1+2}h(k) = h(k_{o2})$. If~$k_{o1}$ does not exist but~$k_{o2}$ does, then $\arg\max_kh(k)$ is one of $\{0,1,k_{o2}+2\}$ (see Fig.~\ref{fig: Lemma 3-9.2-1} for an example). If~$k_{o2}$ does not exist, then $\arg\max_kh(k)$ is either~$0$ or~$1$. The result thus follows.
    \begin{figure}
    	\centering
    	\includegraphics[scale=0.6]{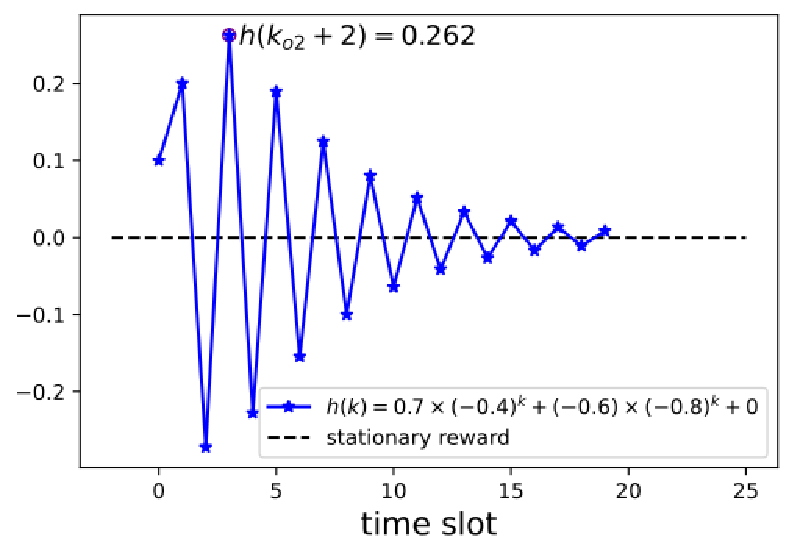}
    	\includegraphics[scale=0.6]{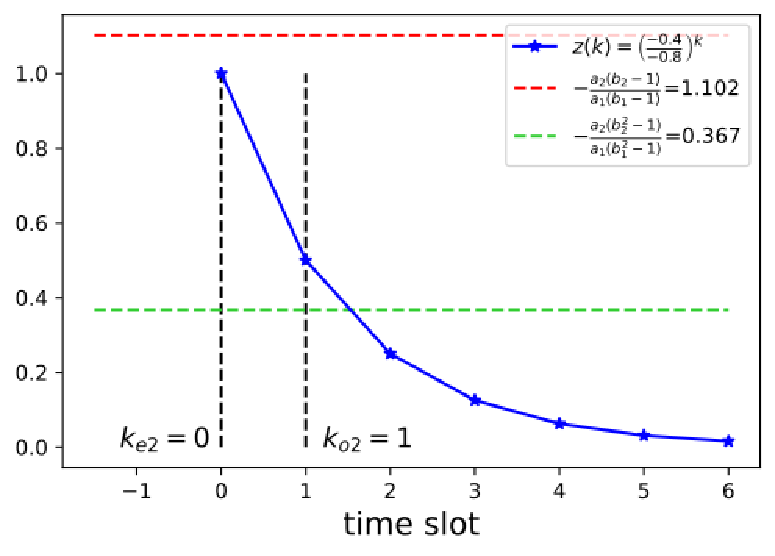}
    	\caption{$h(k)=0.7\times (-0.4)^k - 0.6\times (-0.8)^k,\ g(k) = \left(\frac{-0.4}{-0.8}\right)^k$}
    	\label{fig: Lemma 3-9.2-1}
    \end{figure}
    \item[1.22] $a_2,b_1,b_2<0,a_1<0$: obviously~$h(k)$ achieves its maximum at~$1$ and the result follows.
  \end{itemize}
  \item[2.] $\textbf{P}$ has only real eigenvalues and~$2$ linearly independent eigenvectors: $h(k) = ab^k + ckb^{k-1} + d$.
  \begin{itemize}
  \item[2.1] $b,c>0$: observe that
    $$h(k+1) > h(k) \Leftrightarrow k < \frac{ab - ab^2 - cb}{c(b - 1)}.$$
    Let~$k_1\ge0$ be the maximum integer satisfying the above inequality. If it exists, then~$h(k)$ will keep increasing until $(k_1+1)$ after which it turns to be monotonically decreasing to the stationary reward~$d$. Hence, $k_1 + 1 = \arg\max_k h(k)$ (see Fig.~\ref{fig: Lemma 4-1-1} for an example). If $k_1$ does not exist, then~$h(k)$ is monotonically decreasing and $\arg\max_kh(k) = 0$. The result thus follows.
    \begin{figure}
    	\centering
    	\includegraphics[scale=0.7]{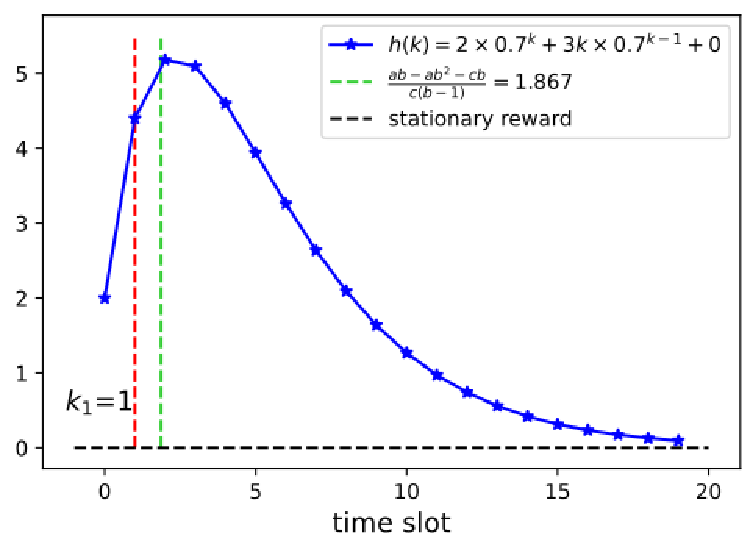}
    	\caption{$h(k) = 2\times 0.7^k +3k\times 0.7^{k-1}$}
    	\label{fig: Lemma 4-1-1}
    \end{figure}
  \item[2.2] $b>0,c<0$: observe that
  $$h(k+1) > h(k) \Leftrightarrow k > \frac{ab - ab^2 - cb}{c(b - 1)}.$$
    Let $k_1\ge0$ be the minimum integer satisfying the above inequality. Clearly~$h(k)$ is monotonically decreasing until $k_1$ after which it keeps increasing to~$d$. So~$h(k)$ achieves its supremum at either~$0$ or~$\infty$ and the result follows.
  \item[2.3] $b<0,c<0$: the proof is similar to that of (1.20) and omitted here (see Fig.~\ref{fig: Lemma 4-4-1} for an example).
  \begin{figure}
	\centering
	\includegraphics[scale=0.6]{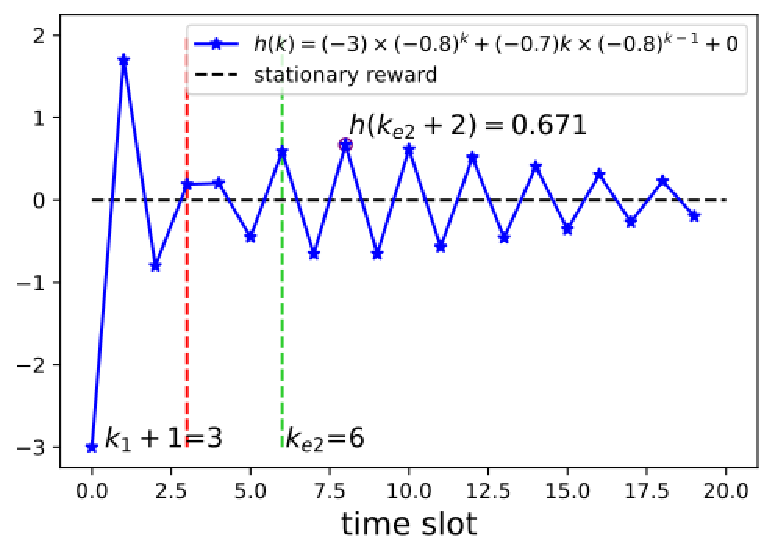}
	\includegraphics[scale=0.6]{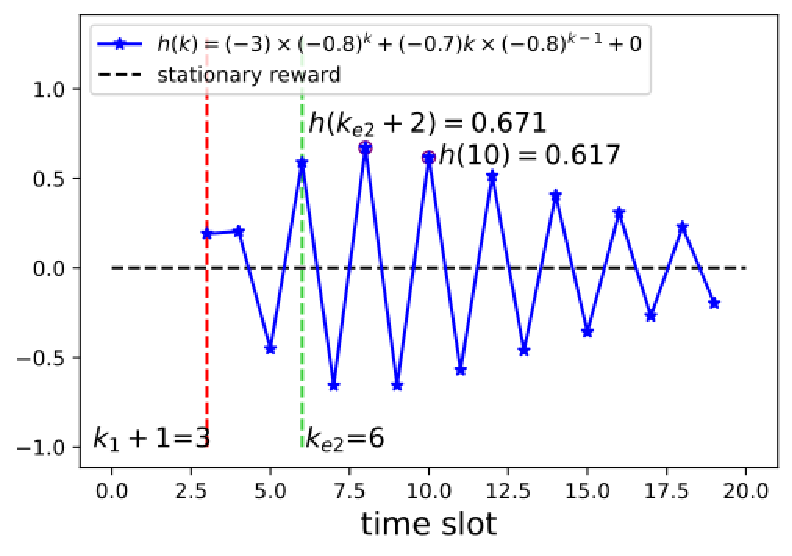}
    \caption{$h(k)=-3\times (-0.8)^k -0.7k\times (-0.8)^{k-1}$}
	\label{fig: Lemma 4-4-1}
  \end{figure}
  \item[2.4] $b<0,c>0$: the proof is similar to that of (1.21) and omitted here.
  \end{itemize}
  \item[3.] $\textbf{P}$ has a pair of conjugate complex eigenvalues: $h(k) = a'A^k\sin(k\theta + b') + c'$.
  \begin{itemize}
  \item[3.1] $d'=\frac{r^*-c'}{a'}>0$: clearly~$h(k)$ will be smaller than~$d'$ as~$k$ becomes sufficiently large and the result follows.
  \item[3.2] $d'=\frac{r^*-c'}{a'}<0$: clearly~$h(k)$ will be larger than~$d'$ as~$k$ becomes sufficiently large and the exhaustion stops in finite time.
  \item[3.3-3.8] These cases follow directly by finding the first~$k\ge0$ such that $\sin(k\theta+b')>0$ and we omit the details here.
  \end{itemize}
\end{itemize}
\Halmos
\endproof

\end{document}